\documentclass{etds}     

\usepackage{chemarr}
\usepackage{amsfonts}
\usepackage{amssymb}
\usepackage{latexsym}
\usepackage{graphicx}
\usepackage{amsbsy}
\usepackage{psfrag}

\newcommand{\D}{\displaystyle}
\newcommand{\restr}[1]{\raisebox{-0.3em}{$\lb|_{#1}\rb.$}} 
\newcommand{\ignore}[1]{}    
\newcommand{\breath}{\medskip} 
\newtheorem{thm}{Theorem}[section]

\newcounter{claimcount}[thm]  
\newtheorem{prop}[thm]{Proposition} 

\newtheorem{lemma}[thm]{Lemma} 

\newtheorem{cor}[thm]{Corollary}

\newcommand{\dfn}{\sf\em} 
\newcommand{\Theorem}[2]{\begin{thm}{\sf #1}  #2 \end{thm}}
\newcommand{\Proposition}[2]{\begin{prop}{\sf #1}  #2 \end{prop}}
\newcommand{\Lemma}[2]{\begin{lemma}{\sf #1}  #2 \end{lemma}}
\newcommand{\Corollary}[2]{\begin{cor}{\sf #1}  #2 \end{cor}} 
\newcommand{\thmfont}[1]{{\sl #1}}    
\newcommand{\example}[1]{        \refstepcounter{thm}                     \begin{list}{} 			{\setlength{\leftmargin}{0em} 			\setlength{\rightmargin}{0em}}        \item {\sc Example \thethm:} #1                   \hfill$\diamondsuit$  \end{list}   			}     
\newcommand{\bthmlist}{ \begin{list}{{\bf(\alph{enumi})}} {\usecounter{enumi} \setlength{\leftmargin}{1em} \setlength{\itemsep}{0.2em} \setlength{\topsep}{0.2em} \setlength{\itemindent}{0em} \setlength{\parsep}{0em} \setlength{\rightmargin}{0em}} } 
\newcommand{\ethmlist}{\end{list}}    
\newcommand{\Claim}[1]{\refstepcounter{claimcount}                \noindent {\sc Claim \theclaimcount: \ }\thmfont{ #1}} 
\newcommand{\bprf}[1][Proof:]{\begin{list}{} 			{\setlength{\leftmargin}{0.7em} 			\setlength{\rightmargin}{0em} 			\setlength{\listparindent}{1em}}                         \item {\em \hspace{-1em}  #1  }} 
\newcommand{\eprf}{\end{list}} 
\newcommand{\bthmprf}{\bprf}
\newcommand{\bclaimprf}{\bprf}
\newcommand{\ethmprf}{ \hfill$\Box$  \eprf  \breath  } 
\newcommand{\eclaimprf}{ \hfill $\Diamond$~{\scriptsize {\tt Claim~\theclaimcount}}\eprf}  
\newcommand{\QED}{\hfill\ensuremath{\Box}}
\newcommand{\qed}{\QED}     
\newcommand{\beq}{\begin{eqnarray*}}
\newcommand{\eeq}{\end{eqnarray*}} 
\newcommand{\beqn}{ \begin{equation} }
\newcommand{\eeqn}{ \end{equation} }
\newcommand{\blist}{\begin{enumerate}}
\newcommand{\elist}{\end{enumerate}} 
\newcommand{\bitem}{\begin{itemize}}
\newcommand{\eitem}{\end{itemize}} 
\newcommand{\bdesc}{\begin{description}}
\newcommand{\edesc}{\end{description}}   
\newcommand{\dB}{{\mathbb{B}}}
\newcommand{\dD}{{\mathbb{D}}}
\newcommand{\dE}{{\mathbb{E}}}
\newcommand{\dF}{{\mathbb{F}}}
\newcommand{\dG}{{\mathbb{G}}}
\newcommand{\dH}{{\mathbb{H}}}
\newcommand{\dK}{{\mathbb{K}}}
\newcommand{\dN}{{\mathbb{N}}}
\newcommand{\dR}{{\mathbb{R}}}
\newcommand{\dS}{{\mathbb{S}}}
\newcommand{\dU}{{\mathbb{U}}}
\newcommand{\dV}{{\mathbb{V}}}
\newcommand{\dW}{{\mathbb{W}}}
\newcommand{\dX}{{\mathbb{X}}}
\newcommand{\dY}{{\mathbb{Y}}}
\newcommand{\dZ}{{\mathbb{Z}}}       
\newcommand{\barA}{{\overline{A}}}
\newcommand{\barV}{{\overline{V}}}
\newcommand{\barpi}{{\overline{\pi }}}
\newcommand{\barsH}{{\overline{\mathcal{ H}}}}
\newcommand{\bC}{{\mathbf{ C}}}
\newcommand{\bG}{{\mathbf{ G}}}
\newcommand{\bK}{{\mathbf{ K}}}
\newcommand{\bS}{{\mathbf{ S}}}
\newcommand{\bT}{{\mathbf{ T}}}
\newcommand{\bU}{{\mathbf{ U}}}
\newcommand{\bV}{{\mathbf{ V}}}
\newcommand{\bW}{{\mathbf{ W}}}
\newcommand{\bX}{{\mathbf{ X}}}
\newcommand{\bY}{{\mathbf{ Y}}}
\newcommand{\ba}{{\mathbf{ a}}}
\newcommand{\bb}{{\mathbf{ b}}}
\newcommand{\bc}{{\mathbf{ c}}}
\newcommand{\bd}{{\mathbf{ d}}}
\newcommand{\bi}{{\mathbf{ i}}}
\newcommand{\bp}{{\mathbf{ p}}}
\newcommand{\bq}{{\mathbf{ q}}}
\newcommand{\bs}{{\mathbf{ s}}}
\newcommand{\bw}{{\mathbf{ w}}}
\newcommand{\bx}{{\mathbf{ x}}}
\newcommand{\bzet }{{\boldsymbol{\zeta}}}
\newcommand{\bphi}{{\boldsymbol{\phi }}}
\newcommand{\sA}{{\mathcal{ A}}}
\newcommand{\sB}{{\mathcal{ B}}}
\newcommand{\sC}{{\mathcal{ C}}}
\newcommand{\sD}{{\mathcal{ D}}}
\newcommand{\sF}{{\mathcal{ F}}}
\newcommand{\sG}{{\mathcal{ G}}}
\newcommand{\sH}{{\mathcal{ H}}}
\newcommand{\sI}{{\mathcal{ I}}}
\newcommand{\sN}{{\mathcal{ N}}}
\newcommand{\sP}{{\mathcal{ P}}}
\newcommand{\sQ}{{\mathcal{ Q}}}
\newcommand{\sT}{{\mathcal{ T}}}
\newcommand{\sV}{{\mathcal{ V}}}
\newcommand{\sW}{{\mathcal{ W}}}
\newcommand{\sX}{{\mathcal{ X}}}
\newcommand{\sZ}{{\mathcal{ Z}}}     
\newcommand{\gA}{{\mathfrak{ A}}}
\newcommand{\gB}{{\mathfrak{ B}}}
\newcommand{\gC}{{\mathfrak{ C}}}
\newcommand{\gD}{{\mathfrak{ D}}}
\newcommand{\gM}{{\mathfrak{ M}}}
\newcommand{\gQ}{{\mathfrak{ Q}}}
\newcommand{\gT}{{\mathfrak{ T}}}
\newcommand{\gW}{{\mathfrak{ W}}}
\newcommand{\gX}{{\mathfrak{ X}}}
\newcommand{\gh}{{\mathfrak{ h}}}
\newcommand{\go}{{\mathfrak{ o}}}
\newcommand{\alp }{\alpha}
\newcommand{\bet }{\beta}
\newcommand{\gam }{\gamma}
\newcommand{\del }{\delta}
\newcommand{\eps }{\epsilon}
\newcommand{\zet }{\zeta}
\newcommand{\sig }{\sigma} 
\newcommand{\omg }{\omega}
\newcommand{\Gam }{\Gamma}
\newcommand{\Del }{\Delta}
\newcommand{\fd}{{\mathsf{ d}}}
\newcommand{\fe}{{\mathsf{ e}}}
\newcommand{\fs}{{\mathsf{ s}}}
\newcommand{\ft}{{\mathsf{ t}}}
\newcommand{\fu}{{\mathsf{ u}}}
\newcommand{\fv}{{\mathsf{ v}}}
\newcommand{\fw}{{\mathsf{ w}}}
\newcommand{\fx}{{\mathsf{ x}}}
\newcommand{\fy}{{\mathsf{ y}}}
\newcommand{\fz}{{\mathsf{ z}}}    
\newcommand{\tl}{\widetilde} 
\newcommand{\tlB}{{\widetilde{B}}}
\newcommand{\tlC}{{\widetilde{C}}}
\newcommand{\tlbX}{{\widetilde{\mathbf{ X}}}}
\newcommand{\tlsC}{{\widetilde{\mathcal{ C}}}}
\newcommand{\tlsG}{{\widetilde{\mathcal{ G}}}}
\newcommand{\tlsH}{{\widetilde{\mathcal{ H}}}}
\newcommand{\tlzet }{{\widetilde{\zeta}}}
\newcommand{\undC}{{\underline{C}}}

\newcommand{\undalp }{{\underline{\alpha}}}
\newcommand{\undbet }{{\underline{\beta}}}
\newcommand{\undgam }{{\underline{\gamma}}}
\newcommand{\undzet }{{\underline{\zeta}}}
\newcommand{\vV}{{\vec{V}}}
\newcommand{\lb}{\left}
\newcommand{\rb}{\right} 
\newcommand{\Array}[2][cccccccccccccccccccccccccccccccccccc] {{\begin{array}{#1}#2\end{array}}}      
\newcommand{\maketall}{\rule[-0.5em]{0em}{1em}}        
\newcommand{\map}{{\longrightarrow}}
\newcommand{\goto}{{\rightarrow}}
\newcommand{\into}{{\map}}
\newcommand{\inject}{\rightarrowtail}
\newcommand{\surject}{\twoheadrightarrow}
\newcommand{\image}[1]{\mathrm{img}\lb(#1\rb)}  
\newcommand{\seilpmi}{{\Longleftarrow}}
\newcommand{\statement}[1]{\lb(  \maketall       \begin{minipage}{40em}       \begin{tabbing}         #1        \end{tabbing}      \end{minipage}  \rb)}     
\newcommand{\oo}{{\infty}}        
\newcommand{\X}{\times}
\newcommand{\x}{\X}
\newcommand{\dirsum}{\oplus}
\newcommand{\Dirsum}{\bigoplus} 
\newcommand{\union}{\cup}
\newcommand{\Union}{\bigcup}
\newcommand{\intsct}{\cap}
\newcommand{\Intsct}{\bigcap}
\newcommand{\disj}{\sqcup}
\newcommand{\Disj}{\bigsqcup}   
\newcommand{\set}[2]{{\left\{ #1 \; ; \; #2 \right\} }} 
\newcommand{\Id}[1]{{\mathbf{ Id}_{{#1}}}}
\newcommand{\eeequals}[1]{\raisebox{-0.9ex}{$\overline{\overline{{\scriptscriptstyle{\mathrm{#1}}}}}$}} 
\newcommand{\leeeq}[1]{\raisebox{-1ex}{${{\D\leq} \atop {\scriptscriptstyle{\mathrm{#1}}}}$}} 
\newcommand{\iiimplies}[1]{\eeequals{#1}\!\!\!\!\!\!\Rightarrow}
\newcommand{\seilpmiii}[1]{\Leftarrow\!\!\!\!\eeequals{#1}}
\newcommand{\iiiff}[1]{\Leftarrow\!\!\!\!\lefteqn{\eeequals{#1}}\Rightarrow}      
\newcommand{\Fix}[1]{{\sf Fix}\lb[#1\rb]}
\newcommand{\shift}[1]{\sig^{#1}}    
\newcommand{\End}[2][]{\mathsf{End}_{#1} \lb(#2\rb)}
\newcommand{\Hom}[2][]{\mathsf{Hom}_{#1} \lb(#2\rb)}      
\newcommand{\Real}{\dR}
\newcommand{\Natur}{\dN}
\newcommand{\Zahl}{\dZ}
\newcommand{\Zahlmod}[1]{{\Zahl_{/#1}}}
\newcommand{\CC}[1]{{\lb[ #1 \rb]}}
\newcommand{\CO}[1]{{\lb[ #1 \rb)}}
\newcommand{\OC}[1]{{\lb( #1 \rb]}}
\renewcommand{\implies}{\ensuremath{\Longrightarrow}}
\renewcommand{\And}{\mbox{\ and \ }} 
\newcommand{\ZD}[1][D]{{\Zahl^{#1}}}
\newcommand{\RD}[1][D]{{\Real^{#1}}}
\newcommand{\AZD}[1][D]{\sA^{\ZD[#1]}}
\newcommand{\AZ}{\sA^{\Zahl}}
\newcommand{\BZD}[1][D]{\sB^{\ZD[#1]}}
\newcommand{\TZD}[1][D]{\sT^{\ZD[#1]}}
\newcommand{\WZD}[1][D]{\sW^{\ZD[#1]}} 
\newcommand{\tlgD}{\widetilde{\gD}} 
\newcommand{\WangTile}[5][]{\mbox{\tiny$\begin{array}{rcl} {}_\ulcorner\! &\!\!\!\!\!\! \mbox{\scriptsize $#3$}\!\!\!\!\!\! &\! {}_\urcorner \\   \mbox{\scriptsize $#2$}      \! &\!\!\!\! {#1}\!\!\!\! &\! \mbox{\scriptsize $#4$} \\ {}^\llcorner\! &\!\!\!\!\!\! \mbox{\scriptsize $#5$}\!\!\!\!\!\!  &\! {}^\lrcorner  \end{array}$}}       
\newcommand{\CornerTile}[2]{{}^{\ulcorner}_{#1} {#2}{}^{\urcorner}_{\lrcorner}} 
\newsavebox{\sqnbox}
\newsavebox{\sqsbox}
\newsavebox{\sqwbox}
\newsavebox{\sqebox} 
\savebox{\sqnbox}{\includegraphics[scale=0.25]{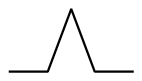}}
\savebox{\sqsbox}{\raisebox{-0.2em}{\includegraphics[scale=0.25]{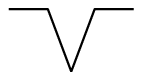}}}
\savebox{\sqebox}{\raisebox{-0.2em}{\includegraphics[scale=0.25]{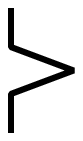}}}
\savebox{\sqwbox}{\raisebox{-0.2em}{\includegraphics[scale=0.25]{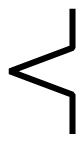}}}
\newcommand{\sqn}{\usebox{\sqnbox}}
\newcommand{\sqs}{\usebox{\sqsbox}}
\newcommand{\sqe}{ \, \usebox{\sqebox}}
\newcommand{\sqw}{\usebox{\sqwbox} \, }               
\newsavebox{\westredpathbox}
\newsavebox{\southredpathbox}
\newsavebox{\westbluepathbox}
\newsavebox{\southbluepathbox} 
\savebox{\westredpathbox}{\,\raisebox{-0.7em}{\includegraphics[scale=0.35]{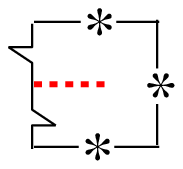}}$\!$}
\savebox{\southredpathbox}{\raisebox{-0.5em}{\includegraphics[scale=0.35]{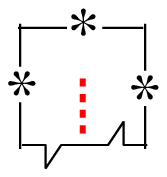}}}
\savebox{\westbluepathbox}{\,\raisebox{-0.7em}{\includegraphics[scale=0.35]{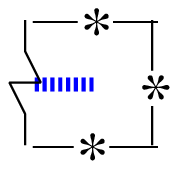}}$\!$}
\savebox{\southbluepathbox}{\raisebox{-0.5em}{\includegraphics[scale=0.35]{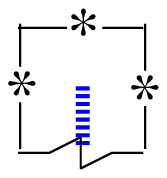}}} 
\newcommand{\westredpath}{\usebox{\westredpathbox}}
\newcommand{\southredpath}{\usebox{\southredpathbox}}
\newcommand{\westbluepath}{\usebox{\westbluepathbox}}
\newcommand{\southbluepath}{\usebox{\southbluepathbox}}  
\newcommand{\Sft}{\gA}
\newcommand{\Remarks}[1]{                      \begin{list}{} 			{\setlength{\leftmargin}{0em} 			\setlength{\rightmargin}{0em}}         \item {\em Remarks:} #1                   \hfill$\diamondsuit$  \end{list}   			} 
\newcommand{\Remark}[1]{                      \begin{list}{} 			{\setlength{\leftmargin}{0em} 			\setlength{\rightmargin}{0em}}         \item {\em Remark:} #1                   \hfill$\diamondsuit$  \end{list}   			} 
\newcommand{\invlim}{\lim_{\longleftarrow}} 
\newcommand{\dirlim}{\lim_{\longrightarrow}}
\newcommand{\tlgA}{\widetilde{\gA}}
\newcommand{\varsig}{\varsigma}
\newcommand{\Nh}{\dH}
\newcommand{\longdownarrow}{\lb.\maketall\rb\downarrow}
\newcommand{\Longdownarrow}{\lb.\maketall\rb\Downarrow}
\newcommand{\longuparrow}{\lb.\maketall\rb\uparrow}
\newcommand{\Kurka}{K\r{u}rka } 
\newcommand{\Mono}{\gM{\scriptstyle\go}}
\newcommand{\Checker}{\gC{\scriptstyle\gh}}
\newcommand{\Ice}{\mathfrak{I{\scriptstyle\!ce}}}
\newcommand{\tlIce}{\widetilde{\Ice}}
\newcommand{\Dom}{\mathfrak{D{\scriptstyle\!om}}}
\newcommand{\tlDom}{\widetilde{\Dom}}
\newcommand{\Pth}{\mathfrak{P{\scriptstyle\!th}}}
\newcommand{\tlPth}{\widetilde{\Pth}}   
\newcommand{\adjacent}[1][]{\stackrel{#1}{\leadsto}}   
\newcommand{\energy}[1][\ba]{\sF_{#1}} 
\newcommand{\unflawed}[1][r]{\dG_{#1}}
\newcommand{\homotopic}[1][\dY]{ \; \raisebox{-1.3ex}{$\widetilde{\widetilde{\,{\scriptscriptstyle{#1}\,}}}$} \; } 
\newcommand{\elhomotopic}[1][\dY]{ \; \raisebox{-1ex}{$\widetilde{\,{\scriptscriptstyle{#1}}\,}$} \; }  
\newcommand{\tilt}[1][\ba]{\rightthreetimes^{C}_{#1}}
\newcommand{\Tilt}[1][\ba]{\rightthreetimes_{#1}} 
\newcommand{\residue}[1][\ba]{\mathrm{Res}_{#1}}
\newcommand{\Residue}[2][r]{\mathrm{Res}^{#1}_{#2} } 
\newcommand{\cgap}[1][\ba]{C_{#1}}     
\newcommand{\Homeo}{\mathrm{H\scriptstyle{omeo}}}    
\newcommand{\onecell}[2]{#1\raisebox{0.1em}{$\scriptscriptstyle\rightarrow$}#2}  
\newcommand{\twocell}[4]{\hspace{-0.2em}\raisebox{-0.5em}{ \psfrag{A}[][]{$\scriptstyle #1$} \psfrag{B}[][]{$\scriptstyle #2$} \psfrag{C}[][]{$\scriptstyle #3$} \psfrag{D}[][]{$\scriptstyle #4$} \includegraphics[scale=0.3]{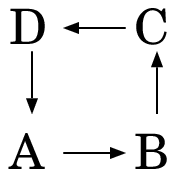}}} 
\newcommand{\threecell}[8]{\hspace{-0.1em}\raisebox{-0.6em}{ \psfrag{A}[][]{$\scriptstyle #1$} \psfrag{B}[][]{$\scriptstyle #2$} \psfrag{C}[][]{$\scriptstyle #3$} \psfrag{D}[][]{$\scriptstyle #4$} \psfrag{E}[][]{$\scriptstyle #5$} \psfrag{F}[][]{$\scriptstyle #6$} \psfrag{G}[][]{$\scriptstyle #7$} \psfrag{H}[][]{$\scriptstyle #8$} \includegraphics[scale=0.3]{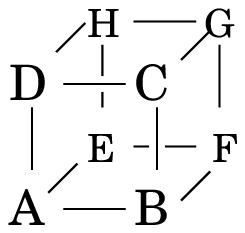}}\, }  
\newcommand{\homoto}{\approx}  
\newcommand{\bkw}[1]{\stackrel{\leftarrow}{#1}\!} 
\newcommand{\Heq}{\sH_{^{\mathrm{eq}}}}
\newcommand{\Ceq}{\sC_{^{\mathrm{eq}}}}
\newcommand{\Zeq}{\sZ_{^{\mathrm{eq}}}}
\newcommand{\Beq}{\sB_{^{\mathrm{eq}}}}
\newcommand{\bCeq}{\bC_{^{\mathrm{eq}}}} 
\newcommand{\rHeq}[1][r]{\,{}_{#1}\!\sH_{^{\mathrm{eq}}}}
\newcommand{\rCeq}[1][r]{\,{}_{#1}\sC_{^{\mathrm{eq}}}}
\newcommand{\rZeq}[1][r]{\,{}_{#1}\!\sZ_{^{\mathrm{eq}}}}
\newcommand{\rBeq}[1][r]{\,{}_{#1}\sB_{^{\mathrm{eq}}}}
\newcommand{\rbCeq}[1][r]{\,{}_{#1}\bC_{^{\mathrm{eq}}}}  
\newcommand{\Hinv}{\sH_{^{\mathrm{inv}}}}
\newcommand{\Cinv}{\sC_{^{\mathrm{inv}}}}
\newcommand{\Zinv}{\sZ_{^{\mathrm{inv}}}}
\newcommand{\Binv}{\sB_{^{\mathrm{inv}}}}
\newcommand{\bCinv}{\bC_{^{\mathrm{inv}}}}  
\newcommand{\rHinv}[1][r]{\,{}_{#1}\!\sH_{^{\mathrm{inv}}}}
\newcommand{\Hdyn}{\sH_{^{\mathrm{dy}}}}
\newcommand{\Zdyn}{\sZ_{^{\mathrm{dy}}}}
\newcommand{\Bdyn}{\sB_{^{\mathrm{dy}}}}  
\newcommand{\Ztwo}{\ddot{\sZ}}   
\newcommand{\blankspace}{\underline{\ \ }} 
\newcommand{\edgematch}[1][d]{\stackrel{#1}{\leadsto}} 
\newcommand{\Adm}[1][r]{\gA_{(#1)}}
\newcommand{\Ext}[1]{\mathrm{Ext}\lb(#1\rb)}
\newcommand{\biota}{\boldsymbol{\iota}} 

\begin{document}

\ETDS{199}{240}{27(\#1)}{2007} 
\title{Algebraic invariants for crystallographic defects in cellular automata}
\runningheads{M. Pivato}{Algebraic invariants for crystallographic defects in cellular automata}
\author{Marcus Pivato}
\address{Dept. of Mathematics \& Computer Science, 
Wesleyan University\footnote{Middletown, CT 06459-0128 USA\qquad
* Peterborough, Ontario, Canada, K9L 1Z6.} \\
and Department of Mathematics, Trent University*}

\email{marcuspivato@trentu.ca}

\recd{10 August 2005 and accepted 29 June 2006}

\begin{abstract}
Let $\AZD$ be the Cantor space of $\ZD$-indexed configurations in a
finite alphabet $\sA$, and let $\shift{}$ be the $\ZD$-action of
shifts on $\AZD$.  A {\dfn cellular automaton} is a continuous,
$\shift{}$-commuting self-map $\Phi$ of $\AZD$, and a {\dfn $\Phi$-invariant
subshift} is a closed, $(\Phi,\shift{})$-invariant subset
$\gA\subset\AZD$.  Suppose $\ba\in\AZD$ is $\gA$-admissible everywhere
except for some small region we call a {\dfn defect}. It has been
empirically observed that such defects persist under iteration of
$\Phi$, and often propagate like `particles' which coalesce or
annihilate on contact.  We construct algebraic invariants
for these defects, which explain their persistence under $\Phi$, and
partly explain the outcomes of their collisions.
Some invariants are based on the cocycles of multidimensional
subshifts; others arise from the higher-dimensional
(co)homology/homotopy groups for subshifts, obtained by generalizing
the Conway-Lagarias tiling groups and the Geller-Propp fundamental
group.

{\footnotesize
\breath

\begin{tabular}{rl}
{\bf MSC:}& 37B50 (primary), 37B15, 37A20 (secondary)\\
{\bf Keywords:}& Cellular automata, subshift, cocycle,
cohomology, tiling group, \\
& defect, kink, domain boundary.
\end{tabular}}
\end{abstract}

  A striking phenomenon in cellular automata is the emergence
of homogeneous `domains' (each exhibiting a
particular spatial pattern), punctuated by {\em defects} (analogous to
`domain boundaries' or `kinks' in a crystalline solid) which evolve
and propagate over time, and occasionally collide.
This phenomenon has been studied empirically in 
\cite{Gra83,Gra84,KrSp88,BoRo91,BNR91,Han,CrHa92,CrHa93a,CrHa93b,CrHa97,CrHM}
and theoretically in 
\cite{Lin84,ElNu,Elo93a,Elo93b,Elo94,CrHR,KuMa00,Kur03,KuMa02,Kur03,Kur05};
see \cite{PivatoDefect0,PivatoDefect1} for a summary.  
The mathematical theory of cellular automaton defect dynamics is still
in its infancy.
Even the term `defect' does 
not yet have a unanimous definition.  Other open questions include:
\blist
  \item   \label{Q:persist} Why do defects persist under the action of cellular
automata, rather than disappearing?  Are there `topological' constraints
imposed by the structure of the underlying domain, which make
defects indestructible?

  \item \label{Q:chem} When defects collide, they often coalesce into
a new type of defect, or mutually annihilate.  Is there a `chemistry'
governing these defect collisions?

  \item \label{Q:inv} Can we assign algebraic invariants to defects,
which reflect {\bf(a)} the `topological constraints' of question
\#\ref{Q:persist} or {\bf(b)} the `defect chemistry' of question
\#\ref{Q:chem}?
\elist
  In a companion paper \cite{PivatoDefect1}, we developed a new
framework for describing defects, and used spectral theory to
get invariants (as in question \#\ref{Q:inv}) for
codimension-one (`domain boundary') defects.  Unfortunately, these
spectral invariants were not applicable to defects of codimension two
or higher (e.g. `holes' in $\Zahl^2$, `strings' in $\Zahl^3$, etc.).
In this paper, we will answer question \#\ref{Q:inv} for such defects,
using methods inspired by algebraic topology.

  This paper is organized as follows: in \S\ref{S:defect} we
review the framework developed in \cite{PivatoDefect1}.  We also
define defect {\em codimension}, and introduce
many examples which recur throughout the paper.
In \S\ref{S:cocycle}, we address question \#\ref{Q:inv} using 
dynamical cohomology, while in
\S\ref{S:homotopy}, we address \#\ref{Q:inv}
 using tiling homotopy/(co)homology groups. 
In \S\ref{S:cohomology} we relate the dynamical cohomology of
\S\ref{S:cocycle} to the tiling cohomology of \S\ref{S:homotopy}.  In
all cases, we are able to use these algebraic invariants to
answer question \#\ref{Q:persist}, and partially answer question
\#\ref{Q:chem}. 
\begin{center}
\begin{tabular}{rcl}
\hspace{-1em}
\begin{minipage}{22em}
\vspace{-11em}
The diagram at right portrays the logical
dependency of these sections.
In particular, notice that \S\ref{S:cocycle} and \S\ref{S:homotopy}
are logically independent of one another, although \S\ref{S:cohomology}
depends upon both.  Our main results are in sections \ref{S:pole}, 
\ref{S:gaps}, \ref{S:projective.homotopy}, and \ref{S:invariant.cohom}.
\end{minipage}
&\qquad&
\psfrag{S1}[][]{\S\ref{S:defect}}
\psfrag{S2.1}[][]{\S\ref{S:cocycle.intro} }
\psfrag{S2.2}[][]{\S\ref{S:pole} }
\psfrag{S2.3}[][]{\S\ref{S:gaps} }
\psfrag{S3.1}[][]{\S\ref{S:canonical.cell.cplx} }
\psfrag{S3.2}[][]{\S\ref{S:cubic.cohom} }
\psfrag{S3.3}[][]{\S\ref{S:wang.homo} }
\psfrag{S3.4}[][]{\S\ref{S:Wang.representation} }
\psfrag{S3.5}[][]{\S\ref{S:projective.homotopy} }
\psfrag{S4.1}[][]{\S\ref{S:equiv.cohom} }
\psfrag{S4.2}[][]{\S\ref{S:invariant.cohom} }
\includegraphics[scale=0.6]{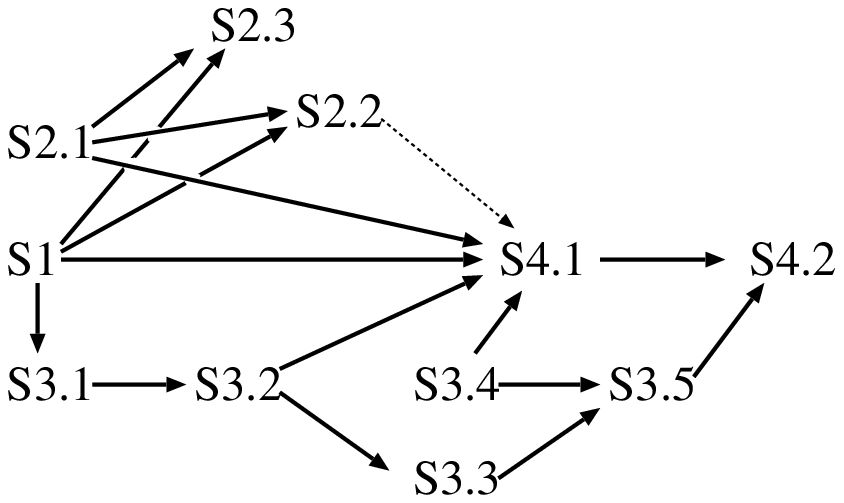}
\end{tabular}
\end{center}

\paragraph*{Preliminaries \& Notation:}
  Let $\sA$ be a finite alphabet.  Let $D\geq 1$, let $\ZD$ be the
$D$-dimensional lattice, and let $\AZD$ be the set of all
$\ZD$-indexed {\dfn configurations} of the form
$\ba=[a_{\fz}]_{\fz\in\ZD}$, where $a_\fz\in\sA$ for all $\fz\in\ZD$.
The {\dfn Cantor metric} on $\AZD$ is defined by
$d(\ba,\bb)=2^{-\Del(\ba,\bb)}$, where
$\Del(\ba,\bb):=\min\set{|\fz|}{a_\fz\neq b_\fz}$.  It follows that
$(\AZD,d)$ is a Cantor space (i.e. a compact, totally disconnected,
perfect metric space).  If $\ba\in\AZD$ and $\dU\subset\ZD$,
then we define $\ba_\dU\in\sA^\dU$ by
$\ba_\dU:=[a_\fu]_{\fu\in\dU}$.  If $\fz\in\ZD$, then strictly speaking,
$\ba_{\fz+\dU}\in\sA^{\fz+\dU}$; however, it will often be convenient
to `abuse notation' and treat $\ba_{\fz+\dU}$ as an element of
$\sA^{\dU}$ in the obvious way.

 For any $\fv\in\ZD$,
we define the {\dfn shift} $\shift{\fv}:\AZD\into\AZD$ by
$\shift{\fv}(\ba)_{\fz} = a_{\fz+\fv}$ for all $\ba\in\AZD$ and
$\fz\in\ZD$.  A {\dfn cellular automaton} is a
transformation $\Phi:\AZD\into\AZD$ that is
continuous and commutes with all shifts.
Equivalently, $\Phi$ is determined by
a {\dfn local rule} $\phi:\sA^\Nh\into\sA$ such that
$\Phi(\ba)_{\fz} = \phi(\ba_{\fz+\Nh})$ for all $\ba\in\AZD$ and
$\fz\in\ZD$ \cite{Hedlund}.   Here, $\Nh\subset\ZD$ is a finite
set which we normally imagine as a `neighbourhood of the origin'.
If $\Nh\subseteq \dB(r):=\CC{-r...r}^D$, we 
say that $\Phi$ has {\dfn radius} $r$.
\ignore{We say that $\Phi$ is a {\dfn nearest neighbour} cellular
automaton (NNCA) if $r=1$.  By replacing $\sA$ with
$\sB=\sA^\Nh$ if necessary, we can `recode' any CA as an NNCA.}

  A subset $\gA\subset\AZD$ is a {\dfn subshift}
\cite{LindMarcus,Kitchens} if $\gA$ is closed in the Cantor topology,
and if $\shift{\fz}(\gA)=\gA$ for all $\fz\in\ZD$.  For any
$\dU\subset\ZD$, we define $\gA_\dU:=\set{\ba_\dU}{\ba\in\gA}$. 
In particular, for any $r>0$, let 
$\Adm := \gA_{\dB(r)}$ be the set of {\dfn admissible
$r$-blocks} for $\gA$.  We say $\gA$ is a {\dfn subshift of
finite type} (SFT) if there is some $r>0$ (the {\dfn radius} of $\gA$)
such that  $\gA$ is entirely described by $\Adm$, in the sense
that $\gA=\set{\ba\in\AZD}{\ba_{\dB(\fz,r)}\in\Adm, \
\forall\fz\in\ZD}$.  If $D=1$, then 
a {\dfn Markov subshift} is an SFT $\gA\subset\sA^\Zahl$ determined
by a set $\gA_{\{0,1\}}\subset\sA^{\{0,1\}}$ of {\dfn admissible transitions}; 
equivalently, $\gA$ is the set of all bi-infinite directed paths in a 
digraph  whose vertices are the elements of $\gA$, with
an edge $a\leadsto b$ iff $(a,b)\in\gA_{\{0,1\}}$. 
If $D=2$, then let $\dE_1:=\{(0,0),(1,0)\}$ and
 $\dE_2:=\{(0,0),(0,1)\}$.  A {\dfn Wang subshift}
 is an SFT $\gA\subset\AZD[2]$
determined by sets $\gA_{\dE_1}\subset\sA^{\dE_1}$ and
$\gA_{\dE_2}\subset\sA^{\dE_2}$ of {\dfn edge-matching conditions}.
Equivalently, $\gA$ is the set of all {\dfn tilings} of the plane $\RD[2]$
by unit square tiles (corresponding to the elements of $\sA$) with
notched edges representing the edge-matching conditions 
\cite[Ch.11]{GrunbaumShephard}.
More generally, if  $D\geq 3$,
then for all $d\in\CC{1...D}$, let $\dE_d:=\{0\}^{d-1}\x\{0,1\}\x\{0\}^{D-d}$.
A {\dfn Wang subshift} is an SFT $\gA\subset\AZD$ determined by sets
$\gA_{\dE_d}\subset\sA^{\dE_d}$ of {\dfn face-matching conditions}
for $d\in\CC{1...D}$.
Equivalently, $\gA$ is the set of tesselations of $\RD$ by unit (hyper)cubes with `notched' faces.

If $\sX$ is any set and $F:\gA\into\sX$ is a function, then $F$ is
{\dfn locally determined} if there is some {\dfn radius} $r\in\Natur$
and some {\dfn local rule} $f:\Adm\into\sX$ such that
$F(\ba)=f(\ba_{\dB(r)})$ for any $\ba\in\gA$.  If $\sX$ is any
discrete space, then $F:\gA\into\sX$ is continuous iff $F$ is locally
determined.  For example, if $\sA$ and $\sB$ are finite sets, then a
(subshift) {\dfn homomorphism} is a continuous, $\shift{}$-commuting
function $\Phi:\BZD\into\AZD$ (e.g. a CA is a homomorphism with
$\sA=\sB$); it follows that $F(\ba)=\Phi(\ba)_0$ is locally
determined.  If $\gB\subset\BZD$ is a subshift of finite type, and
$\Psi:\BZD\into\AZD$ is a homomorphism, then $\gA:=\Psi(\gB)\subset\AZD$
is called a {\dfn sofic shift}. 

  If $\Phi:\AZD\into\AZD$ is a cellular automaton, then we say $\gA$
is  (weakly) {\dfn $\Phi$-invariant} if $\Phi(\gA)\subseteq\gA$
(i.e. $\Phi$ is an {\dfn endomorphism} of $\gA$).  For example,
 if $p\in\Natur$ and $\fv\in\ZD$, then the set
$\Fix{\Phi^p}$ of  {\dfn $(\Phi,p)$-periodic points} and the set
$\Fix{\Phi^p\circ\shift{-p\fv}}$ of {\dfn $(\Phi,p,\fv)$-travelling waves} are
$\Phi$-invariant SFTs.  If $\Phi^\oo(\AZD):=\Intsct_{t=1}^\oo \Phi^t(\AZD)$
is the {\dfn eventual image} of $\Phi$, then $\Phi^\oo(\AZD)$ is a
$\Phi$-invariant subshift (possibly non-sofic), which contains
$\Fix{\Phi^p\circ\shift{-p\fv}}$ for any $p\in\Natur$ and $\fv\in\ZD$.

If $\fy,\fz\in\ZD$, then we write
``$\fy\adjacent\fz$'' if $|\fz-\fy|=1$.  A {\dfn trail} is a sequence
$\zeta=(\fz_1\adjacent\fz_2\adjacent\cdots\adjacent\fz_n)$.  A subset
$\dY\subset\ZD$ is {\dfn trail-connected} if, for any $\fx,\fy\in\dY$,
there is a trail $\fx=\fz_0\adjacent\fz_1\adjacent
\cdots\adjacent\fz_n=\fy$ in $\dY$.

\paragraph{Font conventions:}
Upper case calligraphic letters ($\sA,\sB,\sC,\ldots$) denote 
alphabets or groups.  Upper-case Gothic letters ($\gA,\gB,\gC,\ldots$)
denote subsets of $\AZD$ (e.g. subshifts), lowercase bold-faced letters
($\ba,\bb,\bc,\ldots$) denote elements of $\AZD$, and Roman letters
($a,b,c,\ldots$) are elements of $\sA$ or ordinary numbers.
Lower-case sans-serif ($\ldots,\fx,\fy,\fz$) are elements of $\ZD$,
upper-case hollow font ($\dU,\dV,\dW,\ldots$) are subsets of $\ZD$,
and upper-case bold ($\bU,\bV,\bW,\ldots$) are subsets of
$\RD$. Upper-case Greek letters ($\Phi,\Psi,\ldots$) are functions
on $\AZD$ (e.g. CA), and lower-case Greek letters ($\phi,\psi,\ldots$)
are other functions (e.g. local rules.)

\section{Defects and Codimension\label{S:defect}}

  Let $\gA\subset\AZD$ be any subshift.  If $\ba\in\AZD$, then the
{\dfn defect field} $\energy:\ZD\into\Natur\union\{\oo\}$ is defined
by $\energy(\fz) \ := \ \max\set{r\in\Natur}{\ba_{\dB(\fz,r)} \in \gA_{(r)}}$,
for all $\fz\in\ZD$.
Clearly, $\energy$ is `Lipschitz' in the sense that
$|\energy(\fy)-\energy(\fz)|\leq|\fy-\fz|$.  
The {\dfn defect set} of $\ba$ is the set $\dD(\ba)\subset\ZD$ of local minima of
$\energy$.  See \cite[\S1]{PivatoDefect1} for further discussion.

\example{\label{X:defect.energy}
(a) Suppose $\gA$ is an SFT determined by a set $\gA_{(r)}\subset\sA^{\dB(r)}$
of admissible $r$-blocks, and
let $\dX:=\set{\fz\in\ZD}{\ba_{\dB(\fz,r)}\not\in\gA_{(r)}}$.
Assume for simplicity
that $\gA_{(r-1)}=\sA^{\dB(r-1)}$.  Then 
$\energy(\fz)=r+d(\fz,\dX)$, where
$d(\fz,\dX):=\D \min_{\fx\in\dX} \ |\fz-\fx|$.
In particular, $\energy(\fz)=r$ if and only if $\fz\in\dX$,
and this is the smallest possible value for $\energy(\fz)$.
Thus, $\dD(\ba)=\dX$. 

(b)  Let $\sA=\sB\union\sD$ 
and let $\gA:=\sB^{\ZD}\union\sD^{\ZD}$.  Let $\sC:=\sB\intsct\sD$, 
and let $\sB^*:=\sB\setminus\sC$ and $\sD^*:=\sD\setminus\sC$.
Any $\ba\in\AZD$
is a mixture of $\sB^*$-symbols, $\sC$-symbols, and $\sD^*$-symbols.
If $\fz\in\ZD$ and $a_\fz\in\sB^*$, then
$\energy(\fz)=\min\set{|\fy-\fz|}{a_\fy\in\sD^*}$.  
If $a_\fz\in\sD^*$, then $\energy(\fz)=\min\set{|\fy-\fz|}{a_\fy\in\sB^*}$.
If $a_\fz\in\sC$, then $\energy(\fz)=
\min\{r \ ; \ a_\fx\in\sB^* \And a_\fy\in\sB^* $ $  \ \mbox{for some 
$\fx,\fy\in\dB(\fz,r)$}\}$.  Thus, $\dD(\ba)$ is the set of all points 
which are either on a `boundary' between a $\sB^*$-domain and a $\sD^*$-domain,
or roughly in the middle of a $\sC$-domain.
}
Let $\tlgA:=\set{\ba\in\AZD}{\D\sup_{\fz\in\ZD} \ \energy(\fz) \ = \ \oo}$
be the set of `slightly defective' configurations.   
If $\ba\in\tlgA\setminus\gA$, then we say $\ba$ is {\dfn defective}.
Elements of $\tlgA$ may have infinitely large defects, but also have
arbitrarily large non-defective regions.
Clearly $\gA\subset\tlgA$, and $\tlgA$ is a $\shift{}$-invariant,
dense subset of $\AZD$ (but not a subshift).

For any $R>0$, let
$\unflawed[R](\ba):=\set{\fz\in\ZD}{\energy(\fz)\geq R}$.  Thus,
$\ba\in\tlgA$ iff $\unflawed[R](\ba)\neq\emptyset$ for all $R>0$.  For
example, if $\gA$ is an SFT determined by a set $\gA_{(r)}$ of
admissible $r$-blocks, and
$\dD=\set{\fz\in\ZD}{\ba_{\dB(\fz,r)}\not\in\gA_{(r)}}$ as in Example
\ref{X:defect.energy}(a), then
$\unflawed[R](\ba)=\set{\fz\in\ZD}{d(\fz,\dD)\geq R-r}
=\ZD\setminus\dB(\dD,R-r)$.  Thus, $\dD(\ba)$ encodes all information
about the `defect structure' of $\ba$.  However, if $\gA$ is {\em not}
an SFT [e.g. Example \ref{X:defect.energy}(b)], then in general
$\unflawed[R](\ba) \neq\ZD\setminus\dB(\dD,R')$ for any $R'>0$.  In
this case, $\dD(\ba)$ is an inadequate description of the larger-scale
`defect structures' of $\ba$.  Thus, instead of treating the defect as
a precisely defined subset of $\ZD$, it is better to think of it as a
`fuzzy' object residing in the low areas in the defect field $\energy$.  The
advantage of this approach is its applicability to any kind of
subshift (finite type, sofic, or otherwise).  Nevertheless, most of
our examples will be SFTs, and we may then refer to the specific
region $\dD\subset\ZD$ as `the defect'.

\Proposition{\label{defect.dimension}}
{
 Let $\Phi:\AZD\into\AZD$ be a CA with radius $r>0$.
\bthmlist
  \item Let $\gA\subset\AZD$ be a weakly $\Phi$-invariant subshift. 
Then $\Phi(\tlgA)\subseteq\tlgA$.

  \item If $\ba\in\tlgA$, and $\ba'=\Phi(\ba)$,
then $\energy[\ba'] \geq \energy[\ba]-r$.
Thus, for all $R\in\Natur$, \ 
 $\unflawed[R+r](\ba)\subseteq \unflawed[R](\ba')$.
\ethmlist
}
\bthmprf  {\bf(b)} Let $\fz\in\ZD$ and suppose $\energy(\fz)=R$.  Thus,
$\ba_{\dB(\fz,R)}\in\gA_R$.  But $\gA$ is $\Phi$-invariant;
hence $\ba'_{\dB(\fz,R-r)}\in\gA_{(R-r)}$.  Hence $\energy[\ba'](\fz)\geq R-r$.
Then {\bf(a)} follows from {\bf(b)}.
\ethmprf

 Let $\Phi:\AZD\into\AZD$ be a cellular automaton, and suppose that
$\phi(\gA)=\gA$.  If $\ba\in\tlgA$ then
$\ba$ has a {\dfn $\Phi$-persistent} defect if, for all
$t\in\Natur$, \  $\ba'=\Phi^t(\ba)$ is also defective.  Otherwise $\ba$
has a {\dfn transient} defect ---i.e. one which eventually disappears. 
We say $\ba$  has a {\dfn removable} defect if there is some $r>0$ and some
$\ba'\in\gA$ such that $a'_\fz=a_\fz$ for all $\fz\in\unflawed(\ba)$
(i.e. the defect can be erased by modifying $\ba$ in a finite radius
of the defective region).  Otherwise $\ba$ has an {\dfn essential} defect.

\example{\label{X:removable.finite.defects}
The defect in $\ba$ is {\dfn finite} if 
 $\dD(\ba)$ is finite (or equivalently
$\D \lim_{|\fz|\goto\oo}\ \energy(\fz)=\oo$).

(a) Let $\gA\subset\AZ$.  Then $(\gA,\shift{})$ is topologically mixing
if and only if no finite defect is essential.

(b) $\gA\subset\AZD$ satisfies the {\dfn hole-filling property}
if no finite defect is essential (i.e. every configuration with
a finite defect is `weakly extensible' in the sense of \cite[\S5]{Sch98}).
}

\Proposition{\label{biject.essential.persistent}}
{
Let $\Phi:\AZD\into\AZD$ be a CA and let
$\gA\subset\AZD$ be a $\Phi$-invariant subshift.
If $\Phi:\gA\into\gA$ is bijective, then any essential defect is 
$\Phi$-persistent.

In particular, if $\gA\subseteq\Fix{\Phi}$ or 
 $\gA\subseteq\Fix{\Phi^p\circ\shift{p\fv}}$ 
{\rm(for some $p\in\Natur$ and $\fv\in\ZD$)},
then any essential defect in $\ba$ is $\Phi$-persistent.\qed{\rm\cite{PivatoDefect1}}
}

\subsection{Codimension:\label{S:codimension}}
  Our main goal in the present paper is to develop algebraic
invariants (as described by question \#\ref{Q:inv} from the
introduction) which provide sufficient conditions for the persistence
of defects, even when $\Phi$ is not bijective.  To do this, we must
first assign a `codimension' to defects, but in a somewhat indirect
fashion.  Strictly speaking, the defect set $\dD(\ba)\subset\ZD$ is
discrete, hence of codimension $D$ in $\RD$.  We could `thicken' $\dD$
by replacing each point $\fd\in\dD$ with a unit cube around $\fd$.
However, the cellular automaton $\Phi$, the subshift $\gA$, and other
gadgets we require (e.g. eigenfunctions, cocycles) may have interaction
ranges greater than one (and possibly unbounded), so a unit cube isn't
big enough.  Furthermore, the action of $\Phi$ may locally change the
geometry of the defect, and we are mainly interested in properties
that are invariant under such change (as in the definition of
`essential' defects, above).  Loosely speaking, we will use the word
`projective' to describe `large scale' geometric properties which
remain visible when seen from `far away' (precise definitions will
appear below).

For any $r>0$ and $\ba\in\tlgA$, we say $\ba$ has a {\dfn range $r$
domain boundary} (or a {\dfn range $r$ codimension-one defect}) if
$\unflawed(\ba)$ is trail-disconnected. Domain boundaries divide $\ZD$
into different `domains', which may correspond to different transitive
components of $\gA$ \cite[\S2]{PivatoDefect1}, different eigenfunction
phases \cite[\S3]{PivatoDefect1}, or different cocycle asymptotics
(\S\ref{S:gaps}).  A connected component $\dY$ of $\unflawed(\ba)$ is
called {\dfn projective} if $\dY\intsct\unflawed[R](\ba)\neq\emptyset$
for all $R\geq r$.  (This implies that for any $R\geq 0$, there 
exists $\fy\in\dY$ with $\dB(\fy,R)\subset\dY$. If $\gA$ is of finite
type, then the two conditions are equivalent.)
We say that $\ba$ has a {\dfn projective domain boundary} (or a {\dfn
projective codimension-one defect}) if there is some $R\geq0$ such
that $\unflawed[R](\ba)$ has at least two projective components.
(Hence $\unflawed(\ba)$ is disconnected for all $r\geq R$.)

\begin{figure}
\centerline{\footnotesize
\begin{tabular}{cc}
\includegraphics[scale=0.55]{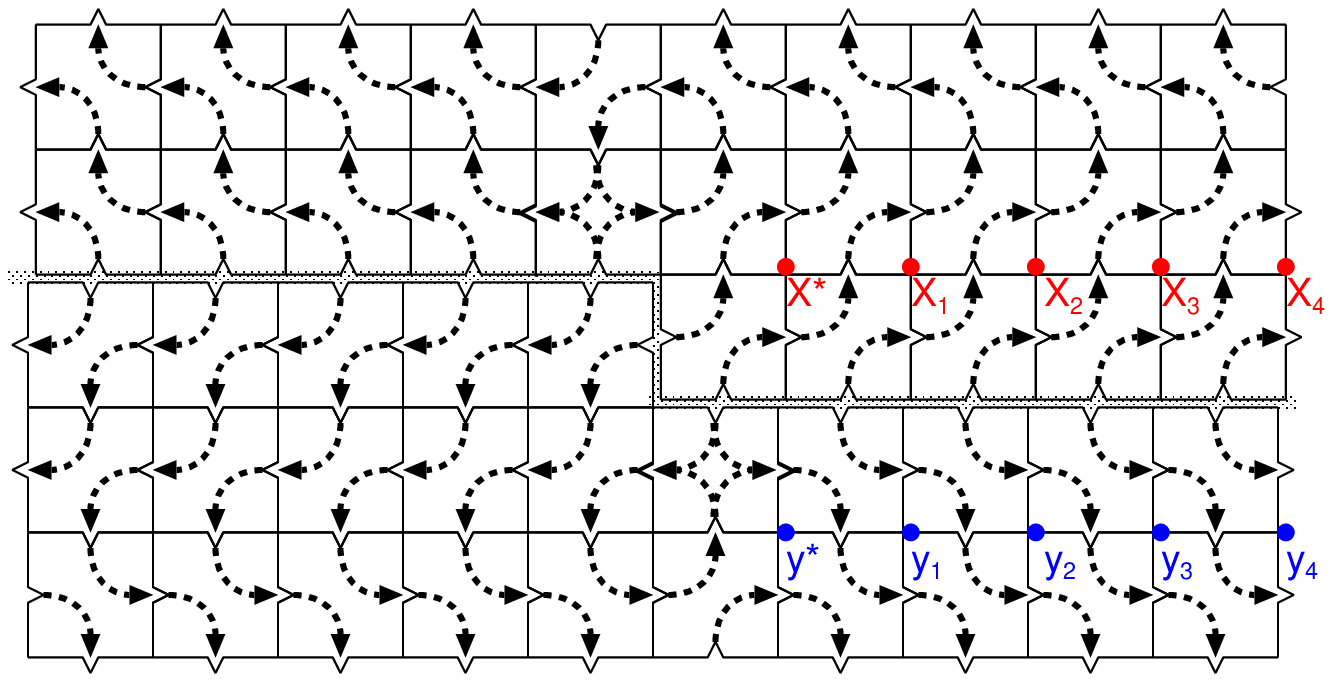} &
\includegraphics[scale=0.6]{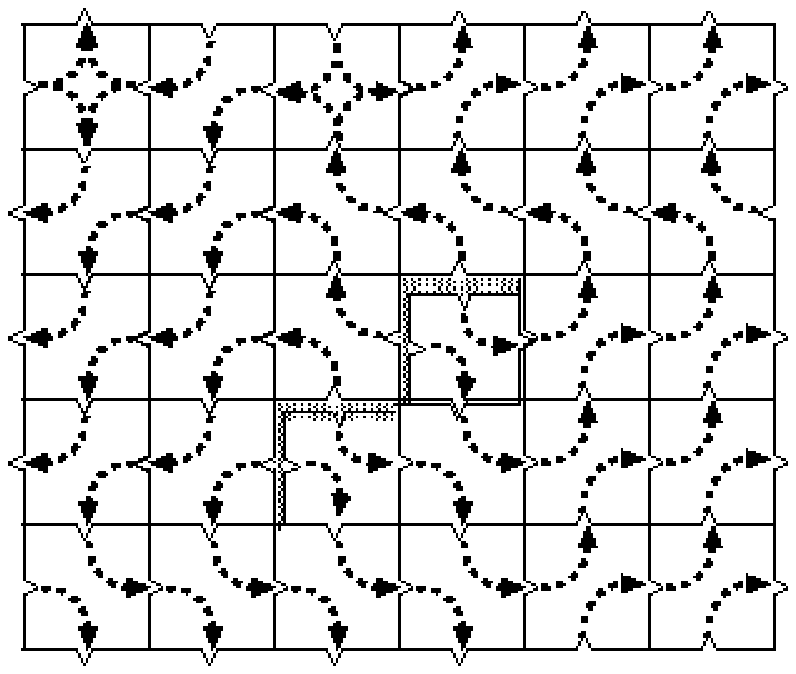} \\
{\bf(A)} & {\bf(D)} \\
\includegraphics[scale=0.6]{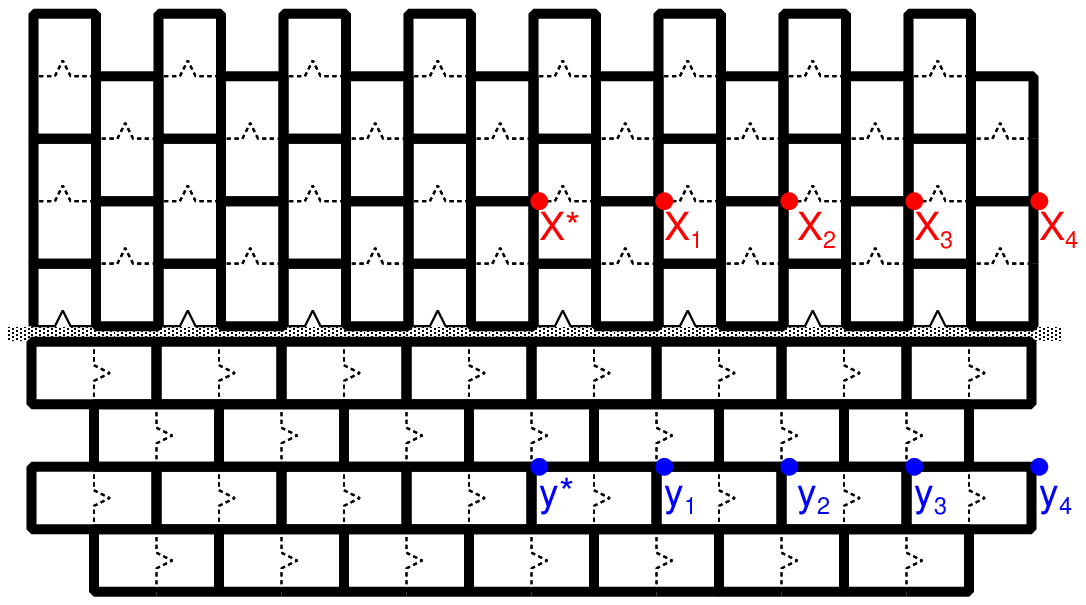}&
\includegraphics[scale=0.55]{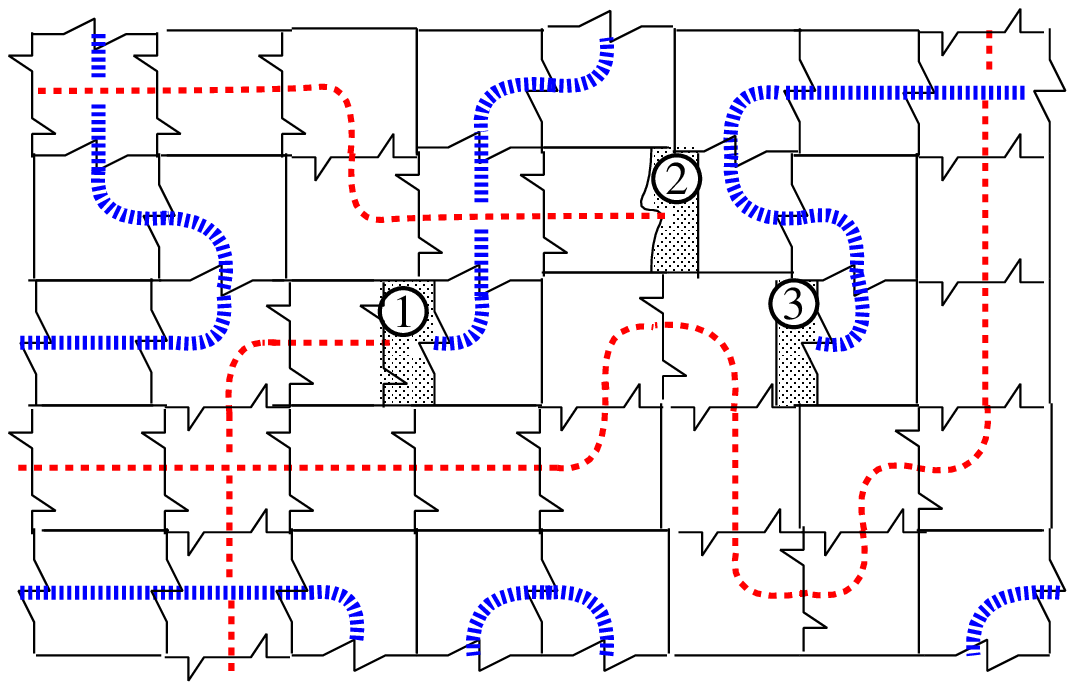}\\
{\bf(B)} & {\bf(E)} \\
\includegraphics[scale=0.6]{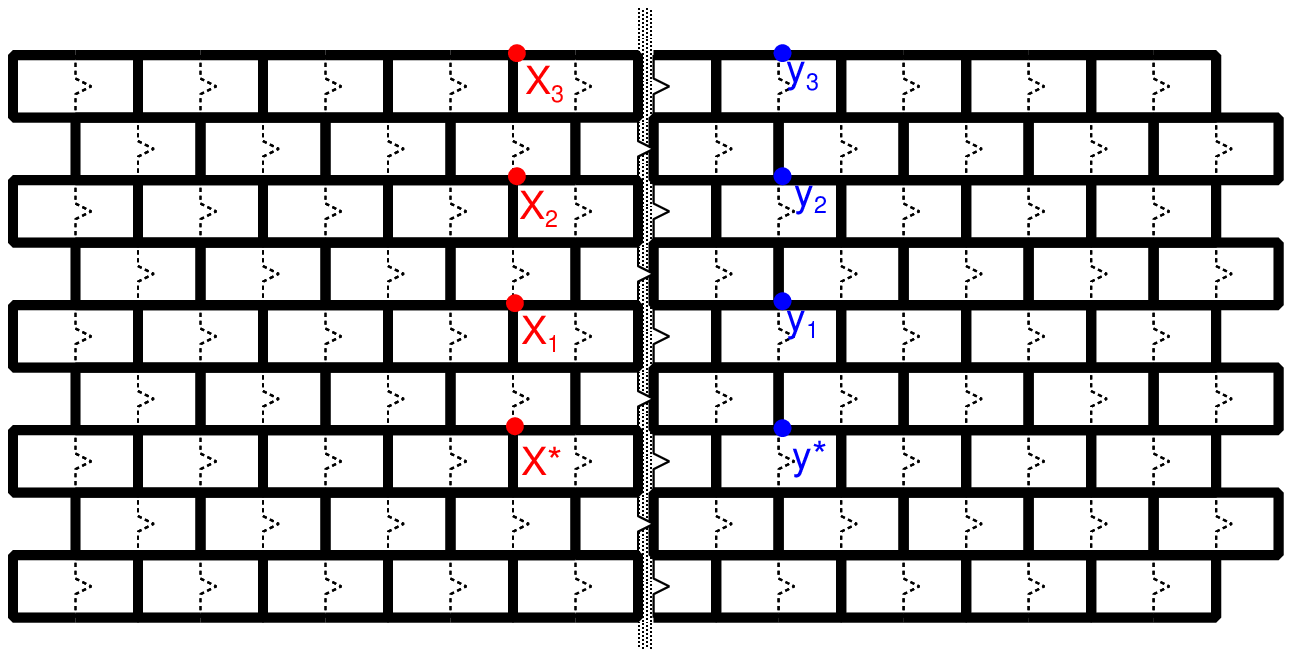}&
\includegraphics[scale=0.6]{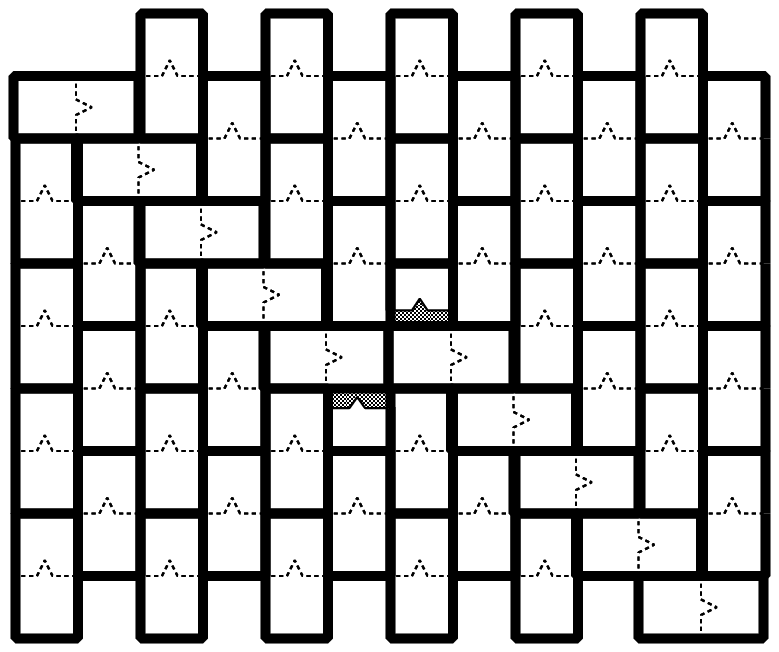} \\
{\bf(C)} & {\bf(F)} \\
\end{tabular}}
\caption{{\footnotesize 
{\bf(A)} \ A gap in $\Ice$;
 see Examples \ref{X:codim.one}(a) and  \ref{X:ice.gap}(a).
{\bf(B,C)} \  Gaps in $\Dom$;
 see Examples \ref{X:codim.one}(b) and  \ref{X:ice.gap}(b,c).
 \label{fig:codim}
{\bf(D)} A `pole' in $\Ice$;  see Examples \ref{X:codim.two}(a) and \ref{X:ice.pole}(a). 
{\bf(E)} Three `poles' in $\Pth$;  see Examples \ref{X:codim.two}(d) and \ref{X:ice.pole}(b).
{\bf(F)} A non-pole in $\Dom$; see Examples \ref{X:codim.two}(b) and \ref{X:ice.pole}(c).
\label{fig:path.tiling}}}
\end{figure}

\begin{figure}
\centerline{\includegraphics[scale=0.25]{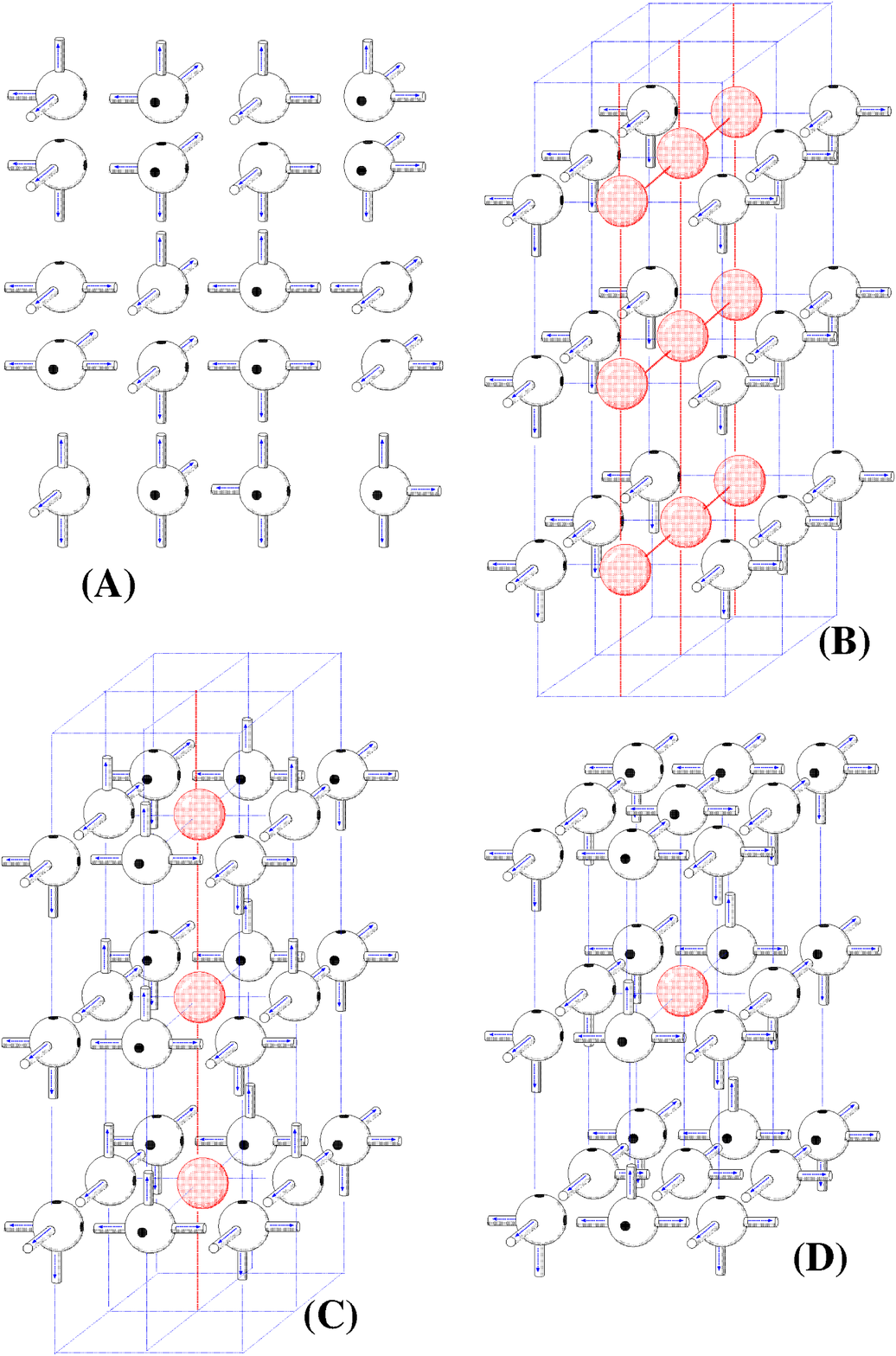}}
\caption{{\footnotesize {\bf(A)} The twenty tiles of the `ice cube'
shift $\gQ\subset\sQ^{\ZD[3]}$ from Examples \ref{X:codim.one}(c)
and \ref{X:codim.two}(c).
{\bf(B)} \ A domain boundary in $\gQ$;  we assume that the same pattern is
continued up, down, north, and south.
(We have `stretched out' the configuration for visibility.
The balls in the defect region are shaded, but left
unspecified, because they don't matter.)
{\bf(C)} \ A codimension-two defect in $\gQ$;
we assume that the same pattern is continued upwards and downwards.
{\bf(D)} \ A codimension-three (`pole') defect in $\gQ$; see also
Example \ref{X.d.pole}(b).}
\label{fig:ice.cube}}
\end{figure}

\example{\label{X:codim.one}
  (a) ({\em Square ice})
  Let $\sI=\lb\{\raisebox{-0.5em}{\includegraphics[scale=0.35]{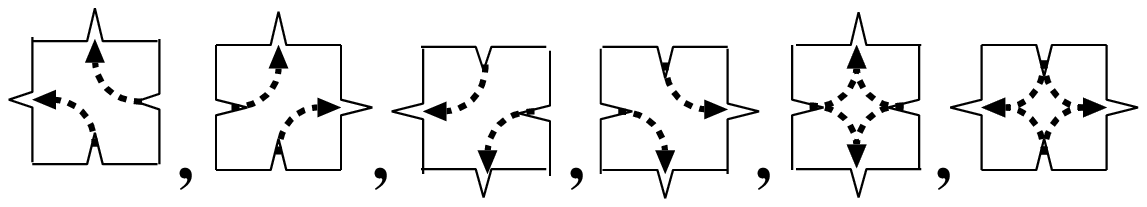}}\rb\}$, 
and let $\Ice\subset\sI^{\ZD[2]}$ be the Wang subshift defined by the obvious
edge-matching conditions.  Figure \ref{fig:codim}(A) shows a domain boundary in $\Ice$.  See also Example \ref{X:ice.gap}(a).

  (b) ({\em Domino Tiling})  Let $\sD:=\lb\{\raisebox{-0.5em}{
\includegraphics[scale=0.65]{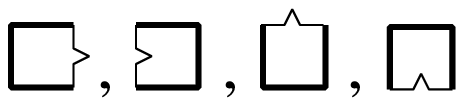}}\rb\}$,
 and let $\Dom\subset\sD^{\ZD[2]}$ 
be the Wang subshift defined by the obvious edge-matching conditions.
Figure \ref{fig:codim}(B,C)  shows two domain boundaries in $\Dom$;
see also Example \ref{X:ice.gap}(b,c).

  (c) ({\em Ice cubes}) Let $\sQ$ be the set of twenty `ball-and-pin'
structures in Figure \ref{fig:ice.cube}(A), and let $\gQ\subset\sQ^{\ZD[3]}$
be the Wang subshift defined by the obvious matching conditions (this is
a three-dimensional version of `square ice').    
Figure \ref{fig:ice.cube}(B) shows a domain boundary in $\gQ$.
}

For any $\fz\in\ZD$, let $\CornerTile{\fz}{ \ }:=\fz+\CC{0,1}^D$; hence $\CornerTile{\fz}{ \ }$ is
a unit cube with one corner at $\fz$, and the other corners at adjacent
points in $\ZD$.  Let $\dK_\fz\subset\ZD$ be the set of corner
points of $\CornerTile{\fz}{ \ }$.  
We adopt the following notational convention:
if $\dY\subset\ZD$ is any subset, then 
let $\bY$ be the minimal closed subset of $\RD$ containing $\dY$ and all
unit cubes whose corners are in $\dY$.  Formally:
\[
\bY\quad:=\quad \Union_{\fz\in\ZD \ \& \ \dK_\fz\subset\dY} \CornerTile{\fz}{ \ }
\]
It follows that $\dY$ is trail-connected iff $\bY$ is path-connected.
In this case,
for any $k\in\CC{2...D}$, we define the {\dfn $k$th homotopy group}
$\pi_k(\dY,\fy):=\pi_k(\bY,\fy)$, for some fixed basepoint $\fy\in\dY$
(different choices of $\fy$ yield isomorphic groups); see 
\cite[\S4.1]{Hatcher}.
If $\ba\in\tlgA$,  and $r>0$,
then $\ba$ has a {\dfn range $r$ codimension-$k$} defect if
$\pi_{k-1}(\unflawed(\ba),\fy)$ is nontrivial for some $\fy\in\unflawed(\ba)$.
If $\unflawed(\ba)$ is disconnected (e.g. by a domain boundary) then
different connected components may have different homotopy groups;
we only require one of these to be nontrivial.

\example{\label{X:codim.two}
(a) Let $\Ice$ be as in Example
\ref{X:codim.one}(a).  Then Figure \ref{fig:path.tiling}(D) shows 
a codimension-two defect in $\Ice$.  See also  Example \ref{X:ice.pole}(a).

  (b) Let $\Dom$ as in
Example \ref{X:codim.one}(b).  Then Figure \ref{fig:path.tiling}(F)
shows  a codimension-two
defect in $\Dom$.  See also  Example \ref{X:ice.pole}(c).

  (c) Let $\gQ$ be as in Example \ref{X:codim.one}(c).
Figure \ref{fig:ice.cube}(C) shows a codimension-two defect in $\gQ$, and
Figure \ref{fig:ice.cube}(D) shows a codimension-three defect.

\centerline{
$\sP \ := \ \lb\{
\raisebox{2.6em}{\includegraphics[scale=0.4,angle=-90]{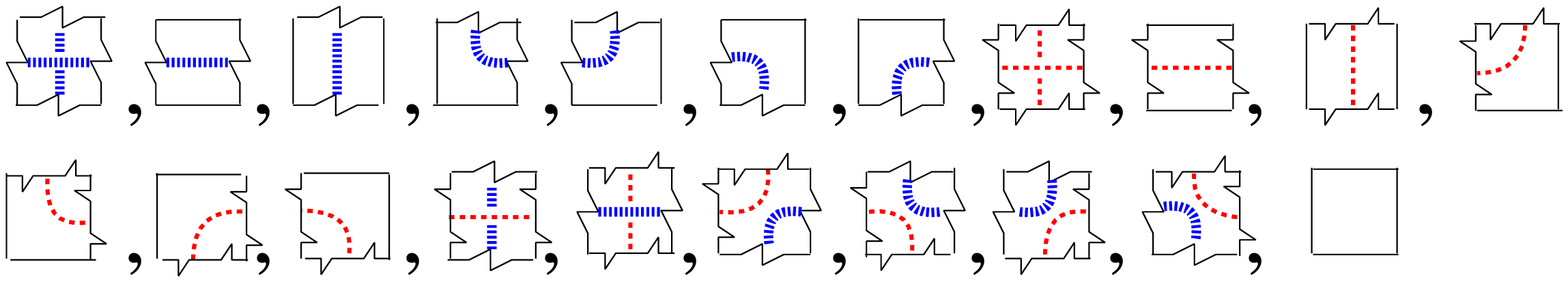}}
\rb\}.$}

  (d)  ({\em Two-coloured, undirected, crossing path tiling})
Let $\sP$ be the set of 21 tiles shown above,
and let $\Pth\subset\sP^{\ZD[2]}$ be the Wang subshift defined by the obvious
edge-matching conditions.  Then $\Pth$-admissible configurations
are tangles of undirected, freely crossing paths in two colours
 \cite[\S3]{Ein}.  Figure \ref{fig:path.tiling}(E)
shows three codimension-two defects in $\Pth$. 
See also Example \ref{X:ice.pole}(b).
}

\subsection{Proper homotopy and projective codimension: \label{S:proj.codim}}
  Let $\bX$ be a topological space and let
$x\in\bX$.  Let $\dS^k\subset\RD[k+1]$ be the unit $k$-sphere, and let
$s\in\dS^k$ be some distinguished point.  We write
$\alp:(\dS^k,s)\into(\bX,x)$ to mean $\alp$ is a continuous function
from $\dS^k$ into $\bX$ and $f(s)=x$.  If $\alp,\bet:(\dS^k,s)\into(\bX,x)$
then we write $\alp \homoto \bet$ to mean that $\alp$ is homotopic to
$\bet$ in a manner which always maps $s$ to $x$; we call this
a {\dfn basepoint-fixing homotopy} (where $x$ is the {\dfn basepoint}).
We then use $\undalp$ to refer to the (basepoint-fixing) homotopy
class of $\alp$.

If $\bet:\CC{0,1}\into\bX$, then $\bkw{\bet}:\CC{0,1}\into\bX$
is defined by $\bkw{\bet}(t)=\bet(1-t)$.  If $\alp:\CC{0,1}\into\bX$,
and $\bet(0)=\alp(1)$, then 
let $\alp\star\bet:\CC{0,1}\into\bX$ be the concatenation of $\alp$ and
$\bet$ (i.e. $\alp\star\bet(t):=\alp(2t)$ if $t\in\CC{0,\frac{1}{2}}$ and 
$\alp\star\bet(t):=\bet(2t-1)$ if $t\in\CC{\frac{1}{2},1}$).   Thus, $\pi_1(\bX,x)$ is the group
of all homotopy classes of loops $\alp:\CC{0,1}\into\bX$
with $\alp(0)=x=\alp(1)$, with operation $\undalp\cdot\undbet
:=\underline{\alp\star\bet}$ \cite[\S1.1]{Hatcher}.
 We can also treat the elements of $\pi_1(\bX,x)$
as homotopy classes of functions $\alp:(\dS^1,s)\into(\bX,x)$.
By generalizing this construction, we can define an abelian
group $\pi_k(\bX,x)$  of homotopy classes of functions 
 $\alp:(\dS^k,s)\into(\bX,x)$.  See \cite[\S4.1]{Hatcher} for details.

  Let $x,y\in\bX$.  If $\bet:\CC{0,1}\into \bX$ is any path with with
$\bet(0)=x$ and $\bet(1)=y$, then $\bet$ yields an isomorphism
$\bet_*:\pi_1(\bX,x)\into\pi_1(\bX,y)$ by $\bet_*(\undalp):=
\underline{\bet\star\alp\star\bkw{\bet}}$.  We can likewise use $\bet$ to
define isomorphisms $\bet_*:\pi_k(\bX,x)\into\pi_k(\bX,y)$ for all
$k\geq 2$.  If $\gam:\CC{0,1}\into\bX$ is another path from $x$
to $y$, and $\gam\homoto \bet$, then $\gam_*=\bet_*$.  However, if
$\gam$ is {\em not} homotopic to $\bet$, then $\gam_*$ and $\bet_*$ may be
different.  Hence, although $\pi_k(\bX,x)\cong\pi_k(\bX,y)$,
this isomorphism is not `canonical'. 

\breath
  
  Let $\gA\subset\AZD$ be a subshift and let $\ba\in\tlgA$,
and suppose, for some $r_0\in\Natur$, that
$\unflawed[r_0](\ba)$ contains a unique projective connected component
$\dY$ (it is necessary, but not sufficient, to assume that $\ba$ has no
projective domain boundaries).  Thus, for all $r>r_0$,
$\dY_r:=\dY\intsct\unflawed(\ba)$ is the unique projective component
of $\unflawed(\ba)$.  A {\dfn proper base ray} is a continuous path $\omg:\CO{0,\oo}\into\bY_{r_0}$ with $\lim_{t\goto\oo} |\omg(t)|=\oo$.
\ignore{such that $\omg^{-1}[\bK]$ is compact in $\CO{0,\oo}$
for every compact subset $\bK\subset\bY_{r_0}$; hence} 
For each $r>r_0$, we define $\pi_k(\dY_r,\omg):=\pi_k(\dY_r,\fy)$,
where $\fy\in\dY_r\intsct\omg\CO{0,\oo}$ is any point.  
This definition is independent of the choice of $\fy$ in the following sense:
if $\fy'\in\dY_r\intsct\omg\CO{0,\oo}$ is another point, then
there is a canonical isomorphism $\pi_k(\dY_r,\fy)\cong\pi_k(\dY_r,\fy')$
given by the segment of $\omg$ between $\fy$ and $\fy'$.

Recall that $\dY_{r+1}\subset\dY_r$; \  the inclusion map
$\iota_r:\dY_{r+1}\hookrightarrow \dY_r$ yields a (canonical)
homomorphism $\iota_r^*:\pi_k(\dY_{r+1},\omg) \into
\pi_k(\dY_r,\omg)$. We define the {\dfn $k$th proper
homotopy group} to be the inverse limit: 
\beqn
\label{projective.homotopy.defn}
 \pi_k(\unflawed[\oo](\ba),\omg)\quad:=\quad \invlim
\lb( \pi_k(\dY_1,\omg) \stackrel{\iota^*_1}{\longleftarrow\!\!\!-}
\pi_k(\dY_2,\omg) \stackrel{\iota^*_2}{\longleftarrow\!\!\!-}
\pi_k(\dY_3,\omg) \stackrel{\iota^*_3}{\longleftarrow\!\!\!-}
\cdots\rb)
\eeqn
(Of course ``$\pi_k(\unflawed[\oo](\ba))$'' is an abuse of notation, because
technically, $\unflawed[\oo](\ba)=\emptyset$.  
See \cite[\S3.F]{Hatcher} or \cite[\S III.9]{Lang} for background on inverse
limits.
The group $\pi_k(\unflawed[\oo](\ba),\omg)$ is analogous 
(but not identical) to the proper
homotopy group of a noncompact topological space; see 
\cite{Brown}, \cite{BrownTucker} or \cite[\S2]{Peschke}.)
We say that $\ba$ has {\dfn projective codimension $(k+1)$ defect} if 
$ \pi_k[\unflawed[\oo](\ba),\omg]$ is nontrivial 
(it follows that $\pi_k[\unflawed(\ba),\omg]$ is nontrivial for all
large enough $r\in\Natur$).
Heuristically,
elements of $\pi_k(\unflawed[\oo](\ba))$  are homotopy classes of `extremely
large' $k$-sphere embeddings in the unflawed part of $\ba$.
Technically, this definition depends upon the homotopy class
of the proper base ray $\omg$; different rays may yield nonisomorphic groups.

\section{Cohomological Defects\label{S:cocycle}}

The main results of this section are Theorems \ref{cohom.defect}
and \ref{persistent.gap} and Proposition \ref{persistent.pole}.

\subsection{Dynamical Cocycles\label{S:cocycle.intro}}

  Let $\Sft\subseteq\AZD$ be a subshift, and let $(\sG,\cdot)$ be a
topological group (usually discrete).  A $\sG$-valued continuous 
({\dfn dynamical}) {\dfn cocycle} for $\Sft$ is a continuous function $C:\ZD\x\Sft\into\sG$
satisfying the {\dfn cocycle equation}
\beqn
\label{cocycle.eqn}
C(\fy+\fz,\ba) \quad =\quad C(\fy,\shift{\fz}(\ba))\cdot  C(\fz,\ba),
\qquad\forall\,\ba\in\AZD \And \forall \,\fy,\fz\in\ZD.
\eeqn
\example{\label{X:cocycle}
(a) If $b:\Sft\into\sG$ is any continuous function, then
the function $C(\fz,\ba):=b(\shift{\fz}(\ba))\cdot b(\ba)^{-1}$ is
a cocycle, and is called a {\dfn coboundary} with {\dfn cobounding function} $b$.

(b) If $h:\ZD\into\sG$ is a homomorphism, then the function
$C(\fz,\ba):=h(\fz)$ is a cocycle.  Conversely, if $C$ is any cocycle
such that $C(\fz,\blankspace)$ is constant for all $\fz\in\ZD$, then $C$
arises from a homomorphism in this manner.  In particular, if
$e_\sG\in\sG$ is the identity, then the constant function $C_e(\fz,\ba)\equiv
e_\sG$
is a cocycle.

(c) Let $\Ice\subset\sI^{\ZD[2]}$ be as in Example
\ref{X:codim.one}(a). Define $c_1,c_2:\sI\into \{\pm1\}$ by
$c_1(\WangTile{*}{*}{*}{\sqn})  :=  +1  =: 
c_2(\WangTile{\sqw\ }{*}{*}{*})$ and $c_1(\WangTile{*}{*}{*}{\sqs})  :=
 -1  =:  c_2(\WangTile{\sqe}{*}{*}{*})$ (`$*$' means
`anything').  We define cocycle $C:\ZD[2]\x\Ice\into\Zahl$ as follows:
If $\bi\in\Ice$ and $\fz=(z_1,z_2)\in\ZD[2]$, then
\beqn
\label{height.func}
C(\fz,\bi) \quad := \quad \sum_{x=0}^{z_1-1} c_1(i_{x,0}) 
+ \sum_{y=0}^{z_2-1} c_2(i_{z_1,y}).
\eeqn 
(d) Let  $\Pth\subset\sP^{\ZD[2]}$ be as in
 Example \ref{X:codim.two}(d).  Define
$c_1,c_2:\sI\into \Zahlmod{2}\dirsum\Zahlmod{2}$ by
$c_1(\westbluepath)  :=  (1,0)  =: 
c_2(\southbluepath)$ and $c_1(\westredpath)  :=
 (0,1)  =:  c_2(\southredpath)$ (where `$*$' means
`anything'),
and extend this to a cocycle
$C:\ZD[2]\x\Pth\into\Zahlmod{2}\dirsum\Zahlmod{2}$ exactly as in
eqn.(\ref{height.func}).

(e) More generally, a {\dfn height function} is any integer-valued\footnote{Sometimes $H$ maps into $\Zahl^n$ \cite{Shef} or other groups \cite{Keny}, and the `height' metaphor is somewhat strained.} cocycle
$H:\ZD\x\gA\into\Zahl$ defined by functions $h_1,\ldots,h_D:\sA\into\Zahl$
via the obvious generalization of eqn.(\ref{height.func}).
Height functions appear in tiling systems like dominos \cite{CNPr,CKPr} and
`ice' tiles \cite{Elo99,Elo03,Elo05}, and in
many lattice models of statistical physics \cite{Baxter}.

(f) Let $\Dom\subset\sD^{\ZD[2]}$ be as in Example \ref{X:codim.one}(b).
Let $\sG:=\Zahlmod{2}*\Zahlmod{2}$ be the group of
finite products like $vhvhv\cdots vhv$, where
 $v$ and $h$ are noncommuting generators with $v^2=e=h^2$.  Define
$c_1,c_2:\sI\into \sG$ by

\centerline{$c_1(\WangTile{|\ }{-}{\,|}{\sqn})  := vhv;\quad  
c_1(\WangTile{*}{*}{*}{-})  :=  h;\quad
c_2(\WangTile{\sqe }{-}{\,|}{-})  := hvh; \ \mbox{and} \ 
c_2(\WangTile{|\ }{*}{*}{*})  := v$.}

and  extend this to a cocycle
$C:\ZD[2]\x\Dom\into\sG$ through the 
multiplicative analogy of eqn.(\ref{height.func}).

(g) If $\sX$ is a topological space, then an {\dfn $\sX$-extension}
of $\gA$ is a continuous $\ZD$-action 
$\Xi:\ZD\x\sX\x\gA\into\sX\x\gA$ such that $(\gA,\shift{})$ is a factor
of $(\sX\x\gA,\Xi)$ via the projection $\pi_\gA:\sX\x\gA\surject\gA$.
Let $\sG:=\Homeo(\sX)$ be the self-homeomorphism group of $\sX$,
topologized as a subspace of the Tychonoff product $\sX^\sX$
(e.g. if $\sX:=\CC{1...n}$, then $\sG=\bS_n$ is a (discrete) permutation group;
this is called an {\dfn $n$-point extension}).
For each $\ba\in\gA$ and $\fz\in\ZD$, let $c(\fz,\ba):=
\pi_\sX\circ\Xi^\fz(\blankspace,\ba):\sX\into\sX$.  Then $c:\ZD\x\gA\into\sG$ is a continuous cocycle.
Conversely, any continuous cocycle $c:\ZD\x\gA\into\sG$ defines an
$\sX$-extension of $\gA$ in the obvious way; 
see e.g. \cite{Zimmer1,Zimmer2,Zimmer3,Kam90,Kam92,Kam93}.
}

\ignore{
(f) Let $\sT$ be a set of Wang tiles, and let $\Gam_0=\Gam_0(\sT)$
be the tiling path group (see \S\ref{}).\marginpar{*}  We define
$c_1,c_2:\sT\into\Gam_0$ by 
$c_1(\WangTile{w}{n}{e}{s}) \ := \   s$ and
$c_2(\WangTile{w}{n}{e}{s})\ :=\  w$.
If $\gT\subset\TZD[2]$ is the set of admissible tilings, then we define
the {\dfn tiling cocycle}  (\cite[\S4]{Sch98},\cite{Ein})
$C:\gT\x\ZD[2]\into\Gam_0$ in a
manner analogous to eqn.(\ref{height.func}).  In other words,
for any
$\fz\in\ZD[2]$ and $\bT\in\TZD$, if $\bT=[\tau_{\fz}]_{\fz\in\ZD}$, then
$C(\fz,\ba) \  := \ \prod_{x=0}^{z_1-1} c_1(a_{x,0}) 
\cdot \prod_{y=0}^{z_2-1} c_1(a_{z_1,y})$.
  The relations defining $\Gam_0$ are exactly what is required to make $C$
satisfy the cocycle eqn.(\ref{cocycle.eqn}).

 If $\Sft\subset\AZD[2]$ is any SFT, and
$\sT$ is a set of Wang tiles representing $\Sft$,
with isomorphism $q:\Sft\into\gT$, we
define a {\dfn tiling cocycle} 
$c':\Sft\x\ZD\into\Gam_0$ by $c'(\bx,\fz) = c(q(\bx),\fz)$ for all
$\bx\in\Sft$.   Any continuous cocycle on $\Sft$ taking
values in a discrete group is a homomorphic image of such a tiling
cocycle \cite[Thm.4.2(b)]{Sch98}.  Conversely, $\Gam_0(\sT)$
is the image of the tiling cocycle.
Hence, cocycles and tiling groups encode the same
information.}

Two continuous cocycles $C$ and $C'$ are {\dfn cohomologous} 
($C\approx C'$) if there
is a continuous {\dfn transfer function} $b:\gA\into\sG$ such that
$C'(\fz,\ba) \ = \ b(\shift{\fz}(\ba))\cdot C(\fz,\ba)\cdot
b(\ba)^{-1}$, for all $\fz\in\ZD$ and $\ba\in\gA$. 
A cocycle $C$ is {\dfn trivial} if $C$ is
cohomologous to a homomorphism.  We will use $\undC$ to denote the
cohomology equivalence class of the cocycle $C$.

\example{\label{shift.cohomologous}
(a) Any coboundary [Example \ref{X:cocycle}(a)] is trivial, because it is cohomologous to the homomorphism $C_e$ [Example \ref{X:cocycle}(b)]. 

(b) Fix $\fy\in\ZD$ and define cocycle $C'(\fz,\ba):=C(\fz,\shift{\fy}(\ba))$. 
Then $C\approx C'$ via the transfer function $b(\ba):=C(\fy,\ba)$.

(c) Let $\Xi,\Xi':\ZD\x\sX\x\gA\into\sX\x\gA$ be two $\sX$-extensions
of $\gA$, with cocycles $C,C':\ZD\x\gX\into\sG:=\Homeo(\sX)$ as in
Example \ref{X:cocycle}(g).  Then $C\approx C'$ via the continuous transfer
function $b:\gA\into\sG$ iff the systems $(\sX\x\gA,\Xi)$ and
$(\sX\x\gA,\Xi')$ are conjugate via the function $(x,\ba)\mapsto 
(b(\ba)(x),\ba)$. Also, $C$ is a homomorphism iff there is a
continuous $\ZD$-action $\xi$ on $\sX$ such that
$(\sX\x\gA,\Xi) = (\sX,\xi)\x(\gA,\shift{})$ [i.e. $C(\fz,\ba)=\xi^\fz$,
for all  $\ba\in\gA$].  Hence, $C$ is trivial
iff $(\sX\x\gA,\Xi)$ is isomorphic to such a Cartesian product.}

 If $(\sG,\cdot)$ is
an abelian group, then the set $\sZ=\Zdyn^1(\gA,\sG)$ of all $\sG$-valued
continuous cocycles is a group under pointwise multiplication.  The
set of trivial cocycles is a subgroup $\sB=\Bdyn^1(\gA,\sG)$.  The
quotient group $\Hdyn^1(\gA,\sG):=\sZ/\sB$ is the (first dynamical) {\dfn
cohomology group} of $\gA$ (with {\dfn coefficients} in $\sG$).  If
$\sG$ is not abelian, then $\sZ$ is not a group, but we still use
$\Hdyn^1(\gA,\sG)$ to denote the set of cohomology equivalence classes
of cocycles in $\sZ$.  (Sadly, some
SFTs (e.g. dominoes) admit nontrivial 
cocycles only in nonabelian groups \cite[Thm.6.6]{Sch98}.)
The cohomology of multidimensional SFTs is closely
related \cite[Thm.4.2(b)]{Sch98} to tiling homotopy groups
(see Example \ref{X:Conway.Lagarias}).  Nontrivial cocycles
represent an algebraic obstruction to the
`hole-filling problem'.  For example, in the full shift $\AZD$, the
hole-filling problem is trivial, and indeed, $\Hdyn^1(\AZD,\bS_n)$
is trivial
\cite{Kam90,Kam92,Kam93}.  More generally,
if $\gA\subset\AZD$ has certain mixing properties, then $\Hdyn^1(\AZD,\sG)$ is
trivial \cite[Thm.3.2, Cor.3.3-3.4]{Sch95a}.

\paragraph{Cocycles along trails:}
Let $\dE:=\set{\fz\in\ZD}{\fz=(0,...,0,\pm1,0,...,0)}$.
Recall that a sequence
$\zeta=(\fz_0,\fz_1,\ldots,\fz_N)\subset\ZD$ is a {\dfn trail}
if $\fz'_n\in\dE$ for all $n\in\CC{1...N}$, where
$\fz_n':=\fz_{n}-\fz_{n-1}$.
Let $r>0$ and let $c:\dE\x\gA_{(r)}\into\sG$ be some function.
We define
\beqn
\label{trail.cocycle}
  c(\zeta,\ba) \quad:=\quad \prod_{n=1}^{N} 
c(\fz'_{n},\ba_{\dB(\fz_{n-1},r)}).
\eeqn
Suppose that, for all $\fe,\fe'\in\dE$, and $\ba\in\gA$,
\beqn
\label{trail.cocycle.eqn}
\begin{array}{rrcl}
{\bf(a)} &  c(\fe',\ba_{\dB(\fe,r)}) \cdot c(\fe,\ba_{\dB(r)}) 
& = & c(\fe,\ba_{\dB(\fe',r)}) \cdot c(\fe',\ba_{\dB(r)}). \\
\quad   
{\bf(b)} &  c(-\fe,\ba_{\dB(\fe,r)}) &=&  c(\fe,\ba_{\dB(r)})^{-1}. 
\end{array}
\eeqn
Then the value of eqn.(\ref{trail.cocycle})
depends only on $\fz_0$ and $\fz_N$, and is 
independent of the particular trail $\zeta$ from 
$\fz_0$ to $\fz_N$.  In particular,
if $\zeta$ is any {\dfn closed trail} (i.e. $\fz_N=0=\fz_0$) then
$c(\zeta,\ba)=C(0,\ba)=e_\sG$.  
For any $\ba\in\gA$
and $\fz\in\ZD$, \ we define $C(\fz,\ba):= c(\zeta,\ba)$, where $\zeta$ is any
trail from $0$ to $\fz$.  The resulting function $C:\ZD\x\gA\into\sG$ is
a continuous cocycle; we say that $C$ is a {\dfn locally determined} cocycle
with {\dfn local rule} $c$ of {\dfn radius} $r$.
If $\sG$ is discrete, then every continuous
 $\sG$-valued cocycle is locally determined in this way.
For instance, 
the cocycles in Examples \ref{X:cocycle}(c,d,e,f)
had radius $r=0$, so that $\gA_{(0)}=\sA$ and the local rule was a
function $c:\dE\x\sA\into\sG$.

\example{\label{X:trail.cocycle}(a) 
Let $C:\ZD[2]\x\Ice\into\Zahl$ be as in Example \ref{X:cocycle}(c).
Any $\bi\in\Ice$ defines a set of
directed `paths' through the plane,  each without beginning
or end.  If $\zeta$ is a trail from $\fy$ to $\fz$ in $\ZD[2]$,
then $C(\zeta,\bi)=\#\{$paths which cut across 
$\zeta$ going left$\} -  \#\{$paths which cut across 
$\zeta$ going right$\}$.   In particular, if $\zeta$ is 
the counterclockwise boundary of a region 
$\dU\subset\ZD[2]$, then $C(\zeta,\bi)=\#\{
\mbox{paths entering $\dU$}\}-\#\{\mbox{paths leaving $\dU$}\}=0$
(because every path which enters $\dU$ must leave).

(b) Let $C:\ZD[2]\x\Pth\into(\Zahlmod{2})^2$ be as in Example \ref{X:cocycle}(d).
Any $\bp\in\Pth$ defines a set of undirected paths in two colours,
say `blue' and `red'.
If $\zeta$ is a trail from $\fy$ to $\fz$ in $\ZD[2]$,
then $C(\zeta,\bp) = (b,r)\in(\Zahlmod{2})^2$,
where $b$ is the parity of blue paths crossing $\zeta$, and
$r$ is the parity of red paths.
}
If $C_1$ and $C_2$ are have local rules $c_1,c_2:\dE\x\gA_{(R)}\into\sG$,
then $C_1\approx C_2$ iff there is some {\dfn local transfer function}
$b:\gA_{(r)}\into\sG$ (for some $r\leq R-1$)
such that:
\beqn
\label{local.cohomology}
\mbox{For any $\fe\in\dE$ and $\ba\in\gA_{(R)}$},\quad
c_2(\fe,\ba) \quad=\quad b(\ba_{\dB(\fe,r)})
\cdot c_1(\fe,\ba)\cdot b(\ba_{\dB(r)})^{-1}.
\eeqn
\paragraph{Fundamental cocycles:} Fix a cocycle $C^*:\ZD\x\gA\into\sG$.
If $\psi:(\sG,\cdot)\into(\sH,\cdot)$ is any group homomorphism, then
$\psi\circ C^*$ is also a cocycle.  The cocycle $C^*$ is called {\dfn
fundamental} \cite{Sch98} if, for any group $(\sH,\cdot)$ and any
cocycle $C\in\Zdyn^1(\gA,\sH)$, there is a homomorphism $\psi:\sG\into\sH$
such that $C$ is cohomologous to $\psi\circ C^*$.  It is unknown whether
every multidimensional subshift possesses a fundamental cocycle, but
fundamental cocycles have been identified for many specific
$\ZD[2]$-shifts, including dominoes \cite[Thm.6.7]{Sch98}, rectangular
polyominoes \cite[Thm.2.7]{Ein}, L-shaped triominoes
\cite[Thm.4.8]{Ein}, three-coloured chessboards \cite[Thm.7.1]{Sch98},
lozenge tilings \cite[Thm.9.1]{Sch98}, coloured path systems
\cite[Thm.3.3]{Ein}, and certain factors of cohomologically trivial
subshifts \cite[Thm.11.1]{Sch98}.  If a fundamental cocycle exists,
then it encodes essentially the same information 
\cite[Thm.5.5]{Sch98} as the projective fundamental group of \cite{GePr}
(see \S\ref{S:projective.homotopy}).

\example{If $C:\ZD\x\gA\into\sG$ is any cocycle, then the {\dfn
$\ZD$-extension} of $C$ is the cocycle $C':\ZD\x\gA\into\ZD\x\sG$ defined
by $C'(\fz; \ba) \ := \ \lb(\fz, C(\fz,\ba)\rb)$,
for any $\fz\in\ZD$ and $\ba\in\gA$.  

(a)  The $\ZD$-extension of Example \ref{X:cocycle}(c)
is a fundamental cocycle for $\Ice$ \cite[Thm.8.1]{Sch98}.

  (b) The $\ZD$-extension of Example \ref{X:cocycle}(d)
is a fundamental cocycle for $\Pth$ \cite[Thm.3.3]{Ein}.

  (c)  The $\ZD$-extension of Example \ref{X:cocycle}(f)
is a fundamental cocycle for $\Dom$ \cite[Thm.2.7]{Ein}.
}

\paragraph{CA vs. Cocycles:}
 The following can be checked through straightforward calculation:

\Proposition{\label{block.map.cohomology}}
{
 Let $\gA\subset\AZD$ and $\gB\subset\BZD$ be subshifts.  Let
$\Phi:\gA\into\gB$ be a subshift homomorphism.

\bthmlist
\item Suppose $C:\ZD\x\gB\into\sG$ is  cocycle on $\gB$, and we define
$\Phi_*C:\ZD\x\gA\into\sG$ by $\Phi_*C(\fz,\ba)=C(\fz,\Phi(\ba))$.
Then $\Phi_*C$ is a cocycle on $\gA$.  

If  $\Phi$ has radius $R$, and $C$ is locally determined
with radius $r$, then $\Phi_*C$ is locally determined with radius $r+R$.

\item Let $C,C'\in\Zdyn^1(\gB,\sG)$.  If $C \approx C'$, then
$\Phi^*C \approx \Phi^*C'$. Thus, $\Phi$ induces a function
$\Phi_*:\Hdyn^1(\gB,\sG)\into\Hdyn^1(\gA,\sG)$.

\item If $(\sG,\cdot)$ is abelian, then 
$\Phi_*:\Hdyn^1(\gB,\sG)\into\Hdyn^1(\gA,\sG)$
is a group homomorphism.
\qed\ethmlist
}
 In particular, if $\Phi:\AZD\into\AZD$ is a cellular automaton,
and $\Phi(\gA)\subseteq\gA$, then Proposition
\ref{block.map.cohomology}(c) yields a group
endomorphism $\Phi_*:\Hdyn^1(\gA,\sG)\into \Hdyn^1(\gA,\sG)$.
[For instance, if $\fy\in\ZD$, then 
$\shift{\fy}_*:\Hdyn^1(\gA,\sG)\into \Hdyn^1(\gA,\sG)$ is the identity,
by Example \ref{shift.cohomologous}(b).] 
 If $\gA$ has an abelian fundamental cocycle, then this
cohomological endomorphism $\Phi_*$ takes a simple form.  To see this,
let $\End{\sG}$ be the set of endomorphisms of $\sG$.
If $(\sG,\cdot)$ is an abelian group, then $\End{\sG}$ is an
abelian group under pointwise multiplication.

\Proposition{\label{fund.cocycle.CA}}
{
Let $\gA\subset\AZD$ have abelian fundamental cocycle 
$C^*\in\Zdyn^1(\gA,\sG^*)$.

\bthmlist
  \item There is a group epimorphism
$\End{\sG^*}\ni \eps\mapsto C_\eps \in \Hdyn^1(\gA,\sG^*)$,
defined by $C_\eps := \eps\circ C^*$.

  \item Let $\Phi:\gA\into\gA$ be a cellular automaton.
Then there is some $\varphi\in\End{\sG^*}$ such that
$\Phi_*(C_\eps) \approx C_{\eps\circ\varphi}$ for all $\eps\in\End{\sG^*}$.
\ethmlist
}
\bthmprf
 {\bf(a)} \ For any $\eps\in\End{\sG^*}$, the function $C_\eps=\eps\circ C^*$
is a cocycle.  The map $(\eps\mapsto C_\eps)$ is a group homomorphism
because $\sG$-multiplication is commutative.  The map $(\eps\mapsto C_\eps)$
is surjective onto $\Hdyn^1(\gA,\sG^*)$ because $C^*$ is fundamental.

 {\bf(b)} \ $\Phi_* C^*$ is a cocycle, so there is some $\varphi\in\End{\sG^*}$
so that $\Phi_* C^*  \approx  \varphi\circ C^*$ 
(because $C^*$ is fundamental). For any
$\eps\in \End{\sG^*}$, note that
  $\Phi_*(C_\eps) = \eps\circ (\Phi_*C^*)$, because, 
for any $\fz\in\ZD$ and $\ba\in\gA$,  \ 
$\Phi_*(C_\eps)(\fz,\ba) =
(C_\eps)(\fz,\Phi(\ba))  =
\eps\circ C^*(\fz,\Phi(\ba)) =
\eps\circ (\Phi_*C^*)(\fz,\ba)$.
But then $\eps\circ (\Phi_*C^*) \ \homotopic[(*)] \
\eps\circ (\varphi\circ C^*) \ \eeequals{(\dagger)} \ C_{\eps\circ\varphi}$.
Here, $(*)$ is because $\Phi_* C^*  \approx  \varphi\circ C^*$ 
(say with transfer function $b$), 
so $\eps\circ (\Phi_*C^*)  \approx \eps\circ (\varphi\circ C^*)$
(with transfer function $\eps\circ b$).  $(\dagger)$ is by definition
of $C_{\eps\circ\varphi}$.
\ethmprf

\paragraph{Trail homotopy:} 
If $\dY\subset\ZD$, and $\zeta=(\fz_1\adjacent\cdots\adjacent\fz_N)$
and $\zeta'=(\fz'_1\adjacent\cdots\adjacent\fz'_{N'})$ are trails in $\dY$, 
then $\zeta'$ is an {\dfn elementary $\dY$-homotope} of $\zeta$ 
(notation: $\zeta\elhomotopic\zeta'$) 
if there is some $n\in\CC{1...N}$ such that
$\fz'_i=\fz_i$ for all $i\in\CO{1...n}$, and one of the following is true:
\bdesc
  \item[(EH1)] $N'=N$ and $\fz'_{i}=\fz_{i}$ for all $i\in\OC{n...N}$, as
follows:\hspace{-3em}
{\footnotesize $\begin{array}[b]{rcl}
\cdots\leadsto\fz'_{n-1}=\fz_{n-1} & \leadsto & \fz_n \\
 \rotatebox{-90}{$\!\!\!\leadsto$} & &  \rotatebox{-90}{$\!\!\!\leadsto$} \\
\fz'_{n} & \leadsto & \fz_{n+1}=\fz'_{n+1}\leadsto\cdots
\end{array}$}

  \item[(EH2)] $\fz_{n+1}=\fz_{n-1}$, \  $N'=N-1$, \ 
and $\fz'_{i-1}=\fz_{i}$ for all $i\in\OC{n...N}$, as follows:

\centerline{\footnotesize
$\begin{array}[c]{rcl}
 & \fz_n \\
& \rotatebox{50}{$\!\!\!\leadsto$} \quad  \rotatebox{-50}{$\!\!\!\leadsto$} \\
\begin{array}[t]{c}\cdots\leadsto \fz_{n-2} \leadsto \\ \ || \\ \cdots\leadsto \fz'_{n-2} \leadsto \end{array}\hspace{-1.5em}
  & \begin{array}[t]{c} \fz_{n-1}=\fz_{n+1}\\ 
|| \qquad  || \\
 \fz'_{n-1}=\fz'_{n} \ \end{array}  & 
\hspace{-1.5em} \begin{array}[t]{c} \leadsto \fz_{n+2} \leadsto\cdots
\\ \!\!\!|| \quad \\ \leadsto \fz'_{n+1} \leadsto\cdots \end{array} 
\end{array}$}

  \item[(EH3)] Same as {\bf(EH2)}, but with $\zeta$ and $\zeta'$ switched.
\edesc
Two trails $\zeta$ and $\zeta'$ are {\dfn homotopic in $\dY$}
(notation: $\zeta\homotopic\zeta'$) 
if there is a sequence of elementary $\dY$-homotopes
 $\zeta=\zeta_0\elhomotopic\zeta_1\elhomotopic\cdots
\elhomotopic\zeta_M=\zeta'$.   This is clearly an equivalence relation.
Assume $\dY$ is connected; then every homotopy class of $\pi_1(\dY)$
can be represented as a trail in $\dY$, and two such trails are
$\dY$-homotopic iff they belong in the same class of  $\pi_1(\dY)$.
Hence we can treat  $\pi_1(\dY)$ as a group of $\dY$-homotopy classes
of $\dY$-trails.

\subsection{Poles and Residues\label{S:pole}}

If $\ba\in\tlgA$, and $\zeta=(\fz_0,\fz_1,\ldots,\fz_N)\subset\unflawed(\ba)$ is
a closed trail in $\unflawed(\ba)$, then
we can define $c(\zeta,\ba)$ as in eqn.(\ref{trail.cocycle})
(see also \cite[p.1489]{Sch98}).   This yields a natural algebraic invariant
for range-$r$ codimension-two defects:

\Proposition{\label{cohom.defect.prop}}
{
Let $\gA\subset\AZD$ be a subshift.
Let $\ba\in\tlgA$ have a range $r$ codimension-two defect.
Let $C\in\Zdyn^1(\gA,\sG)$ be locally determined with radius $r$.
Then:
\bthmlist
\item There is a group homomorphism 
$\Residue{\ba} C\colon\pi_1[\unflawed(\ba)]\rightarrow\sG$
defined $\Residue{\ba} C(\undzet)  :=  C(\zeta,\ba)$. 
\item If $(\sG,\cdot)$ is abelian, and $C\approx C'$, then 
$\Residue{\ba} C \equiv \Residue{\ba} C'$.
\ethmlist
}

The corresponding result for projective codimension-two defects is as follows:

\Theorem{\label{cohom.defect}}
{
Let $\gA\subset\AZD$ be a subshift, and let $(\sG,\cdot)$ be a discrete group.
Let $\ba\in\tlgA$ have a 
projective codimension-two defect.  Let $\unflawed[\oo]:=\unflawed[\oo](\ba)$ and let $\omg:\CO{0,\oo}\into\unflawed[\oo]$
be a proper base ray.
Define $\residue:\Hdyn^1(\gA,\sG)\x\pi_1(\unflawed[\oo],\omg)\into\sG$ 
by $\residue(\undC,\undzet)  :=  C(\zeta,\ba)$.
\bthmlist 
  \item   If  $\residue$ is nontrivial, then $\ba$ has an essential 
codimension-two defect.

  \item If $(\sG,\cdot)$ is abelian, then  $\residue$ is a group homomorphism.
\ethmlist
}

If the homomorphism $\Residue{\ba} C$ in Proposition \ref{cohom.defect.prop}
is nontrivial, we say that $\ba$ has a {\dfn $C$-pole} (of range $r$),
and $\Residue{\ba} C$ is called the {\dfn $C$-residue} of $\ba$, by analogy
with complex analysis.  In this analogy, elements of $\gA$ are like
entire functions, elements of $\tlgA$ with codimension-two defects are
like meromorphic functions, and $C(\zeta,\ba)$ is like a contour
integral.  If the function $\residue$ in Theorem
\ref{cohom.defect} is nontrivial, we say that $\ba$ has a {\dfn
(projective) $\sG$-pole}, and $\residue$ is called the {\dfn
$\sG$-residue} of $\ba$.

\ignore{If $C_\ba$ is nontrivial for any cocycle $C$, we say that
$\ba$ has a {\dfn cohomological defect}, and the homomorphism
$\residue$ is called the {\dfn cohomological signature} of $\ba$.}
\ignore{  This terminology is inspired by Example
\ref{X:trail.cocycle}(a), where the square ice cocycle has a physical
interpretation as a kind of `flux'.  An admissible $\Ice$ tiling
corresponds to a `divergence-free vector field', and a pole
represents a `source' or `sink' of the vector field; the residue
measures the nontrivial divergence at the source/sink.  In the dual
version of square ice \cite[p.1096]{GePr}, an admissible tiling
corresponds to a `curl-free vector field', and a pole is a point of
nonzero curl; the residue measures the intensity of rotation.  If we
interpret the cocycle as a `height function' (Example
\ref{X:cocycle}(e)), then a path around a pole resembles a `spiral
staircase'; the residue measures the `steepness' of this staircase.
}

\example{\label{X:ice.pole}
(a) Let $C:\ZD[2]\x\Ice\into\Zahl$ be as in Example \ref{X:cocycle}(c),
and let $\bi\in\tlIce$
be the configuration shown in Figure \ref{fig:path.tiling}(D), having a
codimension-two defect. 
 If $\zeta$ is any simple, closed clockwise trail
 around this defect, then $C(\zeta,\bi)=8$.  
Observe that $\pi_1(\unflawed[\oo](\bi))\cong\Zahl$ 
is the cyclic group generated by $\undzet$.
For any  $\undzet^n\in \pi_1(\unflawed[\oo])$,
we have $\residue[\bi](\undC,\undzet^n)=8n$.  Thus, 
$\bi$ has a projective pole, and hence, an essential defect.

(b)  Theorem \ref{cohom.defect}(a) is false if we replace 
$\pi_1[\unflawed[\oo](\ba)]$ with $\pi_1[\unflawed(\ba)]$ for some finite
$r>0$. For example, let $C:\ZD[2]\x\Pth\into(\Zahlmod{2})^2$ be as in Example 
\ref{X:cocycle}(d), and let $\bp\in\tlPth$
be the configuration shown in Figure \ref{fig:path.tiling}(E),
having a codimension-two defect with three
components, labelled \textcircled{\footnotesize1}, \textcircled{\footnotesize2}, and \textcircled{\footnotesize3}.
For $k=1,2,3$, let
$\zeta_k$ be a simple clockwise loop going around \textcircled{\footnotesize$k$},
and not around the other two defects.  Then
$C(\zeta_1,\bp)=(1,1)$, \ $C(\zeta_2,\bp)=(0,1)$, \ and
$C(\zeta_3,\bp)=(1,0)$. Hence, each of the
defects \textcircled{\footnotesize1}, \textcircled{\footnotesize2}, and \textcircled{\footnotesize3}
individually is a nontrivial range-1 pole.
  However,
$\pi_1(\unflawed[\oo](\bp))\cong\Zahl$ is the cyclic group generated by a
simple closed curve $\zeta$ that goes around all three defects, and
$\residue[\bp](\undC,\undzet)=C(\zeta,\bp)= C(\zeta_1,\bp) +
C(\zeta_2,\bp) + C(\zeta_3,\bp)=(0,0)$.  Hence, $\residue[\bp](\undC,\blankspace)$ is
trivial.  Indeed, inspection of of Figure \ref{fig:path.tiling}(E)
shows that the defect can be removed through a local change;
hence $\residue[\bp]$ must be trivial by Theorem \ref{cohom.defect}(a).

(c) The converse to Theorem \ref{cohom.defect}(a) is
false.   Triviality of $\residue$ is necessary, 
but {\em not sufficient} to conclude that the defect in $\ba$ is removable.
For example, the codimension-two defect in the
domino tiling of Figure \ref{fig:path.tiling}(F)
is essential, but  $\residue C \equiv e$ for every cocycle $C$.
(This follows from \cite[Example 2.3]{Ein}, 
attributed to Sam Lightwood).}

The proofs of Proposition \ref{cohom.defect.prop}
and Theorem \ref{cohom.defect} depend on the following:

\Lemma{\label{cohom.defect.lemma}}
{
Let $\gA\subset\AZD$, \  $\ba\in\tlgA$, and  $C\in\Zdyn^1(\gA,\sG)$ be as in {\rm Proposition \ref{cohom.defect.prop}}.
Let $\zeta$ be a closed trail in $\unflawed:=\unflawed(\ba)$.
\bthmlist
  \item  If $\zeta$ is homotopic to $\zeta'$ in $\unflawed$, then
$C(\zeta,\ba) = C(\zeta',\ba)$.  
In particular, if $\zeta$ is nullhomotopic
in $\unflawed$, then $C(\zeta,\ba)=e_\sG$.

  \item If $C \approx C'$ via transfer function $b$,
then $C(\zeta,\ba) = b(\ba)C'(\zeta,\ba) b(\ba)^{-1}$
{\rm(assuming $\zeta$ begins and ends at $0$)}. 
Hence, if $(\sG,\cdot)$ is abelian, then 
$C(\zeta,\ba) = C'(\zeta,\ba)$.   
\ethmlist
}
\bthmprf
 It suffices to check {\bf(a)} this when $\zeta$ and $\zeta'$
differ by an elementary homotopy, which can be done by
combining eqn.(\ref{trail.cocycle.eqn}) with
{\bf(EH1)-(EH3)}.  To see {\bf(b)}, substitute eqn.(\ref{local.cohomology})
into eqn.(\ref{trail.cocycle}) to get
$C(\zeta,\ba) = b(\ba)C'(\zeta,\ba) b(\ba)^{-1}$.
\ethmprf

\vspace{-1em}

\addvspace{1em}

\bthmprf[Proof of Proposition \ref{cohom.defect.prop}:] 
{\bf(a)} \ Lemma \ref{cohom.defect.lemma}(a) implies that  $\Residue{\ba} C(\undzet)$ is well-defined,
because $C(\ba,\zeta)$ is determined
by the homotopy class $\undzet$.  If
 $\zeta=\zeta_1 \star \zeta_2$, then eqn.(\ref{trail.cocycle}) implies that
that $C(\zeta,\ba)=C(\zeta_1,\ba)\cdot C(\zeta_2,\ba)$;
hence $\Residue{\ba} C(\undzet)=\Residue{\ba} C(\undzet_1)\cdot\Residue{\ba} C(\undzet_2)$;
hence $\Residue{\ba} C$ is a homomorphism.
\ {\bf(b)} follows from Lemma \ref{cohom.defect.lemma}(b).
\ethmprf

\vspace{-1em}

\addvspace{1em}

\bthmprf[Proof of Theorem \ref{cohom.defect}:] 
{\bf(a)}   By contradiction, suppose there exist $R\geq r$ and
 $\ba'\in\gA$ with
$\ba'_{\unflawed[R]}=\ba_{\unflawed[R]}$, where $\unflawed[R]:=\unflawed[R](\ba)$.
Let $\undzet\in\pi_1(\unflawed[\oo],\omg)$;
then $\undzet$ has a representative
trail $\zeta'$ in $\unflawed[R]$. 
Thus $\residue C(\undzet) \ \eeequals{(*)} \
\residue C(\undzet') :=  C(\zeta',\ba)
\ \eeequals{(\dagger)} \
C(\zeta',\ba') \ \eeequals{(\ddagger)} \ e_\sG$.
Here, $(*)$ is because $\zeta \homoto \zeta'$; \ 
$(\dagger)$ is because $\ba'_{\unflawed[R]}=\ba_{\unflawed[R]}$; \  and
$(\ddagger)$ is by Lemma \ref{cohom.defect.lemma}(a),
 because $\zeta'$ is nullhomotopic in  $\unflawed[R](\ba')=\ZD$.

{\bf(b)}  Fix $\undzet\in\pi_1(\unflawed[\oo],\omg)$.  
We claim  the function
$\Hdyn^1(\gA,\sG) \ni \undC \mapsto C(\zeta,\ba)\in\sG$
is a homomorphism.  First, note that, for any
$\undC\in \Hdyn^1(\gA,\sG)$, with any radius $R>0$, we can find
some representative of $\undzet$ in $\pi_1(\unflawed[R],\omg)$
[because $\undzet\in\pi_1(\unflawed[\oo],\omg)$].
Hence $C(\zeta,\ba)$ is always well-defined.
Furthermore, the value of $C(\zeta,\ba)$
depends only on the cohomology class of $C$, by
Lemma \ref{cohom.defect.lemma}(b).
Now, let $\undC_1,\undC_2\in\Hdyn^1(\gA,\sG)$, and
$\undC:=\undC_1\cdot \undC_2$; it follows from 
eqn.(\ref{trail.cocycle}) that $C(\zeta,\ba) =
C_1(\zeta,\ba) \cdot C_2(\zeta,\ba)$.

 Now fix $\undC\in\Hdyn^1(\gA,\sG)$. 
We claim  the function
$\Residue[\oo]{\ba}{C}:\pi_1(\unflawed[\oo],\omg) \ni \undzet \mapsto C(\zeta,\ba)\in\sG$
is a homomorphism.
Now, $C$  is locally determined (say, with radius $r$), because $\sG$
is discrete.  Thus, Proposition \ref{cohom.defect.prop} yields a homomorphism
$\Residue[R]{\ba} C:\pi_1(\unflawed[R],\omg)\ni\undzet\mapsto C(\zeta,\ba)\in\sG$, for all $R\geq r$.  This yields a commuting cone of homomorphisms:
\[
\psfrag{A1}[][]{\footnotesize$\pi_1(\unflawed[r],\omg)$}
\psfrag{A2}[][]{\footnotesize$\pi_1(\unflawed[r+1],\omg)$}
\psfrag{A3}[][]{\footnotesize$\pi_1(\unflawed[r+2],\omg)$}
\psfrag{Aoo}[][]{\footnotesize$\pi_1(\unflawed[\oo],\omg)$}
\psfrag{f1}[][]{\scriptsize$\iota_r^*$}
\psfrag{f2}[][]{\scriptsize$\iota_{r+1}^*$}
\psfrag{f3}[][]{\scriptsize$\iota_{r+2}^*$}
\psfrag{g1}[][]{\scriptsize$\Residue{\ba} C$}
\psfrag{g2}[][]{\scriptsize$\Residue[r+1]{\ba} C$}
\psfrag{g3}[][]{\scriptsize$\Residue[r+2]{\ba} C$}
\psfrag{goo}[][]{\scriptsize$\Residue[\oo]{\ba} C$}
\psfrag{X}[][]{\footnotesize$\sG$}
\includegraphics[scale=0.8]{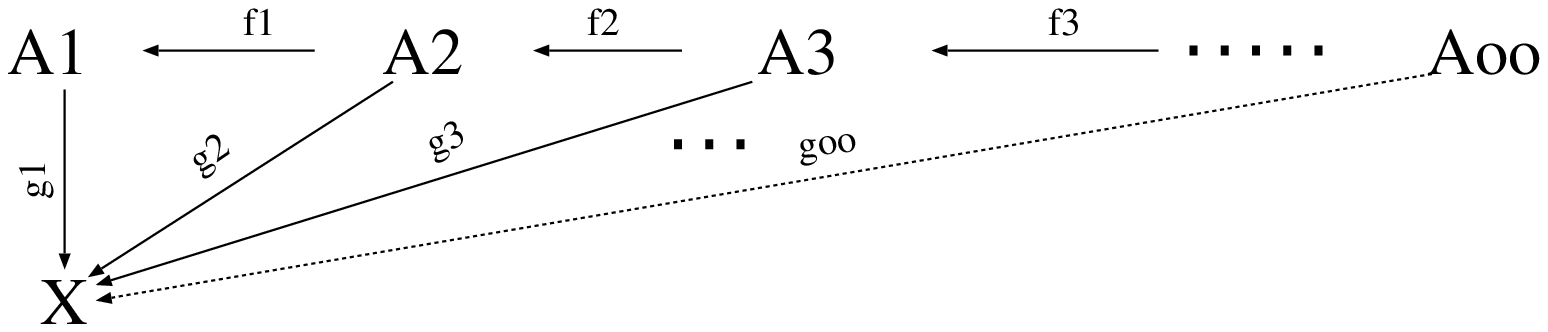}
\]
Each triangle in this cone commutes because of 
Lemma \ref{cohom.defect.lemma}(a).
The cone converges to 
$\Residue[\oo]{\ba}{C}:\pi_1(\unflawed[\oo],\omg)\into\sG$,
so $\Residue[\oo]{\ba}{C}$ is also a homomorphism.
\ethmprf

When are poles persistent under cellular automata?
Recall from Proposition \ref{block.map.cohomology} that any
CA $\Phi:\gA\into\gA$ defines an endomorphism 
$\Phi_*:\Hdyn^1(\gA,\sG)\into \Hdyn^1(\gA,\sG)$.

\Proposition{\label{persistent.pole}}
{
  Let $\Phi:\AZD\into\AZD$ be a CA {\rm(with radius $R>0$)} and
let $\gA\subset\AZD$ be a $\Phi$-invariant subshift.
Let $\ba\in\tlgA$ and let $\bb=\Phi(\ba)$.
\bthmlist
  \item If  $C\in\Zdyn^1(\gA,\sG)$ is locally determined with radius $r>0$,
then  $\Residue[R+r]{\bb}C \ \equiv \ \Residue[R+r]{\ba} (\Phi_*C)$.

  \item If  $\Phi_*:\Hdyn^1(\gA,\sG)\into \Hdyn^1(\gA,\sG)$ is surjective,
then every $\sG$-pole is $\Phi$-persistent.
\ethmlist
Suppose $\gA$ has abelian fundamental cocycle 
$C^*\in\Zdyn^1(\gA,\sG^*)$.
Define
$\residue^*:\pi_1(\unflawed[\oo](\ba))\into\sG^*$ by
$\residue^*(\undzet)=C^*(\zeta,\ba)$ for any $\undzet\in\pi_1(\unflawed[\oo](\ba))$.
\bthmlist
\setcounter{enumi}{2}
  \item If $(\sH,\cdot)$ is any abelian group, and $\ba$ has an $\sH$-pole,
then $\residue^*$ is nontrivial.

  \item Let $\Phi$ and $\varphi$ be as in {\rm Proposition
\ref{fund.cocycle.CA}(b)}.
If $\bb=\Phi(\ba)$, then $\residue[\bb]^* = \varphi\circ\residue^*$.

Hence, if $\varphi^n\circ\residue^*$ is nontrivial for all $n\in\Natur$,
then the defect of $\ba$ is $\Phi$-persistent.

  \item If $\varphi$ is a monomorphism, 
then every pole is $\Phi$-persistent.
\ethmlist
}
\bthmprf {\bf(a)}
If $\zeta$ is any trail in $\unflawed[R+r](\ba)$, then
$C(\bb,\zeta) \ = \ \Phi_*C (\ba,\zeta)$.
{\bf(b)} follows.

{\bf(c)}  \ Suppose there exists $C\in\Zdyn^1(\gA,\sH)$, and 
$\undzet\in\pi_1(\unflawed[\oo](\ba))$ such that 
$\residue(\undC,\undzet)\neq e_\sH$.  There is a homomorphism $\psi:\sG^*\into\sH$
such that $\psi\circ C^*$ is cohomologous to $C$
(because $C^*$ is fundamental). Thus
\[
\psi \lb[C^*(\zeta,\ba)\rb]
\quad \eeequals{(\dagger)} \quad
 (\psi\circ C^*)(\zeta,\ba) 
\quad \eeequals{(*)} \quad   C(\zeta,\ba)
\ \ =: \ \ \residue(\undC,\undzet) 
\ \  \neq\ \  e_\sH.
\]
Here, $(\dagger)$ is by applying homomorphism $\psi$ to eqn.(\ref{trail.cocycle}), and  $(*)$ is by Lemma \ref{cohom.defect.lemma}(b)
(because $\psi\circ C^*\approx C$).
Thus $\residue^*(\undzet):=C^*(\zeta,\ba)$ is nontrivial
in $\sG^*$; hence $\residue^*$ is nontrivial. 

{\bf(d)} Let $\bb=\Phi(\ba)$.  Then for any $\undzet\in\pi_1(\unflawed[\oo](\ba))$,
$\residue[\bb]^*(\undzet) = C^*(\zeta,\bb)=C^*(\zeta,\Phi(\ba))
= \Phi_*C^*(\zeta,\ba) \ \eeequals{(\dagger)}  \ \varphi\circ C^*(\zeta,\ba)
= \varphi\circ\residue^*(\undzet)$, where $(\dagger)$ is by 
Proposition \ref{fund.cocycle.CA}(b)
and Lemma \ref{cohom.defect.lemma}(b). \quad Then {\bf(e)} follows from {\bf(d)}.
\ethmprf

\ignore{  Verifying the $\Phi$-persistence of a codimension-two defect in 
$\ba\in\tlgA$ could involve a lengthy process of: (1) finding
a group $(\sG,\cdot)$ such that the $\sG$-residue $\residue$ is
nontrivial, (2) computing the induced endomorphism
$\Phi^*:\Hdyn^1(\gA,\sG)\into \Hdyn^1(\gA,\sG)$, and then (3) checking the
interaction of $\Phi^*$ with $\residue$.  Sometimes there is a nice
shortcut. }

\ignore{  If $C^*\in\Zdyn^1(\gA,\sG^*)$ is a fundamental cocycle for $\gA$,
and $\ba$ has an essential codimension-two defect, then the
homomorphism $\residue^*$ introduced in 
Proposition \ref{}(c) is called the
 {\dfn fundamental residue} of $\ba$.}

\subsection{Gaps and Tilt \label{S:gaps}}
The domain boundary in Examples \ref{X:codim.one}(c,d)
are not detected by the spectral invariants
of \cite{PivatoDefect1}. 
However, they can be detected using cohomology.
Let $C:\ZD\x\gA\into\sG$ be a locally determined cocycle with radius $r>0$ and
local rule
$c:\dE\x\gA_{(r)}\into\sG$.  For any $\ba\in\gA$,  we define 
$\cgap:\ZD\x\ZD\into\Natur$ by
\beqn
\label{tilt.defn}
  \cgap(\fy,\fz) \ \ :=\ \ 
 C(\fy,\ba)\cdot C(\fz,\ba)^{-1} \ \ = \ \  c(\zeta,\ba),
\eeqn
where $\zeta$ is any trail from $\fz$ to $\fy$, and the
expression $c(\zeta,\ba)$ is interpreted as in eqn.(\ref{trail.cocycle}).
For example, if $C:\ZD\x\gA\into\Zahl$ is a height function
[Example \ref{X:cocycle}(e)] and $\ba\in\gA$, then 
$C$ assigns a `height'  $h(\fz):=C(\fz,\ba)$ to every point
in $\ZD$, so $\cgap(\fy,\fz)=h(\fy)-h(\fz)$ is
the `altitude difference' between $\fy$ and $\fz$.

  Now suppose $\ba\in\tlgA$ has a range $r$ domain boundary and let
$\dY:=\unflawed(\ba)$ have projective connected components
$\dY_1,\ldots,\dY_N$.  Assume that $\ba$ has no codimension-two defects
--hence $\pi_1(\dY_n)=0$ for all $n\in\CC{1...N}$.  Thus, if
$\fy_1,\fy_2\in\dY$ are in the same projective component of $\dY$, then
we can define $\cgap(\fy_1,\fy_2)$ by the right-hand expression in
(\ref{tilt.defn}) (this is path-independent because $\pi_1(\dY_n)=0$).
However, if $\fy_1,\fy_2$ are in {\em different} connected components
of $\dY$, then $\cgap(\fy_1,\fy_2)$ is not well-defined by
eqn.(\ref{tilt.defn}), because there is no trail in $\dY$ connecting
$\fy_1$ to $\fy_2$.  Instead, we will define $\cgap$ up to a constant
as follows: first, for each $n\in\CC{1...N}$, choose a reference point
$\fy_n^*\in\dY_n$ and decree that $\cgap(\fy^*_n,\fy^*_m):=e_\sG$ for
all $n,m\in\CC{1...N}$.  Then for any $\fy_n\in\dY_n$ and
$\fy_m\in\dY_m$, define
\beqn
\label{tilt.defn2}
\cgap(\fy_n,\fy_m)\quad:=\quad\cgap(\fy_n,\fy^*_n)\cdot\cgap(\fy^*_m,\fy_m).
\eeqn
 Let $\sC:=c[\dE\x\gA_{(r)}]\subset\sG$; then $\sC$
is a finite subset of $\sG$, and for any $\fz\in\ZD$ and $\ba\in\gA$,
$C(\fz,\ba)$ is an element of the subgroup generated by $\sC$.
Hence we can assume without loss of generality that $\sC$ generates
$\sG$.
A function $|\blankspace|:\sG\into\CC{0,\oo}$ is a {\dfn pseudonorm} on $\sG$
if:
\beqn
\label{pseudonorm.dfn}
\mbox{For all $g,h\in\sG$,} \qquad
{\bf(a)} \ \  |g\cdot h| \ =  \ |h\cdot g|,
\qquad\And\quad
{\bf(b)} \  \ |g\cdot h|  \ \leq \ |g| + |h|.
\eeqn
(Hence $|\blankspace|$ is constant on each conjugacy class of
$\sG$).    

\example{
(a) If $\sG=\Zahl$, let $|\blankspace|$ be the Euclidean norm.

(b) If $\sG$ is abelian, then $\sG$ is finitely generated (by $\sC$), so
$\sG\cong \Zahl^R \dirsum \Zahlmod{n_1}\dirsum\cdots\dirsum
\Zahlmod{n_K}$ for some $R,n_1,\ldots,n_K\in\Natur$; see
\cite[\S5.2]{DummitFoote} or \cite[\S I.10]{Lang}.  For any
$(z_1,\ldots,z_R)\in\Zahl^R$ and
 $(y_1,\ldots,y_K)\in\Dirsum_{k=1}^K \Zahlmod{n_k}$,
define $|(z_1,\ldots,z_R; y_1,\ldots,y_K) 
:= |z_1|+\cdots+|z_R| + |y_1|_{n_1} + \cdots |y_k|_{n_K}$, where
$|y|_{n_k} :=\min\{|y|,|n_k-y|\}$ for all $k\in\CC{1...K}$.

(c) If $\sG$ is nonabelian, then let $q:\sG\into\tlsG$ be the abelianization
map.  Then $\tlsG$ is finitely generated (by $\tlsC:=q[\sC]$),
so let $|\blankspace|_*:\tlsG\into\Natur$ be as in (b).  Then define
$|g|=|q(g)|_*$ for any $g\in\sG$.
}

\Remarks{A pseudonorm on $\sG$ is equivalent to a  pseudometric
$d:\sG\x\sG\into\CC{0,\oo}$ which is {\dfn bilaterally invariant}
i.e. $d(fg,fh)= d(g,h)=d(gf,hf)$ for all $f,g,h\in\sG$.  
Here $|g|=d(g,e)$ and $d(g,h)=|g h^{-1}|$.  We do {\em not} require that
$d$ be compatible with the topology of $\sG$ (although 
this can always be arranged if $\sG$ is
unimodular; i.e. the left- and right- Haar measures are the same).
However, the following theory is trivial unless  $|\blankspace|$
is unbounded (so if $\sG$ is compact then $d$ {\em shouldn't} be
topologically compatible with $\sG$).}

   We allow  that $|g|=0$ or $g=\oo$ for some $g\neq e$.
However, we require that (i) $\forall \ c\in\sC, \ |c|\leq 1$ and
(ii) $\exists \ c\in\sC$ with $|c|>0$.
(This can always be obtained through renormalization, if $|\blankspace|$ is
nontrivial.)  It follows that $\cgap$ satisfies a Lipschitz-type condition:

\Lemma{\label{gentle.slope}}
{\bthmlist
\item 
 If $\ba\in\gA$, and $\fy,\fz\in\ZD$, then $|\cgap(\fy,\fz)|\leq |\fy-\fz|$.

\item If $\ba\in\tlgA$, and $\fy,\fz$ are in the same connected component
of $\unflawed(\ba)$, then 
 $|\cgap(\fy,\fz)|\leq d_{r,\ba}(\fy,\fz)$,
where $d_{r,\ba}(\fy,\fz)$ is the length of the shortest trail from
$\fy$ to $\fz$ in $\unflawed(\ba)$.\qed
\ethmlist}

Suppose $\ba\in\tlgA$ has a range $r$ domain boundary and let
$\dY_1,\ldots,\dY_N$ be as above.
The {\dfn tilt} of $\dY_n$ relative to $\dY_m$ is then defined:
\beqn
\label{tilt.defn3}
  \tilt(\dY_n,\dY_m)\quad:=\quad \sup_{\fy_n\in\dY_n, \ 
\fy_n\in\dY_m} \ \frac{|\cgap(\fy_n,\fy_m)|}{|\fy_n-\fy_m|}.
\eeqn
If $\tilt(\dY_n,\dY_m)=\oo$, then we say the domain boundary is a {\dfn
$C$-gap}.

\example{\label{X:ice.gap}
(a) Let $C:\ZD[2]\x\Ice\into\Zahl$ be as in Example \ref{X:cocycle}(c),
and let $\bi\in\tlIce$ be
the domain boundary configuration shown in Figure \ref{fig:codim}(A).
Suppose for simplicity that the domain boundary straddles the $x$ axis.
Let $\dX$ and $\dY$ be the north and south connected components,
respectively.  Let $\fx^* := (0,1)\in\dX$ and  $\fy^* := (0,-1)\in\dY$,
and for all $n\in\Natur$, let 
$\fx_n := (n,1)\in\dX$ and  $\fy_n := (n,-1)\in\dY$,
as shown in Figure \ref{fig:codim}(A).
Then $\cgap[\bi](\fx_n,\fx^*)=n$,
while $\cgap[\bi](\fy_n,\fy^*)=-n$;  hence $\cgap[\bi](\fx_n,\fy_n)=2n$,
However, $|\fx_n-\fy_n|=2$ for all $n$; hence $\tilt[\bi](\dX,\dY)=\oo$,
so this defect is a gap.

(b) Let $C:\ZD[2]\x\Dom\into\sG:=\Zahlmod{2}*\Zahlmod{2}$ be the cocycle
from Example \ref{X:cocycle}(f).
Unfortunately, $\sG$ does not
admit any nontrivial pseudonorms (because every nonidentity element belongs
to the same conjugacy class).  However, if $\sZ\subset\sG$ is the cyclic
subgroup generated by $vh$, then $(\sZ,\cdot)\cong (\Zahl,+)$,
and for any $\bd\in\Dom$ and $2\fz\in 2\ZD[2]$, \  $C(2\fz,\bd)\in\sZ$.
Let $\sD_{2}\subset\sD^{2\x2}$ be the alphabet of $\Dom$-admissible
$2\x 2$ blocks, and let $\gD_{2}\subset \sD_{2}^{\ZD[2]}$
 be the `recoding' of $\Dom$ in this alphabet. 
Then $2\ZD[2]$ acts on $\gD_{2}$ by
shifts in the obvious way, and $C$ yields a cocycle
$C':2\ZD[2] \x \gD_{2} \into \sZ\cong \Zahl$.

  Let $\bd\in\tlDom$ be the
domain boundary configuration in Figure \ref{fig:codim}(B)
and let $\bd_{2}$ be its recoding as an element of $\tlgD_2$.
Let $\fx^*:=(0,2)\in\dX\intsct(2\ZD[2])$ and 
$\fy^*:=(0,-2)\in\dY\intsct(2\ZD[2])$,
and for all $n\in\Natur$, let
$\fx_n:=(2n,2)\in\dX\intsct(2\ZD[2])$ and 
$\fy_n:=(2n,-2)\in\dY\intsct(2\ZD[2])$, as shown in Figure \ref{fig:codim}(B).
Then $C'_{\bd_2}(\fx_n,\fx^*) \ = \ (vhvh)^n $,
whereas $C'_{\bd_2}(\fy_n,\fy^*)  = h^{2n} = e_\sG =  (vh)^0$.
Hence $C'_{\bd_2}(\fx_n,\fy_n) = (vhvh)^n \cong 2n\in\Zahl$.
However, $|\fx_n-\fy_n|=2$ for all $n$; 
hence $\tilt(\dX,\dY)=\oo$,
so this defect is a gap.

(c)  Let $\bd'\in\tlDom$ be the
domain boundary configuration in Figure \ref{fig:codim}(C),
and let $\bd'_{2}$ be its recoding as an element of $\tlgD_2$, as in (b).
Let $\dX$ and $\dY$ denote the eastern and western domains 
(assume the boundary straddles the $y$ axis).
Let $\fx^*:=(-2,0)\in\dX\intsct(2\ZD[2])$ and 
$\fy^*:=(2,0)\in\dY\intsct(2\ZD[2])$,
and for all $n\in\Natur$, let
$\fx_n:=(-2,2n)\in\dX\intsct(2\ZD[2])$ and 
$\fy_n:=(2,2n)\in\dY\intsct(2\ZD[2])$,  as shown in Figure \ref{fig:codim}(C).
Then $C'_{\bd'_2}(\fx_n,\fx^*) \ = \ (hvhv)^n $,
whereas $C'_{\bd'_2}(\fy_n,\fy^*) \ = \ (vhvh)^n$.
Hence $C'_{\bd'_2}(\fx_n,\fy_n) = (hvhv)^n (vhvh)^{-n} = (vhvh)^{-2n}
\cong -4n \in \Zahl$.
However, $|\fx_n-\fy_n|=2$ for all $n$;  
hence $\tilt(\dX,\dY)=\oo$, so this defect is a gap.

(d) If $C$ is a cocycle into a finite group,
 then there can be no $C$-gaps.
For example, the cocycle $C:\ZD\x\Pth\into(\Zahlmod{2})^2$ of Example 
\ref{X:cocycle}(d) admits no gaps.
}
If $(\sG,\cdot)$ is a group, then
a {\dfn $\sG$-gap} is a $C$-gap for some $\undC\in\Hdyn^1(\gA,\sG)$.
We say the gap is {\dfn sharp} if, for all $R\geq r \geq 0$, there
is some constant $K=K(R,r)\in\Natur$ such that, for any $\fy\in\unflawed(\ba)$,
there exists some $\fx\in\unflawed[R](\ba)$ which is trail-connected
to $\fy$ in $\unflawed(\ba)$,  with $d_{r,\ba}(\fx,\fy) \leq K$.
Heuristically, this means that the defect field $\energy$ does not have
arbitrarily large `flat' areas, and that the gap does not ramify into
lots of `tributaries'.  For example, if $\gA$ is a subshift of finite
type, and the defect set $\dD(\ba)$ is confined to a thickened hyperplane
[as in Examples \ref{X:ice.gap}(a,b,c)], then the gap is sharp.

\Theorem{\label{persistent.gap}}
{
\bthmlist
  \item Sharp $\sG$-Gaps are essential defects.

  \item If $\Phi\colon\AZD\into\AZD$ is a CA with $\Phi(\gA)\subseteq\gA$,
and $\Phi_*:\Hdyn^1(\gA,\sG)\into\Hdyn^1(\gA,\sG)$
is surjective, then any $\sG$-gap is $\Phi$-persistent.
\ethmlist
}

To prove Theorem \ref{persistent.gap}, we use:

\Lemma{\label{persistent.gap.lemma}}
{
\bthmlist
  \item The existence of a gap does not depend on the choice
of $\{\fy_1^*,\ldots,\fy_N^*\}$.

  \item Gaps depend only on cohomology classes.
If $C\approx C'$, then any $C$-gap is also a $C'$-gap.
\ethmlist
}
\bthmprf
 {\bf(a)} Suppose we defined $\cgap^\dagger$ as in eqn.(\ref{tilt.defn2}), but
using a different set $\{\fy_1^\dagger,\ldots,\fy_N^\dagger\}$. 
For all $n\in\CC{1...N}$,
let $c_n:=\cgap(\fy_n^*,\fy_n^\dagger)$
(this well-defined by eqn.(\ref{tilt.defn})
because they are in the same connected component  $\dY_n$).
If $\fy_n\in\dY_n$ and $\fy_m\in\dY_m$, 
then $\cgap(\fy_n,\fy_m) \leq \cgap^\dagger(\fy_n,\fy_m) + c_n + c_m$,
because
\beq
\cgap(\fy_n,\fy_m) &\eeequals{(\diamond)}&
\cgap(\fy_n,\fy^*_n)\cdot\cgap(\fy^*_m,\fy_m)
\\&\eeequals{(\dagger)}&
\cgap(\fy_n,\fy^\dagger_n)\cdot  \cgap(\fy_n^\dagger,\fy_n^*) 
 \cdot \cgap(\fy^*_m,\fy^\dagger_m) \cdot \cgap(\fy^\dagger_m,\fy_m),\\
\mbox{so} \ \ 
\lb|\cgap(\fy_n,\fy_m)\rb| &=&
\lb|\cgap(\fy_n,\fy^\dagger_n)\cdot  \cgap(\fy_n^\dagger,\fy_n^*) 
 \cdot \cgap(\fy^*_m,\fy^\dagger_m) \cdot \cgap(\fy^\dagger_m,\fy_m)\rb|
\\&\eeequals{(\ddagger)}&
\lb| \cgap(\fy_n^\dagger,\fy_n^*) \cdot \cgap(\fy_n,\fy^\dagger_n)
  \cdot \cgap(\fy^\dagger_m,\fy_m)\cdot \cgap(\fy^*_m,\fy^\dagger_m)\rb|
\\&\leeeq{(*)}&
\lb| \cgap(\fy_n^\dagger,\fy_n^*)\rb| + \lb| \cgap(\fy_n,\fy^\dagger_n)
  \cdot \cgap(\fy^\dagger_m,\fy_m)\rb| + \lb|\cgap(\fy^*_m,\fy^\dagger_m)\rb|
\\&\eeequals{(\diamond)}&
c_n + \lb|\cgap^\dagger(\fy_n,\fy_m)\rb| + c_m.
\eeq
Here, $(\diamond)$ is by eqn.(\ref{tilt.defn2}),
$(\dagger)$ is by eqn.(\ref{tilt.defn}),
$(\ddagger)$ is by eqn.(\ref{pseudonorm.dfn}a), and
$(*)$ is by eqn.(\ref{pseudonorm.dfn}b).

Likewise, $\cgap^\dagger(\fy_n,\fy_m) \leq \cgap(\fy_n,\fy_m) + c_n + c_m$
(by symmetric reasoning).  Thus,
\[
\frac{\cgap^\dagger(\fy_n,\fy_m)-c_n-c_m}{|\fy_n-\fy_m|}
\quad\leq\quad
\frac{\cgap(\fy_n,\fy_m)}{|\fy_n-\fy_m|}
\quad\leq\quad
\frac{\cgap^\dagger(\fy_n,\fy_m)+c_n+c_m}{|\fy_n-\fy_m|}.
\]
Substitute into eqn.(\ref{tilt.defn3}) to see that
$\tilt(\dY_n,\dY_m) = \oo$ if and only if
$\Tilt^\dagger(\dY_n,\dY_m)=\oo$.

{\bf(b)} \ If $C'\approx C$, then there is a local transfer
function $b:\gA_{(r)}\into\sG$ satisfying eqn.(\ref{local.cohomology}).
Thus, for any $\fy_n\in\dY_n$ and $\fy_m\in\dY_m$,
$C'_\ba(\fy_n,\fy^*_n) = b(\ba_{\dB(\fy_n,r)})\cdot C_\ba(\fy_n,\fy^*_n)
\cdot b(\ba_{\dB(\fy_n^*,r)})^{-1}$.
Now, $b[\gA_{(r)}]$ is finite because $\gA_{(r)}$ is finite.
Thus $B:=
\max\set{|b(\ba)|}{\ba\in\gA_{(r)}}\union\set{|b(\ba)^{-1}|}{\ba\in\gA_{(r)}}$
is finite.  Furthermore,
part (a) says that we can assume without loss of generality
that $\fy_n^*$ and $\fy_m^*$ are chosen such that
$b(\ba_{\dB(\fy_n^*,r)})=b(\ba_{\dB(\fy_m^*,r)})$.  Thus,
for any  $\fy_n\in\dY_n$ and $\fy_m\in\dY_m$,
\beq
\lefteqn{C'_\ba(\fy_n,\fy_m) \quad\eeequals{(\diamond)}\quad 
C'_\ba(\fy_n,\fy^*_n)  \cdot C'_\ba(\fy^*_m,\fy_m)}
\quad\\&=&
b(\ba_{\dB(\fy_n,r)})\cdot C_\ba(\fy_n,\fy^*_n)
\cdot b(\ba_{\dB(\fy_n^*,r)})^{-1}\cdot
b(\ba_{\dB(\fy^*_m,r)})\cdot C_\ba(\fy^*_m,\fy^m)
\cdot b(\ba_{\dB(\fy_m,r)})^{-1}
\\&=&
b(\ba_{\dB(\fy_n,r)})\cdot C_\ba(\fy_n,\fy^*_n)
\cdot C_\ba(\fy^*_m,\fy^m)
\cdot b(\ba_{\dB(\fy_m,r)})^{-1}
\  \eeequals{(\diamond)} \ 
b(\ba_{\dB(\fy_n,r)})\cdot C_\ba(\fy_n,\fy_m)
\cdot b(\ba_{\dB(\fy_m,r)})^{-1}\\
\lefteqn{\hspace{-1em}\mbox{Thus,} \ \ 
|C'_\ba(\fy_n,\fy_m)|
\quad= \quad
\lb|b(\ba_{\dB(\fy_n,r)})\cdot C_\ba(\fy_n,\fy_m)
\cdot b(\ba_{\dB(\fy_m,r)})^{-1}\rb|}
\\ & \leeeq{(*)}&
\lb|b(\ba_{\dB(\fy_n,r)})\rb|+ \lb|C_\ba(\fy_n,\fy_m)\rb|
+\lb| b(\ba_{\dB(\fy_m,r)})^{-1}\rb|
\quad\leq\quad
\lb|C_\ba(\fy_n,\fy_m)\rb| + 2B.
\eeq
Here, $(\diamond)$ is by eqn.(\ref{tilt.defn2}) and
$(*)$ is by eqn.(\ref{pseudonorm.dfn}b).
By symmetric reasoning, we can show that
$|C'_\ba(\fy_n,\fy_m)| \ \leq \ \lb|C_\ba(\fy_n,\fy_m)\rb| + 2B$.
Hence, just as in part (a), we conclude that
$\tilt(\dY_n,\dY_m) = \oo$ if and only if
$\Tilt^{C'}(\dY_n,\dY_m)=\oo$.
\ethmprf

\bthmprf[Proof of Theorem \ref{persistent.gap}:]
{\bf(a)} \  Let $\dY:=\unflawed(\ba)$ have a sharp gap, but
suppose (by contradiction) that the defect in $\ba$ is removable.
Thus, there is some $R\geq r$ and  $\bb\in\gA$
such that  $\bb_{\dX}=\ba_{\dX}$, where $\dX:= \unflawed[R](\ba)$.
By Lemma \ref{persistent.gap.lemma}{(a)}, we can 
assume without loss of generality
that $\{\fy_1^*,\ldots,\fy_N^*\}$ are in $\dX$, and 
were chosen so that
$\cgap[\bb](\fy^*_n,\fy^*_m)=0$ for all $n,m\in\CC{1...M}$.
It follows that $\cgap(\fx_1,\fx_2)=\cgap[\bb](\fx_1,\fx_2)$ 
for all $\fx_1,\fx_2\in\dX$.  Let $K=K(R,r)$ be the constant
arising from the sharpness of the gap.

\Claim{For all $\fy_1,\fy_2\in\dY$,\quad $\cgap(\fy_1,\fy_2)\leq
|\fy_1-\fy_2|+4K$.}
\bclaimprf
Suppose $\fy_1$ (resp. $\fy_2$) is in projective component $\dY_1$
(resp. $\dY_2$) of $\dY$.
For $n=1,2$, there exists
$\fx_n\in\dX\intsct\dY_n$ such that $d_{r,\ba}(\fx_n,\fy_n)\leq K$
(by definition of `sharpness'). Then
\beq
\lb|\cgap(\fy_1,\fy_2)\rb|& \leeeq{(\dagger)}&
\lb| \cgap(\fy_1,\fx_1)\rb| + \lb|\cgap(\fx_1,\fx_2)\rb|
+ \lb|\cgap(\fx_2,\fy_2)\rb| \quad \leeeq{(*)} \quad
  K + \lb|\cgap[\bb](\fx_1,\fx_2)\rb| + K 
\\ &\leeeq{(\diamond)}& K + |\fx_1 -\fx_2| + K
\ \ \leeeq{(\bigtriangleup)}\  \
 K + K + |\fy_1 -\fy_2| + K + K.
\eeq
$(\dagger)$ is eqn.(\ref{pseudonorm.dfn}b),
$(*)$ is  Lemma \ref{gentle.slope}(b) applied to $\fx_1,\fy_1\in\dY_1$
and $\fx_2,\fy_2\in\dY_2$ in $\ba$,
and $(\diamond)$ is  Lemma \ref{gentle.slope}(a) applied to $\fx_n,\fx_m$
in $\bb$.  $(\bigtriangleup)$ is the triangle inequality.
\eclaimprf
 For any $n,m\in\CC{1...M}$,
Claim 1 and eqn.(\ref{tilt.defn3}) 
imply that $\tilt(\dY_n,\dY_m)\leq 1$,
which means that $\ba$ has no gap between $\dY_n$ and $\dY_m$.
Contradiction.

{\bf(b)} \ Suppose $\ba\in\tlgA$ has a $C$-gap for some
 $C\in\Zdyn^1(\gA,\sG)$, and let $\ba':=\Phi(\ba)$.
  There exists $C'\in\Zdyn^1(\gA,\sG)$
with $C_1:=\Phi_*(C')\approx C$ 
(because $\Phi_*:\Hdyn^1(\gA,\sG)\into\Hdyn^1(\gA,\sG)$ is surjective).
Thus, Lemma \ref{persistent.gap.lemma}{(b)} says that $\ba$ also has a $C_1$-gap.  But
for any $\fz\in\ZD$, \ 
$C_1(\fz,\ba) = C'(\fz,\Phi(\ba))= C'(\fz,\ba')$.  Thus,
$\ba$ has a $C_1$-gap iff $\ba'$ has a $C'$-gap.
\ethmprf

\section{Homotopy Defects for Subshifts of Finite Type \label{S:homotopy}}

We will introduce homotopy/homology groups for Wang tile systems and
subshifts of finite type, generalizing the constructions of
\cite{Cola,GePr}.  Nontrivial elements of these groups represent
codimension-$(d+1)$ `obstructions'  to the hole-filling problem, and
can be used to characterize codimension-$(d+1)$
defects.  This section's main results are Theorem
\ref{homotopy.comm.square} and Corollary \ref{persistent.homotopy.defect}.

\subsection{Canonical cell complex of $\RD$: \label{S:canonical.cell.cplx}}
For any $\fz\in\ZD$, let $\CornerTile{\fz}{\ } := \fz+\CC{0,1}^D$ be
the unit cube in $\RD$ with its minimal-coordinate corner at $\fz$.
For all $d\in\CC{1...D}$, let $\fe_d:=(0,...,0,1,0,...,0)$ be the $d$th
unit vector.  Suppose $\fz=(z_1,\ldots,z_D)$; for all $d\in\CC{1...D}$,
let $\partial_d^{-}\ \CornerTile{\fz}{\ }:=\set{\bx\in\CornerTile{\fz}{\ }}{x_d=z_d}$ and $\partial_d^{+}\
\CornerTile{\fz}{\ }:=\set{\bx\in\CornerTile{\fz}{\ }}{x_d=z_d+1}$ be
the two  `faces' of $\CornerTile{\fz}{\ }$ of codimension one
which are
orthogonal to $\fe_d$.  For any $k\in\CC{2...D}$, we define
codimension-$k$  faces by intersecting $k$ of these codimension-one
faces (e.g. if $\bd:=\{d_1,\ldots,d_k\}\subset\CC{1...D}$ and $\bs:=(s_1,\ldots,s_d)\in\{\pm1\}^k$, then
 $\partial_\bd^\bs \ \CornerTile{\fz}{\ }:= 
\partial_{d_1}^{s_1}\, \CornerTile{\fz}{\ }\intsct\cdots\intsct\partial_{d_k}^{s_k}\,\CornerTile{\fz}{\ }$).  This yields a natural $D$-cell decomposition $\bY$
of $\RD$, whose zero-skeleton $\bY^0$ is $\ZD$, and whose $D$-skeleton
$\bY^D$ is the set of $D$-cubes $\set{\CornerTile{\fz}{\ }}{\fz\in\ZD}$.
For any $d\in\CO{1...D}$, the $d$-skeleton $\bY^d$ is the set of all
$d$-dimensional edges/faces/etc. of the cubes in
$\set{\CornerTile{\fz}{\ }}{\fz\in\ZD}$.  We will call this the
{\dfn canonical cell complex} for $\RD$.

\subsection{Review of Cubic (co)homology:  \label{S:cubic.cohom}}
For any $d\in\CC{0...D}$, let $\Zahl[\bY^d]$ be the free abelian group
generated by $\bY^d$,  i.e. the group of 
formal $\Zahl$-linear combinations
of $d$-cells in $\bY$. Elements of $\Zahl[\bY^d]$ are called {\dfn $d$-chains}.
Since $\bY^{-1}=\emptyset=\bY^{D+1}$, 
we also formally define  $\Zahl[\bY^{-1}]:=\{0\}=:\Zahl[\bY^{D+1}]$.
Let $\partial_d:\Zahl[\bY^d]\into\Zahl[\bY^{d-1}]$
be the {\dfn boundary homomorphism} from cubic homology theory.  For example:

  $\bullet$ \  For any $\fz\in\ZD$, let $\dot\fz\in\bY^0$ be the
corresponding zero-cell (i.e. vertex).  Then $\partial_0(\dot\fz):=0$.  

  $\bullet$ \  If $\fy,\fz\in\ZD$ are adjacent,
and $(\onecell{\fy}{\fz})\in\bY^1$ is the one-cell (i.e. oriented edge) from 
$\dot\fy$ to $\dot\fz$, then 
$\partial_1(\onecell{\fy}{\fz}):=\dot\fz-\dot\fy$. 

  $\bullet$ \  If 
$\twocell{\fw}{\fx}{\fy}{\fz}\in\bY^2$ is the two-cell (i.e. oriented square)
whose four corner vertices are $\dot\fw,\dot\fx,\dot\fy,\dot\fz\in\bY^0$, then
$\partial_2 \ (\twocell{\fw}{\fx}{\fy}{\fz} )\ := \ 
(\onecell{\fw}{\fx}) + (\onecell{\fx}{\fy}) + (\onecell{\fy}{\fz}) +
 (\onecell{\fz}{\fw})$.

  $\bullet$ \  If $\threecell{\fs}{\ft}{\fu}{\fv}{\fw}{\fx}{\fy}{\fz}$ is
a three-cell (i.e. oriented cube), then
$\partial_3 (\threecell{\fs}{\ft}{\fu}{\fv}{\fw}{\fx}{\fy}{\fz})=
\twocell{\fs}{\ft}{\fu}{\fv} + \twocell{\ft}{\fx}{\fy}{\fu}+
\twocell{\fx}{\fw}{\fz}{\fy} + \twocell{\fw}{\fs}{\fv}{\fz}
+ \twocell{\fv}{\fu}{\fy}{\fz} + \twocell{\fx}{\fw}{\fs}{\ft}$.

In general, $\partial_d$ can be computed by decomposing a cubic $d$-cell as
a formal sum of $d$-simplices, computing the boundary of
the resulting simplicial $d$-chain as in standard simplicial homology \cite[\S2.1]{Hatcher}, and then expressing the result as a sum of $(d-1)$-cubes.

Let $(\sG,+)$ be an abelian group, and let $\sC^d=\sC^d(\bY,\sG)$ be the set of
($d$-dimensional, $\sG$-valued) {\dfn cochains}:  i.e. homomorphisms
$C:\Zahl[\bY^d]\into\sG$ (or equivalently, arbitrary functions
$C:\bY^d\into\sG$).  We define the {\dfn coboundary}
homomorphism $\del_d:\sC^d\into\sC^{d+1}$ as follows:
If $C\in\sC^d$, then $\del_d C\in\sC^{d+1}$ is defined by
$\del_d C(\zeta):= C(\partial_{d+1} \zeta)$ for any $\zeta\in\Zahl[\bY^{d+1}]$.
We say $C$ is a {\dfn cocycle}
if $\del_d C\equiv 0$ (i.e.
$C(\partial_{d+1} \zeta)=0$ for any $\zeta\in\Zahl[\bY_{d+1}]$).
Let $\sZ^d:=\ker(\del_d)$ be the group of cocycles.
We say $C$ is a {\dfn coboundary} if there is some
{\dfn cobounding function} $b\in\sC^{d-1}$ such that $C=\del_{d-1} b$.
Let $\sB^d:=\del_{d-1}(\sC^{d-1})$
 be the group of coboundaries. 
Two cocycles $C,C'\in\sZ^d$
are {\dfn cohomologous} (notation: $C\approx C'$)
if there is some coboundary $B:=\del_{d-1}b$
so that $C' = C + B$.  Let $\undC$ denote the cohomology class of $C$.

\example{\label{X:cubic.cocycle}
If $d=1$, then $C:\bY^1\into\sG$ is a cocycle iff, 
for any two-cell $\twocell{\fw}{\fx}{\fy}{\fz}$, we have
$C(\onecell{\fw}{\fx}) + C(\onecell{\fx}{\fy}) + C(\onecell{\fy}{\fz}) +
 C(\onecell{\fz}{\fw})=0$.   Equivalently, 
$C(\onecell{\fw}{\fx}) + C(\onecell{\fx}{\fy}) = C(\onecell{\fw}{\fz}) +
 C(\onecell{\fz}{\fy})$.  By induction, this is equivalent to saying:
for any $\fw,\fz\in\ZD$ and any two chains
$\zeta,\zeta'\in\Zahl[\bY^1]$ with
$\partial_1(\zeta)=(\dot\fy-\dot\fw)=\partial_1(\zeta')$ (i.e. any two `paths' from
$\dot\fw$ to $\dot\fy$), we have
$C(\zeta)=C(\zeta')$.  Thus, $C$ defines a function
$\ddot{C}:\ZD\x\ZD\into\sG$, by $\ddot{C}(\fw,\fy):=C(\zeta)$, where
$\zeta\in\Zahl[\bY^1]$ is any 1-chain with $\partial_1(\zeta)=(\dot\fy-\dot\fw)$.
The function $\ddot{C}$ is a {\dfn two-point cocycle}, by which we mean:
\beqn
\label{X:cubic.cocycle.eqn}
\forall \ \fw,\fx,\fy\in\ZD,\quad
{\bf(a)} \ \ddot{C}(\fy,\fw)=-\ddot{C}(\fw,\fy)
\quad\And\quad
{\bf(b)} \ \ddot{C}(\fy,\fw)=\ddot{C}(\fy,\fx) + \ddot{C}(\fx,\fw)
\eeqn
Conversely, any two-point cocycle $\ddot{C}:\ZD\x\ZD\into\sG$ 
defines a cocycle $C:\bY^1\into\sG$ by  $C(\onecell{\fx}{\fy})
:= \ddot{C}(\fx,\fy)$.
  
Also, $C\in\sB^1$ iff
there exists $b:\ZD\into\sG$ such that, for any one-cell
 $(\onecell{\fy}{\fz})$, we have $C(\onecell{\fy}{\fz}) = b(\fz)-b(\fy)$.
Thus, $C\approx C'$ iff 
there exists $b:\ZD\into\sG$ with
$C'(\onecell{\fy}{\fz}) = b(\fz) +  C(\onecell{\fy}{\fz}) - b(\fy)$.
}
\ignore{
(b) If $d=0$, then $C:\ZD\into\sG$ is a {\dfn cocycle} iff
for any one-cell $\onecell{\fy}{fz}$, we have $C(\fz)-C(\fy)=0$;
i.e. $C(\fy)=C(\fz)$; i.e. $C$ is constant on $\ZD$.  The
set $\sB^0$ is trivial (because there are no $-1$-cells).
}

\subsection{Homotopy/Homology groups for Wang Tiles:\label{S:wang.homo}}
Let $\sW$ be a set of Wang tiles, and let $\gW\subset\WZD$ be the
corresponding Wang subshift.  If $w_1,w_2\in\sW$
and  $d\in\CC{1...D}$,  then we will write
``$w_1\edgematch w_2$'' to mean that the face 
$\partial_d^+ w_1$ is compatible with the face $\partial_d^- w_2$.
The {\dfn tile complex} of $\gW$
is defined as follows:  for each $\fz\in\ZD$ and $w\in\sW$, let
$\CornerTile{\fz}{w}$ be a $D$-cell, which we imagine
as a $D$-dimensional unit cube `labelled' by $w$, with
its minimal-coordinate corner `over' $\fz$.  Let 
\[
\tlbX\quad:=\quad\Disj_{\fz\in\ZD} \Disj_{w\in\sW} \CornerTile{\fz}{w}.
\]
For all $d\in\CC{1...D}$, let  
$\partial_d^{\pm} \  \CornerTile{\fz}{w}$
be the two $(D-1)$-dimensional faces of $\CornerTile{\fz}{w}$
which are orthogonal to $\fe_d$.  
For any $\fz\in\ZD$ and $\fz':=\fz+\fe_d$, and any
$w,w'\in\sW$, we will `glue together' the faces
 $\partial_d^{+}\ \CornerTile{\fz}{w}$
and $\partial_d^{-}\ \CornerTile{\fz'}{w'}$ if and only if
$w\edgematch w'$.  Let $\sim$ be the equivalence relation
on $\tlbX$ which instantiates all these gluing operations.  Then
 $\bX:=\tlbX/\!\sim$ is a $D$-dimensional cell
complex with the following properties:
\bitem
\item[(a)] The $D$-cells of $\bX$ are in bijective correspondence
with the cubes comprising $\tlbX$.

\item[(b)] There is a natural continuous surjection $\Pi:\bX\into\RD$,
such that $\Pi\restr{}:\CornerTile{\fz}{w} \into\CornerTile{\fz}{ \ }$,
is a homeomorphism for each $\fz\in\ZD$ and $w\in\sW$. 
(Here, $\CornerTile{\fz}{ \ }$ is as in \S\ref{S:canonical.cell.cplx}.)

\item[(c)]   $\Pi$ is a {\dfn cellular map.}
For each $d\in\CC{0...D}$, let $\bY^d$ be the $d$-skeleton of $\RD$ from
\S\ref{S:canonical.cell.cplx},  let
 $\bX^d$ be the $d$-skeleton of $\bX$,
and let $\Pi^d:=\Pi\restr{\bX^d}$. Then $\Pi^d:\bX^d\surject \bY^d$
surjectively.
In particular, $\Pi^0:\bX^0\surject\bY^0=\ZD$ is a surjection
(and in many cases, a bijection).

\item[(d)] Any continuous {\dfn section} of $\Pi$ (i.e. a function
$\varsigma:\RD\into\bX$ such that $\Pi\circ\varsigma=\Id{\RD}$) assigns
a unique $w\in\sW$ to each $D$-cube in $\RD$, and thereby defines a
tiling $\bw_\varsigma\in\gW$.  Conversely, any admissible tiling
$\bw\in\gW$ determines a continuous section $\varsigma_\bw$ of $\Pi$.

\item[(e)] For each $\fz\in\ZD$, let $\Xi^{\fz}:\bX\into\bX$ be the
self-homeomorphism of $\bX$ induced by translating all cells by $\fz$
in the obvious way [i.e. for any $\fx\in\ZD$ and $w\in\sW$, let
 $\Xi^\fz\lb(\CornerTile{\fx}{w}\rb) :=
\CornerTile{\fy}{w}$, where $\fy:=\fx+\fz$]. This defines a homeomorphic $\ZD$-action on $\bX$.
Then $\Pi\circ\Xi^{\fz}=\shift{\fz}\circ\Pi$, and
for any $\bw\in\gW$, and $\fy\in\RD$,
$\varsigma_{\shift{\fz}(\bw)}(\fy) = \Xi^{-\fz}\circ \varsigma_\bw(\fy+\fz)$.
\eitem
Fix $x\in\bX^0$. For any $d\in\CC{1...D}$, let $\barpi_d(\sW)
:=\pi_d(\bX,x)$ be the  {\dfn $d$th homotopy group} of $\sW$.
Let $(\sG,+)$ be an abelian group, 
and let $\barsH_d(\sW,\sG):=\sH_d(\bX,\sG)$ and
$\barsH^d(\sW,\sG):=\sH^d(\bX,\sG)$
be the {\dfn $d$th homology group} and
{\dfn cohomology group} respectively (with coefficients in $\sG$).
We will briefly review how to construct these groups in terms of
the cellular structure of $\bX$.
Any element of $\barpi_d(\sW)$ can be represented as a continuous
function $\xi:(\dS^{d},s)\into(\bX^{d},x)$ (where $(\dS^{d},s)$
is as in \S\ref{S:proj.codim}),  and two such functions 
represent the same element of $\barpi_d(\sW,x)$ if and only if they are
homotopic in $(\bX^{d+1},x)$ (in a basepoint-fixing way); see
\cite[Corollary 4.12]{Hatcher}.
Let $\sC_d:=\Dirsum_{x\in\bX_d} \sG$ be the group
of {\em $d$-dimensional $\sG$-chains}
(i.e.  functions $\bX^d\into\sG$ with only finitely many
nontrivial entries).
There is a natural `boundary' homomorphism $\partial_d:\sC_d
\into \sC_{d-1}$ (see \S\ref{S:cubic.cohom}), and
$\barsH_d(\sW,\sG) := \sZ_d/\sB_d$,
where $\sZ_d:=\ker(\partial_d)$ and $\sB_d:=\image{\partial_{d+1}}$.
Let $\sC^d:=\sG^{\bX_d}$ be the 
group of {\em $d$-dimensional $\sG$-cochains}
i.e. all functions $\bX^d\into\sG$, or equivalently, all homomorphisms
$\Zahl[\bX^d]\into\sG$.
There is a natural `coboundary' homomorphism $\del_d:
\sC^d \into \sC^{d+1}$ defined by $\del_d(\eta)=\eta\circ\partial_d$.
Then $\barsH^d(\sW,\sG) := \sZ^d/\sB^d$, where 
$\sZ^d:=\ker(\del_d)$ and $\sB^d:=\image{\del_{d-1}}$.

\begin{figure}
\centerline{\includegraphics[scale=1]{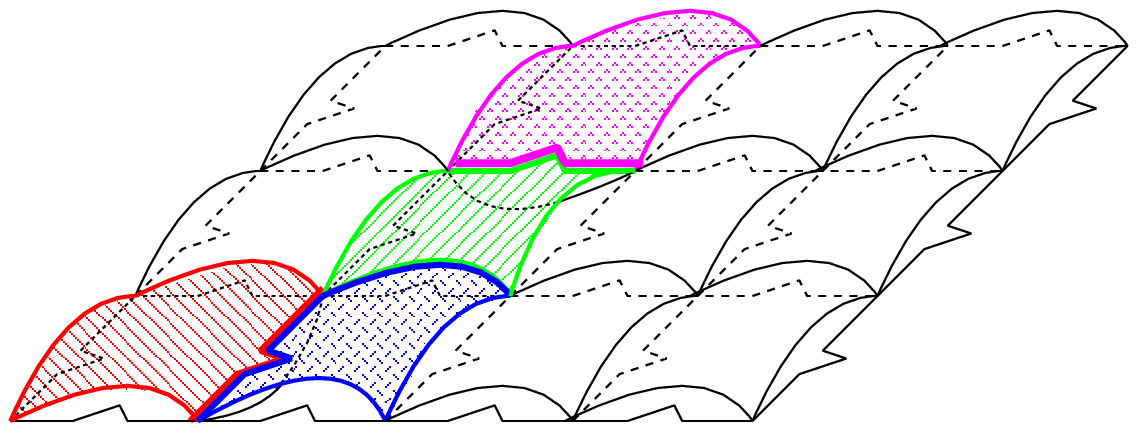}}
\caption{{\footnotesize A fragment of the tile complex $\bX$ for $\Dom$.
The elements of $\bX^0$ are in bijective correspondence with $\Zahl^2$.
Between each pair of adjacent vertices in $\bX_0$, we adjoin two edges:
one `straight' edge, and one `notched' edge;  we define $\bX^1$ 
to be the union of all these edges. 
In every square bounded by four vertices
of $\bX_0$, there are four distinct 2-cells, each of which has exactly
three `straight' edges and one `notched' edge (four such 2-cells are
depicted in the figure).  We define $\bX^2$ to be the union of all
such 2-cells.  
\label{fig:domino.tile.complex}}}
\end{figure}

\example{\label{X:Conway.Lagarias}
(a) Let $D=2$, so that $\sW$ is a set of two-dimensional square
tiles with edge-matching conditions (e.g. edge `colours'). 
Then $\barpi_1(\sW)$ is the Conway-Lagarias
{\dfn tile homotopy group} of \cite[\S3]{Cola}; see also \cite{Thu},
\cite{Pro} or \cite[\S4]{Re03}.

To see this, note that $\bX$ is a two-dimensional cell-complex
obtained by taking a collection $\tlbX$ of $\sW$-labelled unit squares and
gluing them along their edges in accordance with the edge-matching
conditions of $\sW$.  For example, Figure \ref{fig:domino.tile.complex}
shows a fragment of the tile complex for the domino tiles $\Dom$ from
Example \ref{X:codim.one}(b).
Let $\sH$ (resp. $\sV$) be the set of `colours'
of horizontal (resp. vertical) tile edges in $\sW$.  Assume that
$\sH$ and $\sV$ are disjoint, and let $\sC$ be the free group
generated by $\sH\disj\sV$.  Any element of $\barpi_1(\sW)$ corresponds
to a continuous function $\xi:\dS^1\into\bX^1$ ---i.e. a closed
continuous path along the edges of the tile complex, beginning and
ending at zero.  The function $\xi$ defines an element
$c_1^{\pm1} c_2^{\pm1} \cdots c_k^{\pm1}\in\sC$, where $c_1$ is the
colour of the first edge traversed by $\xi$ and we put $c^{+1}$ if
$\xi$ heads east or north along this edge, and $c^{-1}$ if $\xi$ heads
west or south.  Likewise $c_2$ the colour of the second edge, with the
same sign convention, and so on.  The word $c_1^{\pm1} c_2^{\pm1} \cdots c_k^{\pm1}$ satisfies two constraints:
 \bitem
 \item[{[i]}]  If $N$, $E$, $S$, and $W$ are the total \# of northward,
eastward, southward, and westward edges (as indicated by the colours and
sign conventions), then $N=S$ and $E=W$.

 \item[{[ii]}] $c_1$ and $c_k$ must be the colours of edges coming into
or out of the vertex $x$. 
\eitem
(If $\Pi^0:\bX^0\into\ZD$ is bijective, as is the
case in \cite{Cola,Re03,Thu,Pro}, then condition [ii] is trivial.)
Let $\sD\subset\sC$ be the subgroup of elements satisfying
[i] and [ii].  Let $\sN$ be the normal subgroup of $\sC$
generated by all words of the form $sen^{-1}w^{-1}$, where 
$w$, $n$, $e$, and $s$ are the four edge colours of any tile in $\sW$.
Let $\sG:=\sD/\sN$.  Then $\sG$ is the Conway-Lagarias group,
and $\sG\cong\barpi_1(\sW)$ (because a nullhomotopy of $\xi$
is equivalent to an algebraic reduction of
$c_1^{\pm1} c_2^{\pm1} \cdots c_k^{\pm1}$ to an element 
of $\sN$).

Similarly, $\barsH_1(\sW,\Zahl)$ is the Conway-Lagarias
{\dfn tile homology group} \cite[\S5]{Cola}, \cite[\S2]{Re03}.

(b)  Let $\sI$ be as in Example \ref{X:codim.one}(a).
We apply (a) to show that $\barpi_1(\sI)\cong\Zahl$.
In this case, $\sH:=\{A,V\}$ and $\sV:=\{\sqw,\sqe\}$, and $\sN$ is the normal
subgroup generated by 
\beqn
\label{X:Conway.Lagarias.ice.eqn}
\{A\sqw A^{-1}\sqw ^{-1}, \ \  A\sqe A^{-1}\sqe ^{-1}, \ \  V\sqw V^{-1}\sqw ^{-1}, \ \
V\sqe V^{-1}\sqe ^{-1},  \  \ V\sqw A^{-1}\sqe ^{-1}, \  \  A\sqe V^{-1}\sqw ^{-1}\}. 
\eeqn
Let $\barA,\barV,\bar\sqe $ and $\bar\sqw $ be the corresponding generators
of $\sC/\sN$. The
first four generators in eqn.(\ref{X:Conway.Lagarias.ice.eqn})
make $\sC/\sN$ abelian, so we will
switch to additive notation.  Thus, any element of $\sC/\sN$ has
the form $a\barA + v\barV + p\bar\sqe  + q\bar\sqw $ for some $a,v,p,q\in\Zahl$,
and $\sD/\sN$ is the subgroup of elements satisfying $a=-v$ and $p=-q$
(from condition [i]).
Now the last two generators in eqn.(\ref{X:Conway.Lagarias.ice.eqn})
both say that $\bar\sqw =\barA-\barV + \bar\sqe $.
Thus, $a\barA + v\barV + p\bar\sqe  + q\bar\sqw  = (a+q)\barA + (v-q) \barV + (p+q)\bar\sqe  
= (a+q)\barA + (-a-q) \barV = (a+q)(\barA-\barV)$.  Thus, if 
we define $G:=\barA-\barV$, then
$\barpi_1(\sI)\cong \sD/\sN$ is the cyclic group generated by $G$;
hence  $\barpi_1(\sI)\cong \Zahl$.
}

\subsection{Wang representations:\label{S:Wang.representation}}

We can represent any subshift of finite type using Wang tiles.
To do this requires some notation.  For any $r\in\Natur$, recall
that $\dB(r):=\CC{-r...r}^D\subset\ZD$ is the cube of sidelength $(2r+1)$.
For any $d\in\CC{1...D}$, we define the `top' and `bottom' faces of this cube
in the $d$th dimension by:
\[
\partial_d^+\, \dB(r) := \CC{-r...r}^{d-1} \x \OC{-r...r} \x  \CC{-r...r}^{D-d-1}
\ \And \ 
\partial_d^-\, \dB(r) := \CC{-r...r}^{d-1} \x \CO{-r...r} \x  \CC{-r...r}^{D-d-1}
\]
\ignore{For any $k\in\CC{2...D}$, we define
codimension-$k$  faces by intersecting $k$ of these codimension-one
faces.  If $\bd:=\{d_1,\ldots,d_k\}\subset\CC{1...D}$
and $\bs:=(s_1,\ldots,s_d)\in\{\pm1\}^k$, then we define
$\partial_\bd^\bs \ \dB(r) :=
\partial_{d_1}^{s_1}\, \dB(r)\intsct\cdots\intsct\partial_{d_k}^{s_k}\,\dB(r)$.}  Now, let $\gA\subset\AZD$ be a subshift of finite type defined by a set
$\gA_{(R)}\subset\sA^{\dB(R)}$ of admissible $R$-blocks.  For any
$r\geq R$, the {\dfn radius $r$ Wang representation} of $\gA$ is
defined as follows:  let $\sW_r:=\gA_{(r)}$, and for any 
$\ba,\bb\in\sW_r$, and any $d\in\CC{1...D}$, we allow
$\ba \edgematch  \bb$ if and only if
$\ba_{\partial_d^+ \dB(r)} \ = \ \bb_{\partial_d^- \dB(r)}$.
For example, suppose $D=2$ and $r=1$; and suppose
\[
\ba= {\scriptsize\lb[\begin{array}{ccc}
 a_{-1,1} & a_{0,1} & a_{1,1}\\
 a_{-1,0} & a_{0,0} & a_{1,0}\\
 a_{-1,-1} & a_{0,-1} & a_{1,-1}
\end{array} \rb]}
\quad\And\quad
\bb= {\scriptsize\lb[\begin{array}{ccc}
 b_{-1,1} & b_{0,1} & b_{1,1}\\
 b_{-1,0} & b_{0,0} & b_{1,0}\\
 b_{-1,-1} & b_{0,-1} & b_{1,-1}
\end{array} \rb]}
\]
Then $\ba  \edgematch[1] \bb$ iff
$ {\scriptsize\lb[\begin{array}{cc}
  a_{0,1} & a_{1,1}\\
  a_{0,0} & a_{1,0}\\
  a_{0,-1} & a_{1,-1}
\end{array} \rb]}
 =  
 {\scriptsize\lb[\begin{array}{ccc}
 b_{-1,1} & b_{0,1} \\
 b_{-1,0} & b_{0,0} \\
 b_{-1,-1} & b_{0,-1} 
\end{array} \rb]}$,
and
 $\ba  \edgematch[2] \bb$ iff
${\scriptsize\lb[\begin{array}{ccc}
 a_{-1,1} & a_{0,1} & a_{1,1}\\
 a_{-1,0} & a_{0,0} & a_{1,0}\\
\end{array} \rb]}
  =   {\scriptsize\lb[\begin{array}{ccc}
 b_{-1,0} & b_{0,0} & b_{1,0}\\
 b_{-1,-1} & b_{0,-1} & b_{1,-1}
\end{array} \rb]}$.

\subsection{Projective Homotopy/Homology groups for SFTs 
\label{S:projective.homotopy}}
  Using the Wang representation of \S\ref{S:Wang.representation}, we can define homotopy/homology
groups for any subshift of finite type as in \S\ref{S:wang.homo}.
There are two problems, however:
 \bitem
  \item[{[i]}]  There are many different Wang tile representations for
any SFT, and none of them is `canonical'.  Different Wang representations
may yield non-isomorphic groups.

  \item[{[ii]}]  Wang tile representations (and hence, the corresponding
homotopy/homology groups) do not behave well under subshift homomorphisms (e.g.
cellular automata).
\eitem
  To obviate these problems, we use an inverse limit
which encompasses `all possible' Wang tile
representations within a single algebraic structure.

\breath

Throughout this section, let $\gA\subset\AZD$ be a subshift of finite
type of radius $R$.  For any $r\geq R$, let $\sW_r:=\gA_{(r)}$, and
let $\gW_r\subset\sW_r^{\ZD}$ be the {Wang representation} of $\gA$
from \S\ref{S:Wang.representation}.  Let
$\bX_r=(\bX_r^0,\ldots,\bX_r^D)$ be the corresponding tile complex;
hence the $D$-cells of $\bX_r$ have the form $\CornerTile{\fz}{\bb}$,
where $\fz\in\ZD$ and $\bb\in\gA_{(r)}$ is an $\gA$-admissible
$\dB(r)$-block.  Let $\Pi_r:\bX_r\into\RD$ be the natural projection
map; then there is a natural bijective correspondence between:
\bitem
  \item $\gA$-admissible configurations in $\AZD$.
  \item $\sW_r$-tilings of $\RD$ 
(satisfying the relevant edge-matching constraints).
   \item Continuous sections $\varsig:\RD\into\bX_r$ of $\Pi_r$.
\eitem
Fix $\ba\in\gA$, and let $\ba_r := \ba_{\dB(r)}$,
so that $\CornerTile{0}{\ba_r}$ is a $D$-cell in $\bX_r$.
Let $x_r=x_r(\ba)$ be the unique element of the singleton set
$\Pi_r^{-1}\{0\}\intsct \CornerTile{0}{\ba_r}$; 
then $x_r$ is a corner vertex of $\CornerTile{0}{\ba_r}$. Define
$\pi^r_d(\gA,\ba):=\pi_d(\bX_r,x_r)$.  For example, if $\sW$ is a set of
Wang tiles and $\gW\subset\WZD$ is the corresponding Wang subshift,
then $\pi^0_d(\gW,\bw) = \barpi_d(\sW)$.

There is a natural
continuous surjection $\zeta_r:\bX_{r+1}\into\bX_r$ where,
for any $\ba\in\gA_{(r+1)}$, if $\ba':=\ba_{\dB(r)}$, then
for any $\fz\in\ZD$, $\zeta\restr{}:\CornerTile{\fz}{\ba}\into\CornerTile{\fz}{\ba'}$ is a homeomorphism.  Furthermore, $\zeta_r$ is a {\dfn cellular map}: i.e.
 $\zeta(\bX_{r+1}^d)\subseteq \bX_r^d$ for all $d\in\CC{0...D}$.
Also, $\zeta_r(x_{r+1})=x_r$.
This induces homomorphisms $\pi_d\zeta_r:\pi^{r+1}_d(\gA,\ba)\into\pi^r_d(\gA,\ba)$ for all $d\in\CC{1...D}$.  
We define
\[
  \pi_d(\gA,\ba)\quad:=\quad \invlim 
\lb( \pi^R_d(\gA,\ba) \xleftarrow{\pi_d\zeta_R}
\pi^{R+1}_d(\gA,\ba) \xleftarrow{\pi_d\zeta_{R+1}}
\pi^{R+2}_d(\gA,\ba) \xleftarrow{\pi_d\zeta_{R+2}}
\cdots\rb)
\]
[Note:  $\pi_d(\gA,\ba)$ is {\em not} the homotopy group of $\gA$ as
a (zero-dimensional) topological space.]
For example, if $d=1$, then $\pi_1(\gA)$ is the
 {\dfn projective fundamental group} of \cite{GePr} (the cell-complex
$\bX_r$ is a dual version of the cellular realization of `scenery space'
 on \cite[p.1101]{GePr}).   We would like this definition
to be independent of the choice of `basepoint' $\ba$.
This will be true if $\gA$ is {\em projectively connected} 
in the sense of \cite[p.1098]{GePr}, but we also have the following
criteria.  Recall that a subshift $\gA\subset\AZD$ is {\dfn
topologically weakly mixing} if the Cartesian product system
$(\gA\x\gA, \shift{}\x\shift{})$ is topologically transitive.

\Proposition{\label{basepoint.free}}
{
Fix $\ba,\bb\in\gA$.
\bthmlist
  \item Suppose $\Pi_r^0:\bX_r^0\into\ZD$ is injective for all 
large enough $r\in\Natur$. Then  for  any $d\geq 1$,
 there is a canonical isomorphism
 $\pi_d(\gA,\ba) \cong \pi_d(\gA,\bb)$.
\ethmlist
 Suppose $(\gA,\shift{})$ is topologically weakly mixing.  Then:
\bthmlist
\setcounter{enumi}{1}
  \item  If $\pi_1(\gA,\ba)$ is abelian, 
then  there is a canonical isomorphism
 $\pi_1(\gA,\ba) \cong \pi_1(\gA,\bb)$.

  \item If $\pi_1(\gA,\ba)$ is trivial,
then there are canonical isomorphisms
$\pi_d(\gA,\ba) \cong \pi_d(\gA,\bb)$
for all $d\in\Natur$.
\ethmlist
}
\bthmprf {\bf(a)} If $\Pi_r^0:\bX_r^0\into\ZD$ is injective, then
$(\Pi_r^0)^{-1}\{0\}$ is a singleton, which means
$x_r(\ba)=x_r(\bb)$. Hence,  
 $\pi^r_d(\gA,\ba) := \pi_d[\bX_r,x_r(\ba)] = \pi_d[\bX_r,x_r(\bb)] =:
\pi^r_d(\gA,\bb)$.  If this holds for all large $r\in\Natur$ 
then clearly $\pi_d(\gA,\ba)= \pi_d(\gA,\bb)$.  This proves {\bf(a)}.
For {\bf(b,c)} we need:

\Claim{If $(\gA,\shift{})$ is topologically weakly mixing, then
for all $r\geq R$, the space $\bX_r$ is path connected.}
\bclaimprf  
Fix $\bc\in\gA_{(r)}$, and let $[\bc]:=\set{\bd\in\gA}{
\bd_{\dB(r)}=\bc}$ be the corresponding cylinder set.  
Likewise, let $[\ba_{\dB(r)}]$ and $[\bb_{\dB(r)}]$ be the
cylinder sets defined by $\ba_{\dB(r)}$ and $\bb_{\dB(r)}$.
The system
$(\gA\x\gA,\shift{}\x\shift{})$ is transitive (because 
$(\gA,\shift{})$ is weakly mixing), so there is some $\fz$
such that $(\shift{\fz}[\ba_{\dB(r)}] \x \shift{\fz}[\bb_{\dB(r)}])
\intsct ([\bc]\x[\bc])\neq\emptyset$.  Hence,
there exist $\ba',\bb'\in\gA$ such that
$\ba'_{\dB(r)}=\ba_{\dB(r)}$ and $\bb'_{\dB(r)}=\bb_{\dB(r)}$,
but also $\ba'_{\dB(\fz,r)}=\bc=\bb'_{\dB(\fz,r)}$.
Clearly, $x_r(\ba')=x_r(\ba)$ and $x_r(\bb')=x_r(\bb)$.

Let $\varsig^r_{\ba'},\varsig^r_{\bb'}:\RD\into\bX_r$ 
be the continuous sections of $\Pi_r$ defined by $\ba'$ and
$\bb'$.  Let $\gam:\CC{0,1}\into\Real$ be a continuous path, with
$\gam(0)=0$ and $\gam(1)=\fz$.  If $\alp := \varsigma^r_{\ba'}\circ\gam$,
then $\alp(0)=x_r(\ba)$ and $\alp(1)=\Xi^{\fz}(x_r(\bc))$.
Likewise, if  $\bet := \varsigma^r_{\bb'}\circ\gam$.
then $\bet(0)=x_r(\bb)$ and $\bet(1)=\Xi^{\fz}(x_r(\bc))$.
If $\gam_r:=\bkw{\bet}\star\alp:\CC{0,1}\into\bX_r$, then
$\gam_r(0)=x_r(\ba)$ and $\gam_r(1)=x_r(\bb)$,
as desired.
\eclaimprf
 Let $r\geq R$.
 Any path $\gam:\CC{0,1}\into\bX_r$ from $x_r(\ba)$ to
$x_r(\bb)$ (such as in Claim 1)
yields an isomorphism $\gam_*:\pi_d(\gA,\ba)\into\pi_d(\gA,\bb)$
(see \S\ref{S:proj.codim}).
If $\eta:\CC{0,1}\into\bX_r$ is another path from $x_r(\ba)$ to $x_r(\bb)$,
and $\gam\homoto\eta$, then $\gam_* =\eta_*$. 

{\bf(c)} \ If $\pi_1^r(\gA,\ba)=\pi_1(\bX_r,x_r(\ba))$ is trivial, then any two
such such paths $\gam$ and $\eta$ are homotopic.  Hence in this case
there is a canonical isomorphism $I^r_d:\pi^r_d(\gA,\ba)\into\pi^r_d(\gA,\bb)$,
which is independent of the choice of path.
This yields a commuting ladder with canonical isomorphisms for rungs:
\[
\Array{
\pi^{R}_d(\gA,\ba) &\xleftarrow{\pi_d\zeta_{R}}&\pi^{R+1}_d(\Sft,\ba) &\xleftarrow{\pi_d\zeta_{R+1}}&\pi^{R+2}_d(\gA,\ba) &\xleftarrow{\pi_d\zeta_{R+2}}\cdots\hspace{-0.5em} &\mbox{\footnotesize which yields }\hspace{-0.5em} & \pi_d(\gA,\ba)\\
I^R_d\Longdownarrow &      & I^{R+1}_d\Longdownarrow &      & I^{R+2}_d\Longdownarrow && \begin{array}{cc}\mbox{\footnotesize a canonical}\\ \mbox{\footnotesize isomorphism}\end{array} & \Longdownarrow\\
\pi^{R}_d(\gA,\bb) &\xleftarrow[\pi_d\zeta_{R}]{}&\pi^{R+1}_d(\gA,\bb) &\xleftarrow[\pi_d\zeta_{R+1}]{}& \pi^{R+2}_d(\gA,\bb)
&\xleftarrow[\pi_d\zeta_{R+2}]{}\cdots\hspace{-0.5em} & \mbox{\footnotesize of colimits:}\hspace{-0.5em} & \pi_d(\gA,\bb). \\}
\]
{\bf(b)} \ If $\pi_1(\bX_r,x_r(\ba))$ is not trivial, and
$\gam$ and $\eta$ are non-homotopic paths from 
$x_r(\ba)$ to $x_r(\bb)$, then in general $\gam_*\neq\eta_*$.
Hence, $\eta_*^{-1}\circ\gam_*$ will be a nontrivial
automorphism of $\pi_1(\bX_r,x_r(\ba))$.  Indeed, if $\alp:=\bkw{\eta}\!\!\!\star \gam$,
then $\alp$ is a closed
loop based at $x_r(\ba)$, hence $\undalp\in\pi_1(\bX_r,x_r(\ba))$,
and the automorphism 
$\eta_*^{-1}\circ\gam_* = (\bkw{\eta}\!\!\!\star\gam)_* = \alp_*:
\pi_1(\bX_r,x_r(\ba)) \into \pi_1(\bX_r,x_r(\ba))$
is simply the inner automorphism 
$\alp_*(\undbet) = \undalp^{-1} \cdot \undbet \cdot \undalp$.  
But if $\pi_1(\bX_r,x_r(\ba))$ is abelian, then all inner automorphisms
are trivial; hence $\alp_*=\Id{}$, hence $\gam_*=\eta_*$
after all.  Thus, the isomorphism  
$I^r_1:\pi^r_1(\gA,\ba)\into\pi^r_1(\gA,\bb)$ once again well-defined
independent of the choice of path from $x_r(\ba)$ to $x_r(\bb)$.
Now proceed as in {\bf(b)}.
\ethmprf

  If any of the conditions of Proposition \ref{basepoint.free}
 is satisfied, then 
we say that $\gA$ is {\dfn basepoint-free} in codimension $d+1$.
We will then write ``$\pi_d(\gA)$'' to mean ``$\pi_d(\gA,\ba)$'',
where $\ba\in\gA$ is arbitrary.

\example{Let $\Ice$ be as in  Example \ref{X:codim.one}(c).
Then $\pi_1(\Ice)=\Zahl$ \cite[Theorem 3]{GePr}.}

\ignore{We will sketch a direct proof.  

\Claim{For any $r>0$, $\pi_1^r(\Ice)=\Zahl$.}
\bclaimprf
  The case $r=1$ was just Example \ref{X:Conway.Lagarias}(b).
For any $r\geq 1$, \  $\pi_1^r(\Ice)$ is the Conway-Lagarias
homotopy group of the tile-set 
\eclaimprf}

The (co)homology groups do not require a basepoint.
For any $d\in\CC{0...D}$, any $r\geq R$,
and any abelian group $(\sG,+)$,
we define $\sH^r_d(\gA,\sG):=\sH_d(\bX_r,\sG)$
and $\sH^d_r(\gA,\sG):=\sH^d(\bX_r,\sG)$
(see \S\ref{S:wang.homo}).
For example, if $\sW$ is a set of
Wang tiles and $\gW\subset\WZD$ is the corresponding Wang subshift,
then $\sH^0_d(\gW,\sG) = \barsH_d(\sW,\sG)$ and
$\sH^d_0(\gW,\sG) = \barsH^d(\sW,\sG)$.
The cellular maps $\zeta_r$ induce
homomorphisms $\sH_d\zeta_r\colon\sH^{r+1}_d(\gA)\rightarrow\sH^r_d(\gA)$ for all $d\in\CC{0...D}$.
We define
\[
  \sH_d(\gA,\sG)\quad:=\quad \invlim
\lb( \sH^R_d(\gA,\sG) \xleftarrow{\sH_d\zeta_R}
\sH^{R+1}_d(\gA,\sG) \xleftarrow{\sH_d\zeta_{R+1}}
\sH^{R+2}_d(\gA,\sG) \xleftarrow{\sH_d\zeta_{R+2}}
\cdots\rb)
\]
The functions $\zeta_r$ also induce (contravariant)
homomorphisms $\sH^d\zeta_r\colon\sH_{r}^d(\gA,\sG)\rightarrow\sH_{r+1}^d(\gA,\sG)$ for all $d\in\CC{0...D}$.
We define
\beqn
\label{tile.cohom.defn}
  \sH^d(\gA,\sG) \ := \ \dirlim 
\lb( \sH_R^d(\gA,\sG) \xrightarrow{\sH^d\zeta_R}
\sH_{R+1}^d(\gA,\sG) \xrightarrow{\sH^d\zeta_{R+1}}
\sH_{R+2}^d(\gA,\sG) \xrightarrow{\sH^d\zeta_{R+2}}
\cdots\rb)
\eeqn
(See \cite[\S3.F]{Hatcher} or \cite[\S III.9]{Lang} for background on direct
limits.)  We then  have the following 
generalizations of \cite[Theorem 1, \S4]{GePr}:

\Proposition{\label{CA.induced.homotopy.map}}
{
\bthmlist
\item Let $\Sft\subseteq\AZD$ and $\gB\subseteq\BZD$ be $D$-dimensional
SFTs, and let
$\Phi\colon\Sft\rightarrow\gB$ be a subshift homomorphism. 
 Let $\ba\in\gA$ and let
 $\bb:=\Phi(\ba)$.
Then  $\Phi$ induces group homomorphisms 
$\pi_d\Phi\colon\pi_d(\Sft,\ba)\into\pi_d(\gB,\bb)$
for all $d\in\Natur$,
and homomorphisms $\sH_d\Phi\colon\sH_d(\Sft,\sG)\into \sH_d(\gB,\sG)$ and 
$\sH^d\Phi\colon\sH^d(\gB,\sG)\into \sH^d(\gA,\sG)$, for all $d\in\CC{0...D}$.

\item  In particular, if $\Phi\colon\AZD\into\AZD$ is a CA
and $\Phi(\Sft)\subseteq\Sft$, then 
  $\Phi$ induces a group homomorphism 
$\pi_d\Phi\colon\pi_d(\Sft,\ba)\into\pi_d(\Sft,\bb)$
and group endomorphisms
$\sH_d\Phi\colon\sH_d(\Sft,\sG)\into \sH_d(\Sft,\sG)$
and 
$\sH^d\Phi\colon\sH^d(\Sft,\sG)\into \sH^d(\Sft,\sG)$.
\ethmlist
}
\bthmprf {\bf(a)} \
  Suppose $\gA$ has radius $R_A$, and $\gB$ has radius $R_B$,
and let $R:=\max\{R_A,R_B\}$.  If $\Phi$ has radius $q$, 
then for any $r>R$,
 $\Phi$ induces a natural cellular map
$\Phi_*:\bX_{r+q}(\gA)\into\bX_{r}(\gB)$, 
such that, for any $\ba\in\gA_{(r+q)}$ and $\bb:=\Phi(\ba)\in\gB_{(r)}$,
and any $\fz\in\ZD$, the restriction $(\Phi_*)\restr{}:\CornerTile{\fz}{\ba}
\into\CornerTile{\fz}{\bb}$ is a homeomorphism.
Thus, $\Phi_*(x_{r+q}(\ba)) = x_r(\bb)$.
Thus we get a homomorphism 
$\pi_k^r\Phi:\pi_d[\bX_{r+q}(\gA),x_{r+q}(\ba)]\into\pi_d[\bX_{r}(\gB),x_r(\bb)]$ ---i.e.
 a function $\pi_k^r\Phi:\pi^{r+q}_d(\gA,\ba)\into\pi^r_d(\gB,\bb)$.  This yields
a commuting ladder of homomorphisms:
\[
\Array{
\pi^{r+q}_d(\gA,\ba) &\hspace{-1em}\xleftarrow{\pi_d\zeta_{r+q}}\hspace{-1em}&\pi^{r+q+1}_d(\gA,\ba) &\hspace{-1em}\xleftarrow{\pi_d\zeta_{r+q+1}}\hspace{-1em}&\pi^{r+q+2}_d(\gA,\ba) &\hspace{-1em}\xleftarrow{\pi_d\zeta_{r+q+2}}\cdots\hspace{-0.5em} &\mbox{\footnotesize which yields a}\hspace{-0.5em}& \pi_d(\gA,\ba)\\
{\scriptstyle \pi_k^r\Phi}\longdownarrow &      & {\scriptstyle\pi_k^{r+1}\Phi}\longdownarrow &      & 
{\scriptstyle\pi_k^{r+2}\Phi}\longdownarrow && \mbox{\footnotesize homomorphism}& 
{\scriptstyle\pi_k\Phi}\longdownarrow\\
\pi^{r}_d(\gB,\bb) &\hspace{-1em}\xleftarrow[\pi_d\zeta_{r}]{}\hspace{-1em}&\pi^{r+1}_d(\gB,\bb) &\hspace{-1em}\xleftarrow[\pi_d\zeta_{r+1}]{}\hspace{-1em}& \pi^{r+2}_d(\gB,\bb)
&\hspace{-1em}\xleftarrow[\pi_d\zeta_{r+2}]{} \ \cdots\hspace{-0.5em} & \mbox{\footnotesize of colimits:}\hspace{-0.5em} & \pi_d(\gB,\bb). \\}
\]
\ignore{\[
\Array{
\cdots&\into&\pi^{R+q+2}_d(\gA) &\into&\pi^{R+q+1}_d(\gA) &\into&\pi^{R+q}_d(\gA) & \quad&\mbox{\footnotesize which yields a} &\quad& \pi_d(\gA)\\
      &     & \longdownarrow &      & \longdownarrow &      & \longdownarrow && \mbox{\footnotesize homomorphism} && \longdownarrow\\
\cdot&\into&\pi^{r+2}_d(\gB) &\into&\pi^{r+1}_d(\gB) &\into& \pi^{r}_d(\gB)
&& \mbox{\footnotesize of colimits:} && \pi_d(\gB). \\}
\]}
The (co)homology group proof is analogous.  {\bf(b)} follows from {\bf(a)}.
\ethmprf

 For any $\ba\in\tlgA$ and $r\geq R$,
there is a continuous section $\varsigma^r_\ba:\unflawed(\ba)\into\bX_r$ such that $\Pi_r\circ\varsigma^r_\ba = \Id{\unflawed(\ba)}$.
If $k\in\Natur$ and $\ba$  has a range-$r$, codimension-$(k+1)$ defect,
and $0\in\unflawed(\ba)$, then
$\pi_{k}[\unflawed(\ba),0]$ is nontrivial.
Then $\varsigma^r_\ba$ induces a group homomorphism $\pi^r_k\ba:
\pi_{k}[\unflawed(\ba),0]\into\pi_k(\bX_r,x_r(\ba))=\pi_k^r(\gA,\ba)$,
defined by $\pi^r_k\ba(\undgam) :=  \underline{\varsigma^r_\ba\circ\gam}$
for all $\undgam\in\pi_{k}[\unflawed(\ba),0]$.

\example{Recall from Example \ref{X:Conway.Lagarias}(b) that
$\pi_1^1(\Ice)\cong\barpi_1(\sI)\cong\Zahl$, and is generated by the
element $G=\barV-\barA$.
Let $\bi\in\tlIce$ be as in Figure \ref{fig:codim}(D).  
Then $\unflawed[1](\bi)$ is a punctured plane,
so $\pi_1[\unflawed[1](\bi)]\cong\Zahl$.  
If $\zeta$ is a path in $\unflawed[1](\bi)$ that goes once 
counterclockwise around the defect, 
then $\zeta$ generates $\pi_1[\unflawed[1](\bi)]$, and
 $\pi^1_1 \bi (\zeta) = 2\barV + 2\bar\sqe
-2\barA - 2\bar\sqw = 4G$.  Hence, for all $n\in\Zahl$,
$\pi^1_1 \bi (\zeta^n) = 4n G$, so $\pi^1_1 \bi$ is equivalent
to the function $\Zahl\ni n \mapsto 4n \in \Zahl$.}

If $(\sG,+)$ is an abelian group, then 
for any $k\in\CC{0...D}$, 
we likewise get homomorphisms
$\sH^r_k\ba:\sH^r_k[\unflawed(\ba),\sG]\into \sH^r_k(\gA,\sG)$
and $\sH_r^k\ba:\sH_r^k(\gA,\sG)\into\sH_r^k[\unflawed(\ba),\sG]$.

\Theorem{\label{homotopy.comm.square}}
{
Let $\ba\in\tlgA$ have a defect of projective codimension $(k+1)$.
\bthmlist
\item For any abelian group $(\sG,+)$ there are
homomorphisms
$\sH_k\ba\colon \sH_k[\unflawed[\oo](\ba),\sG]\into \sH_k(\gA,\sG)$
and
$\sH^k\ba\colon \sH^k(\gA,\sG)\into  \sH^k[\unflawed[\oo](\ba),\sG]$.

 \item Suppose $\gA$ is basepoint-free in codimension $(k+1)$, and 
let $\omg:\CO{0,\oo}\into\RD$ be a proper base ray.  Then
$\ba$ induces a homomorphism
$\pi_k\ba\colon \pi_k[\unflawed[\oo](\ba),\omg]\into \pi_k(\gA)$.

 \item Let $\Phi\colon\AZD\into\AZD$ be a CA
with $\Phi(\gA)\subseteq\gA$, and let $\Phi(\ba)=\bb$.  Then we have
homomorphisms $\sH_k\iota$ and $\sH^k\iota$ and 
commuting diagrams:
\[  \begin{array}{ccc}
\sH_k[\unflawed[\oo](\ba),\sG] & \xrightarrow{\ \sH_k\iota \ } & \sH_k[\unflawed[\oo](\bb),\sG]\\
{\scriptstyle \sH_k\ba} \longdownarrow & & \longdownarrow {\scriptstyle \sH_k\bb} \\
\sH_k(\gA,\sG) & \xrightarrow{\sH_k\Phi} & \sH_k(\gA,\sG)\\
\end{array}
\qquad\mbox{and}\qquad
 \begin{array}{ccc}
\sH^k[\unflawed[\oo](\ba),\sG] & \xleftarrow{\ \sH^k\iota \ } & \sH^k[\unflawed[\oo](\bb),\sG]\\
{\scriptstyle \sH^k\ba} \longuparrow & & \longuparrow {\scriptstyle \sH^k\bb} \\
\sH^k(\gA,\sG) & \xleftarrow{\sH^k\Phi} & \sH^k(\gA,\sG)\\
\end{array}\]
Assuming the hypothesis of {\bf(b)}, we have a 
homomorphism  $\pi_k\iota$ and a commuting diagram:
\beqn
\label{homotopy.comm.square.e1}
  \begin{array}{ccc}
\pi_k[\unflawed[\oo](\ba),\omg] & \xrightarrow{\ \pi_k\iota \ } & \pi_k[\unflawed[\oo](\bb),\omg]\\
{\scriptstyle \pi_k\ba} \longdownarrow & & \longdownarrow {\scriptstyle \pi_k\bb} \\
\pi_k(\gA) & \xrightarrow{\pi_k\Phi} & \pi_k(\gA)\\
\end{array}
\eeqn
\ethmlist
}
\bthmprf We will prove the statements for $\pi_k$. 
 The (co)homological versions are analogous.

{\bf(b)} \ Let $R$ be the radius of $\gA$, and
recall that $\unflawed[R](\ba)\supseteq\unflawed[R+1](\ba)\supseteq\unflawed[R+2](\ba)\supseteq\cdots$.  For each $r\geq R$, 
the inclusion map
$\alp_r:\unflawed[r+1](\ba)\hookrightarrow\unflawed(\ba)$
induces a canonical homomorphism $\alp^*_r:\pi_k[\unflawed[r+1](\ba),\omg]
\into\pi_k[\unflawed(\ba),\omg]$.  
Let $w_r\in\omg\CO{0,\oo}\intsct\unflawed(\ba)$
(so  $\pi_k[\unflawed(\ba),\omg] \cong 
\pi_k[\unflawed(\ba),w_r]$ canonically), let
$a_r := \varsigma^r_\ba(w_r)\in\bX_r$, and define
$\pi^r_k(\gA,a_r):=\pi_k(\bX_r,a_r)$.  There is a canonical
homomorphism $\zet^*_r:\pi^{r+1}_k(\gA,a_{r+1})
\into\pi^r_k(\gA,a_r)$, because $\gA$ is basepoint-free.
This yields a commuting ladder of homomorphisms,
which defines homomorphism of inverse limits:
\beqn
\label{homotopy.comm.square.e2}
\Array{
\pi_d[\unflawed[R](\ba),\omg] &\hspace{-0.7em}\xleftarrow{\alp^*_{R}}\hspace{-0.7em}&\pi_d[\unflawed[R+1](\ba),\omg] &\hspace{-0.7em}\xleftarrow{\alp^*_{R+1}}\hspace{-0.7em}&\pi_d[\unflawed[R+2](\ba),\omg] &\hspace{-0.7em}\xleftarrow{\alp^*_{R+2}}\cdots\hspace{-0.7em} & & \pi_d[\unflawed[\oo](\ba),\omg]\\
{\scriptstyle \pi_k^R\ba}\longdownarrow &      & {\scriptstyle \pi_k^{R+1}\ba}\longdownarrow &      & {\scriptstyle \pi_k^{R+2}\ba}\longdownarrow &&  & {\scriptstyle \pi_k\ba}\longdownarrow\\
\pi^{R}_d(\gA,a_R) &\hspace{-0.7em}\xleftarrow[\zet^*_{R}]{}\hspace{-0.7em}&\pi^{R+1}_d(\gA,a_{R+1}) &\hspace{-0.7em}\xleftarrow[\zet^*_{R+1}]{}\hspace{-0.7em}& \pi^{R+2}_d(\gA,a_{R+2})
&\hspace{-0.7em}\xleftarrow[\zet^*_{R+2}]{}\cdots\hspace{-0.7em} &  & \pi_d(\gA). \\}
\eeqn
{\bf(c)} 
For any $r\geq R$, the inclusion map
$\bet_r:\unflawed[r+1](\bb)\hookrightarrow\unflawed(\bb)$
induces a canonical homomorphism $\bet^*_r:\pi_k[\unflawed[r+1](\bb),\omg]
\into\pi_k[\unflawed(\bb),\omg]$.  Assume that $w_r\in\unflawed(\bb)$,
let $b_r := \varsigma^r_\ba(w_r)\in\bX_r$, and define
$\pi^r_k(\gB,b_r):=\pi_k(\bX_r,b_r)$.
This yields a commuting ladder like
eqn.(\ref{homotopy.comm.square.e2}), only with $\bb$ instead
of $\ba$, $\beta_r^*$ instead of $\alp_r^*$, 
 and $b_r$ instead of $a_r$. 
Suppose $\Phi$ has radius $q>0$.  If $\Phi_*:\bX_{r+q}\into\bX_r$ is the
cellular map induced by $\Phi$, then $\Phi_*(a_{q+r})=b_r$ [because
$\Phi(\ba_{\dB(q+r)})=\bb_{\dB(r)}$].
This yields homomorphisms
 $\pi^r_d\Phi:\pi^{r+q}_k(\gA,a_{r+q}) \into \pi^r_k(\gA,b_r)$,
for all $r\geq R$, as in the proof of Proposition
\ref{CA.induced.homotopy.map}.

For any $r\geq R$, Proposition \ref{defect.dimension}(b) says
$\unflawed[r+q](\ba)\subseteq\unflawed(\bb)$. 
The inclusion map
$\iota_r:\unflawed[r+q](\ba)\hookrightarrow\unflawed(\bb)$
induces a (canonical) homomorphism $\iota^*_r:\pi_k[\unflawed[r+q](\ba),\omg]
\into \pi_k[\unflawed(\bb),\omg]$, and we have a commuting diagram
\beqn
\label{homotopy.comm.square.e3}
  \begin{array}{ccc}
\pi_k[\unflawed[r+q](\ba),\omg] & \xrightarrow{\iota_r^*} & \pi_k[\unflawed(\bb),\omg]\\
{\scriptstyle \pi_k^{r+q}\ba}\longdownarrow & & \longdownarrow {\scriptstyle \pi_k^{r}\bb} \\
\pi^{r+q}_k(\gA,a_{r+q}) & \xrightarrow{\pi^r_d\Phi} & \pi^r_k(\gA,b_r)\\
\end{array}
\eeqn
Combining the commuting ladders (\ref{homotopy.comm.square.e2})
for $\ba$ and $\bb$, along with copies of
the square (\ref{homotopy.comm.square.e3}) for each $r\in\Natur$, we obtain
a `commuting girder' of homomorphisms:

\centerline{
\psfrag{A1}[][]{$\pi_d[\unflawed[q+R](\ba)]$}
\psfrag{B1}[][]{$\pi_d[\unflawed[R](\bb)] $}
\psfrag{C1}[][]{$\pi^{q+R}_d(\gA)$}
\psfrag{D1}[][]{$\pi^{R}_d(\gA)$}
\psfrag{f1}[][]{$\alp^*_{q+R}$}
\psfrag{g1}[][]{$\bet^*_R$}
\psfrag{h1}[][]{$\zeta^*_{q+R}$}
\psfrag{i1}[][]{$\zeta^*_{R}$}
\psfrag{j1}[][]{$\iota^*_R$}
\psfrag{k1}[][]{$\pi^R_d\Phi $}
\psfrag{l1}[][]{$\pi^R_d\bb$}
\psfrag{m1}[][]{$\pi^{q+R}_d\ba$}
\psfrag{A2}[][]{$\pi_d[\unflawed[q+R+1](\ba)]$}
\psfrag{B2}[][]{$\pi_d[\unflawed[R+1](\bb)] $}
\psfrag{C2}[][]{$\pi^{q+R+1}_d(\gA)$}
\psfrag{D2}[][]{$\pi^{R+1}_d(\gA)$}
\psfrag{f2}[][]{$\alp^*_{q+R+1}$}
\psfrag{g2}[][]{$\bet^*_{R+1}$}
\psfrag{h2}[][]{$\zeta^*_{q+R+1}$}
\psfrag{i2}[][]{$\zeta^*_{R+1}$}
\psfrag{j2}[][]{$\iota^*_{R+1}$}
\psfrag{k2}[][]{$\pi^{R+1}_d\Phi $}
\psfrag{l2}[][]{$\pi^{R+1}_d\bb$}
\psfrag{m2}[][]{$\pi^{q+R+1}_d\ba$}
\psfrag{A3}[][]{$\pi_d[\unflawed[q+R+2](\ba)]$}
\psfrag{B3}[][]{$\pi_d[\unflawed[R+2](\bb)] $}
\psfrag{C3}[][]{$\pi^{q+R+2}_d(\gA)$}
\psfrag{D3}[][]{$\pi^{R+2}_d(\gA)$}
\psfrag{f3}[][]{$\alp^*_{q+R+2}$}
\psfrag{g3}[][]{$\bet^*_{R+2}$}
\psfrag{h3}[][]{$\zeta^*_{q+R+2}$}
\psfrag{i3}[][]{$\zeta^*_{R+2}$}
\psfrag{j3}[][]{$\iota^*_{R+2}$}
\psfrag{k3}[][]{$\pi^{R+2}_d\Phi $}
\psfrag{l3}[][]{$\pi^{R+2}_d\bb$}
\psfrag{m3}[][]{$\pi^{q+R+2}_d\ba$}
\psfrag{Aoo}[][]{$\pi_d[\unflawed[\oo](\ba)]$}
\psfrag{Boo}[][]{$\pi_d[\unflawed[\oo](\bb)] $}
\psfrag{Coo}[][]{$\pi_d(\gA)$}
\psfrag{Doo}[][]{$\pi_d(\gA)$}
\psfrag{joo}[][]{$\iota^*$}
\psfrag{koo}[][]{$\pi_d\Phi $}
\psfrag{loo}[][]{$\pi_d\bb$}
\psfrag{moo}[][]{$\pi_d\ba$}
\includegraphics[angle=-90,scale=0.75]{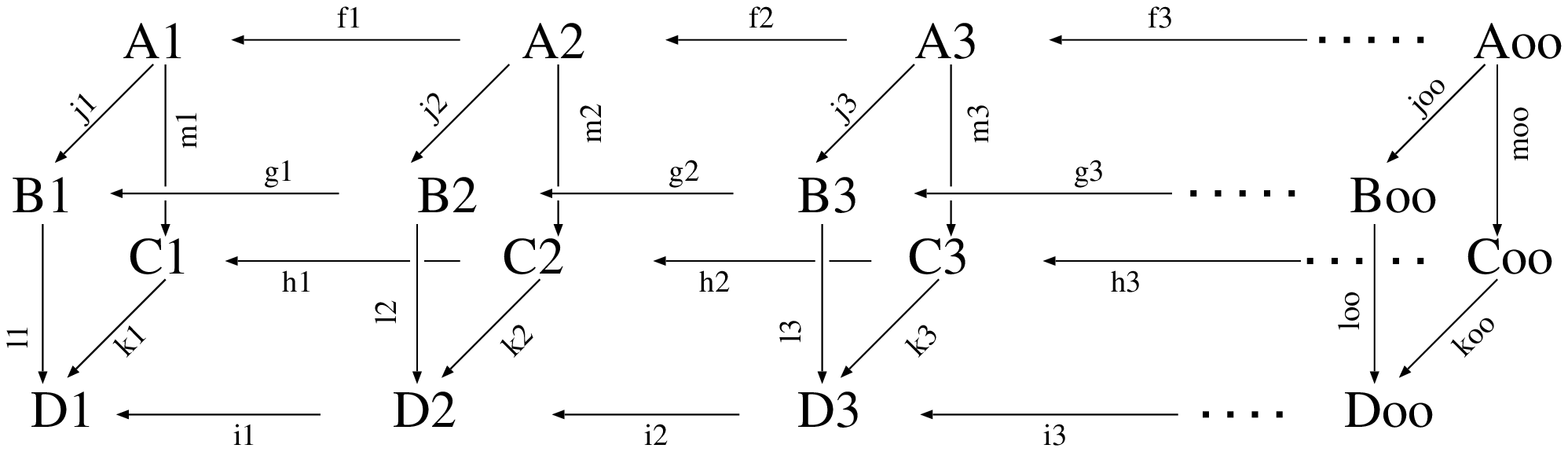}}

which yields the commuting
square (\ref{homotopy.comm.square.e1}) of colimit homomorphisms.
\ethmprf

  We call $\pi_k\ba$ (resp. $\sH_k\ba$ or $\sH^k\ba$) 
 the {\dfn $k$th homotopy} (resp. {\dfn (co)homology})
{\dfn signature} of $\ba$; if it is nontrivial,
we say $\ba$ has a {\dfn homotopy} (resp. ({\dfn co}){\dfn homology}) {\dfn defect} of codimension $(k+1)$.  The next result is analogous to
Proposition \ref{persistent.pole}(b):

\Corollary{\label{persistent.homotopy.defect}}
{
 Let $\Phi:\AZD\into\AZD$ be a CA and with $\Phi(\gA)\subseteq\gA$.
\bthmlist  
\item Suppose $\gA$ is basepoint-free in codimension $(k+1)$.
If $\pi_k\Phi:\pi_k(\gA)\into\pi_k(\gA)$ is injective,
then every homotopy defect of codimension $(k+1)$ is $\Phi$-persistent.

 \item If $\sH_k\Phi:\sH_k(\gA,\sG)\into\sH_k(\gA,\sG)$ is injective,
then every homology defect of codimension $(k+1)$ is $\Phi$-persistent.

 \item If $\sH^k\Phi:\sH^k(\gA,\sG)\into\sH^k(\gA,\sG)$ is surjective,
then every cohomology defect of codimension $(k+1)$ is $\Phi$-persistent.
\ethmlist
}
\bthmprf {\bf(a)} If $\pi_k\ba$ is nontrivial, and $\pi_k\Phi$ is a monomorphism,
then $\pi_k\Phi\circ\pi_k\ba$ must also be nontrivial.
Then diagram (\ref{homotopy.comm.square.e1}) says that
 $\pi_k\bb\circ \pi_k\iota$
must be nontrivial, hence $\pi_k\bb$ is nontrivial, hence
$\bb$ has a homotopy defect in codimension $(k+1)$.  The proofs of {\bf(b,c)}
are similar.
\ethmprf

\Remark{Using the machinery of Appendix \ref{S:homo.alg} (below), we could
compute $\sH^k(\gA,\sG)$ as follows. 
For each $r\geq R$,
the function  $\zeta_r:\bX_{r+1}\into\bX_r$ is surjective,
so the contravariant homomorphism 
$\sC^d\zeta_r\colon\sC^d(\bX_r,\sG)\rightarrow\sC^d(\bX_{r+1},\sG)$
(defined $C\mapsto C\circ \zeta_r$) is injective.  
Thus, we can identify $\sC^d(\bX_r,\sG)$ as a subgroup of
$\sC^d(\bX_{r+1},\sG)$ in a natural way.  This yields an ascending chain
$\sC^d(\bX_R,\sG)\subseteq \sC^d(\bX_{R+1},\sG) \subseteq \sC^d(\bX_{R+2},\sG)
\subseteq \cdots$.  Let $\sC^d_\oo := \Union_{r=R}^\oo \sC^d(\bX_r,\sG)$.
The system of coboundary maps 
$\{\partial_d:\sC^d(\bX_r,\sG) \into \sC^{d+1}(\bX_r,\sG)\}_{r=R}^\oo$
then defines a 
coboundary map $\partial_d:\sC_\oo^d \into \sC_\oo^{d+1}$.  This
yields a  chain complex $\bC_\oo := \{\sC_\oo^d,\partial_d\}_{d=0}^\oo$,
which is the colimit of the system of chain complexes
$(\bC_R\stackrel{\bzet_R}{\into} \bC_{R+1}\stackrel{\bzet_{R+1}}{\into}\cdots)$
(where  $\bC_r := \{\sC^d(\bX_r,\sG),\partial_d\}_{d=0}^\oo$
for each $r\geq R$).  Lemma \ref{cohomology.limit.lemma} 
(in Appendix \ref{S:homo.alg}) then implies that
 $\sH^k(\gA,\sG) = \sH^k(\bC_\oo)$.

Unfortunately, this argument doesn't dualize  to homology or homotopy
groups.
}

\section{Equivariant vs. Invariant Cohomology \label{S:cohomology}}

By relating the dynamical cohomology of \S\ref{S:cocycle.intro} to the
tiling cohomology of \S\ref{S:projective.homotopy}, we can
generalize the results of
\S\ref{S:pole} to defects of higher codimensions.
The main results of this section are Theorems
\ref{two.cocycles.are.one}, \ref{persistent.dpole},
and \ref{three.cohomology.theorem}, and Proposition \ref{projective.d.pole.essential}.

\subsection{Equivariant Cohomology of Subshifts: \label{S:equiv.cohom}}
Let $\bY=(\bY^0,\ldots,\bY^D)$ be as in
\S\ref{S:canonical.cell.cplx}.
Let
$(\sG,+)$ be an abelian topological group, and for all $d\in\CC{1...D}$,
give $\sC^d(\bY,\sG)
\cong \sG^{\bY^d}$ the Tychonoff product topology.
For any $\fz\in\ZD$, let $\Upsilon^\fz_d:\bY^d\into\bY^d$ be the 
obvious translation function.  For example,
we define $\Upsilon^\fz_0:\bY^0\into\bY^0$  by $\Upsilon^\fz_0(\dot\fx):=\dot\fy$
where $\fy:=\fx+\fz$,
and we define $\Upsilon^\fz_1:\bY^1\into\bY^1$  by 
$\Upsilon^\fz_1(\onecell{\fx}{\fy}):=\onecell{\fx'}{\fy'}$ where $\fx':=\fx+\fz$ and 
$\fy':=\fy+\fz$.
We extend this to a function $\Upsilon^\fz_d:\Zahl[\bY^d]\into\Zahl[\bY^d]$
by linearity.   Clearly, $\partial_{d}\circ\Upsilon_d^\fz
\ = \ \Upsilon_{d-1}^\fz\circ\partial_d$.
We define $\Upsilon^d_\fz:\sC^d\into\sC^d$ by 
$\Upsilon^d_\fz(C)(\zeta)=C\circ\Upsilon^\fz_d(\zeta)$ for any $C\in\sC^d$
and $\zeta\in\Zahl[\bY^d]$. Thus,
$\del_{d}\circ\Upsilon^d_\fz \ = \ \Upsilon^{d+1}_\fz\circ\del_d$.

Let $\gA\subset\AZD$ be a subshift.  A 
($d$-dimensional, $\sG$-valued, continuous)
{\dfn equivariant cochain} is a continuous function $C:\gA\into\sC^d$ 
so that, for any $\ba\in\gA$ and $\fz\in\ZD$, 
$C(\shift{\fz}(\ba)) = \Upsilon^d_\fz\circ C(\ba)$.  Equivalently,
an equivariant cochain is a continuous 
function $C:\Zahl[\bY^d]\x\gA\into\sG$ such that:
\beqn
\label{cochain.eqn}
\mbox{For all $\zeta\in\Zahl[\bY^d]$, \  $\ba\in\gA$, \ and $\fz\in\ZD$,}\qquad
C(\zeta,\shift{\fz}(\ba))\quad=\quad C(\Upsilon^\fz_d(\zeta),\ba).
\eeqn

\example{\label{X:eq.cochain}
(a) ($d=0$) \   An equivariant zero-cochain is equivalent to a 
continuous function 
$B:\ZD\x\gA\into\sG$ such that $B(\fx,\shift{\fy}(\ba))=B(\fx+\fy,\ba)$
for any $\fx,\fy\in\ZD$ and $\ba\in\gA$.

(b)  ($d=1$) \ An equivariant 1-cochain is equivalent to a continuous function 
$C:\bY^1\x\gA\into\sG$ such that, for any $\ba\in\gA$ and
$\fx,\fy,\fz\in\ZD$, if
$\fx':=\fx+\fz$ and $\fy':=\fy+\fz$, then
$C(\onecell{\fx}{\fy};\shift{\fz}(\ba)) \ = \ 
C(\onecell{\fx'}{\fy'};\ba)$.
}

 Let $\Ceq^d(\gA,\sG)$ be the abelian group of
all equivariant cochains (with pointwise addition).  We define the {\dfn coboundary}
homomorphism $\del_d:\Ceq^d(\gA,\sG)\into\Ceq^{d+1}(\gA,\sG)$ by applying
the cubic coboundary map (from \S\ref{S:cubic.cohom}) pointwise. 
That is, for any
 $\zeta\in\Zahl[\bY^{d+1}]$
and $\ba\in\gA$,  let $\del_d C(\zeta,\ba) := C(\partial_{d+1} \zeta,\ba)$.
We say that $C$ is an {\dfn equivariant cocycle} if $\del_d C=0$, and
we say that $C$ is an {\dfn equivariant coboundary} if
$C=\del_{d-1} b$ for some $b\in \Ceq^{d-1}(\gA,\sG)$.

\example{\label{X:ice.cube.cocycle}
 Let $\gQ\subset\sQ^{\ZD[3]}$ be the `ice-cube' shift from
Example \ref{X:codim.one}(c).  Define $C:\bY^2\x\gQ\into\Zahl$
as follows. Any two-cell $\twocell{\fw}{\fx}{\fy}{\fz}$ 
is an oriented square frame, with a unique normal vector $\vV$
defined by the right-hand
rule.  For any $\bq\in\gQ$ (seen as a `ball-and-pin' assembly),
there is a unique `pin' passing through $\twocell{\fw}{\fx}{\fy}{\fz}$,
from one the two adjoining cubes into the other. 
Define $C(\twocell{\fw}{\fx}{\fy}{\fz}, \bq) = +1$ if this pin 
is parallel to  $\vV$, and $-1$ if it is antiparallel.
Then $C$ is a 2-dimensional equivariant cocycle.
}
Let $\Zeq^d(\gA,\sG):=\ker(\del_d)$ and
$\Beq^d(\gA,\sG):=\image{\del_{d-1}}$ be the subgroups of cocycles and
coboundaries respectively; then $\Beq^d(\gA,\sG) \subseteq
\Zeq^d(\gA,\sG)$, and the quotient $\Heq^d(\gA,\sG):=
\Zeq^d(\gA,\sG)/\Beq^d(\gA,\sG)$ is the {\dfn $d$th equivariant
cohomology group} of $\gA$ (with coefficients in $\sG$).  Two cocycles
$C_1$ and $C_2$ are {\dfn cohomologous} ($C_1\approx C_2$) if they
project to the same coset in $\Heq^d(\gA,\sG)$ ---i.e. if $C_1 = C_2 +
\del_{d-1} B$ for some $B\in\Ceq^{d-1}(\gA,\sG)$.

\example{\label{X:eq.cobound}
If  $C_1,C_2\in\Zeq^1(\gA,\sG)$, then
 $C_1 \approx C_2$ iff there is an equivariant zero-cochain
$B:\bY^0\x\gA\into\sG$ such that, 
for any $(\onecell{\fx}{\fy})\in\bY^1$ and $\ba\in\gA$, \ \ 
$C_2[(\onecell{\fx}{\fy}),\ba] \ = \ 
B(\dot\fy,\ba)+ C_1[(\onecell{\fx}{\fy}),\ba] -B(\dot\fx,\ba)$.
}
Equivariant cocycles are the natural generalization of the `dynamical'
cocycles from \S\ref{S:cocycle.intro}.  To see this, recall that a
(continuous, $\sG$-valued) {\dfn dynamical cocycle} on $\gA$ is a
continuous function $\tlC:\ZD\x\gA\into\sG$ satisfying the cocycle
equation (\ref{cocycle.eqn}), which, in additive notation, reads:
$\tlC(\fy+\fz,\ba) \ = \ \tlC[\fy,\shift{\fz}(\ba)] +
\tlC(\fz,\ba)$.
Two dynamical cocycles $\tlC_1,\tlC_2$ are {\dfn cohomologous} 
if there is some $\tlB:\gA\into\sG$ such that $\tlC_2(\fz,\ba)
\ = \ \tlB(\shift{\fz}(\ba)) + \tlC_1(\fz,\ba) - \tlB(\ba)$ for any
$\fz\in\ZD$ and $\ba\in\gA$.   Let $\Zdyn^1(\gA,\sG)$
be the additive group of dynamical cocycles, and 
let $\Hdyn^1(\gA,\sG)$ be the dynamical cohomology group.
The main result of \S\ref{S:equiv.cohom} is:

\Theorem{\label{two.cocycles.are.one}}
{
  Let $\gA\subset\AZD$ be a subshift and 
let $(\sG,+)$ be an abelian group.  
  There are canonical isomorphisms $\Zeq^1(\gA,\sG)\cong \Zdyn^1(\gA,\sG)$,
and 
$\Heq^1(\gA,\sG)\cong \Hdyn^1(\gA,\sG)$.
}

To prove Theorem \ref{two.cocycles.are.one},
we will introduce another intermediate notion of one-dimensional cocycle.
A (continuous, $\sG$-valued) {\dfn two-point cocycle} on $\gA$
is a continuous function $\ddot{C}:\ZD\x\ZD\x\gA\into\sG$ such that,
for any fixed $\ba\in\gA$,

  {\bf(\"{C}1)} \  The function $\ddot{C}(\blankspace,\blankspace;\ba):
\ZD\x\ZD\into\sG$ is a two-point cocycle in the sense of
eqn.{\rm(\ref{X:cubic.cocycle.eqn})} in Example \ref{X:cubic.cocycle}.

  {\bf(\"{C}2)} \ For any $\fw,\fy,\fz\in\ZD$,
 we have  $\ddot{C}[\fw,\fy;\shift{\fz}(\ba)] =  \ddot{C}(\fw+\fz, \fy+\fz;\ba)$.

\noindent Let $\Ztwo^1(\gA,\sG)$ be the group of two-point cocycles.
Then Theorem \ref{two.cocycles.are.one} follows from:

\Lemma{\label{three.cocycles}}
{Let $\gA\subset\AZD$ be a subshift, and let $(\sG,+)$ be an abelian group.
\bthmlist
  \item
   Let $C\in\Zeq^1(\gA,\sG)$.
For any $\fy,\fz\in\ZD$ and $\ba\in\gA$, define
$\ddot{C}(\fy,\fz;\ba):=C(\zeta,\ba)$,
where $\zeta\in\Zahl[\bY^1]$ is any chain such that $\partial_1(\zeta)=
\dot\fz-\dot\fy$.  Then $\ddot{C}$ well-defined, and
$\ddot{C}\in\Ztwo^1(\gA,\sG)$.

  \item Let $\ddot{C}\in\Ztwo^1(\gA,\sG)$.
Define $\tlC(\fz,\ba):=\ddot{C}(0,\fz;\ba)$ for all
 $\fz\in\ZD$ and $\ba\in\gA$.
Then $\tlC\in\Zdyn^1(\gA,\sG)$.

  \item Let $\tlC\in\Zdyn^1(\gA,\sG)$.
Define $C(\onecell{\fx}{\fy};\ba) :=  \tlC[(\fy-\fx), \shift{\fx}(\ba)]$
for all $(\onecell{\fx}{\fy})\in\bY^1$ and $\ba\in\gA$,
and extend to $C:\Zahl[\bY^1]\x\gA\into\sG$ by linearity.
Then $C\in\Zeq^1(\gA,\sG)$.

  \item  Any two of the following statements imply the third:
\qquad  {\rm[i]} $\ddot{C}$ comes from $C$ via {\rm(a)}.
\qquad {\rm[ii]} $\tlC$ comes from $\ddot{C}$ via {\rm(b)}.
\qquad {\rm[iii]} $C$ comes from $\tlC$ via {\rm(c)}.

  \item If $C_1,C_2\in\Zeq^1(\gA,\sG)$,
and $\tlC_1,\tlC_2\in\Zdyn^1(\gA,\sG)$ are
related to $C_1$ and $C_2$ as in {\rm(c)}, then $(C_1\approx C_2)
\iff (\tlC_1\approx \tlC_2)$.
\ethmlist
}
\bthmprf
{\bf(a)} \ Reason as in Example \ref{X:cubic.cocycle} to get {\bf(\"{C}1)}.
Then use eqn.(\ref{cochain.eqn}) to get {\bf(\"{C}2)}.

{\bf(b)} \
$\tlC(\fy+\fz,\ba)
\ \ \eeequals{(\ddagger)} \ \ \ddot{C}(0,\fy+\fz;\ba)
\ \eeequals{(*)} \  \ddot{C}(\fz,\fy+\fz;\ba) + \ddot{C}(0,\fz;\ba)
\ \eeequals{(\dagger)} \  
\ddot{C}[0,\fy;\shift{\fz}(\ba)] + \ddot{C}(0,\fz;\ba)$
$ \ \ \eeequals{(\ddagger)} \ \  \tlC[\fy,\shift{\fz}(\ba)] + \tlC(\fz,\ba)$.
Here $(\ddagger)$ is the definition of $\tlC$ in part (b).
$(*)$ is by {\bf(\"{C}1)} and eqn.(\ref{X:cubic.cocycle.eqn}b),
while $(\dagger)$ is by  {\bf(\"{C}2)}.

{\bf(c)}  Let $\twocell{\fs}{\ft}{\fu}{\fv}$ be a two-cell
and fix $\bb\in\gA$.
To show that $C(\bb)$ is a one-cycle, it suffices to check that
$C(\onecell{\fs}{\ft};\bb) + C(\onecell{\ft}{\fu};\bb)
\ = \ 
C(\onecell{\fs}{\fv};\bb) + C(\onecell{\fv}{\fu};\bb)$ (see Example \ref{X:cubic.cocycle}).  But
\beq
\lefteqn{C(\onecell{\fs}{\ft};\bb) + C(\onecell{\ft}{\fu};\bb)}\\
&\eeequals{(\ddagger)}&
\tlC[\ft-\fs;\shift{\fs}(\bb)] + \tlC[\fu-\ft;\shift{\ft}(\bb)] 
\quad=\quad
\tlC[\ft-\fs;\shift{\fs}(\bb)] + \tlC[\fu-\ft;\shift{\ft-\fs}(\shift{\fs}(\bb))]\\&\eeequals{(*)}&
\tlC[\fu-\fs;\shift{\fs}(\bb)]
\quad\eeequals{(\dagger)}\quad
\tlC[\fv-\fs;\shift{\fs}(\bb)] + \tlC[\fu-\fv;\shift{\fv-\fs}(\shift{\fs}(\bb))]\\&=&
\tlC[\fv-\fs;\shift{\fs}(\bb)] + \tlC[\fu-\fv;\shift{\fv}(\bb)] 
\quad\eeequals{(\ddagger)}\quad
C(\onecell{\fs}{\fv};\bb) + C(\onecell{\fv}{\fu};\bb),\quad\mbox{as desired.}
\eeq
Here, $(\ddagger)$ is the definition of $C$ in part (c).
 $(*)$ is by  eqn.(\ref{cocycle.eqn})
with  $\ba:=\shift{\fs}(\bb)$,  $\fz:=\ft-\fs$ and $\fy := \fu-\ft$,
while
$(\dagger)$ is by  eqn.(\ref{cocycle.eqn}) with
 $\ba:=\shift{\fs}(\bb)$,  $\fz:=\fv-\fs$ and $\fy := \fu-\fv$.

To see that $C$ is equivariant, let
 $\ba\in\gA$ and $\fx,\fy,\fz\in\ZD$. If $\fx':=\fx+\fz$ and 
$\fy':=\fy+\fz$, then\hfill
$
C[\onecell{\fx}{\fy};\shift{\fz}(\ba)]
\hfill \ \eeequals{(\dagger)} \ \hfill 
\tlC[\fy-\fx;\shift{\fx}(\shift{\fz}(\ba))] 
\hfill = \hfill \tlC[(\fy+\fz)-(\fx+\fz);\shift{\fx+\fz}(\ba)] \hfill
= \hfill$
$\tlC[\fy'-\fx';\shift{\fx'}(\ba)] 
\ \ \eeequals{(\dagger)} \ \ C(\onecell{\fx'}{\fy'};\ba)$,
as desired.  Here $(\dagger)$ is the definition of $C$ in part (c).

{\bf(d)}  Straightforward calculation.

{\bf(e)} 
Given any continuous function $\tlB:\gA\into\sG$, 
define $B:\bY^0\x\gA\into\sG$ by 
$B(\dot\fz,\ba):=\tlB(\shift{\fz}(\ba))$ for all $\fz\in\ZD$ and
$\ba\in\gA$.  Then $B$ is an equivariant zero-cochain, i.e.
$B[\dot\fx,\shift{\fy}(\ba)] = B(\dot\fz,\ba)$ (where $\fz=\fx+\fy$).
Conversely, any equivariant zero-cochain arises in this manner:
given $B:\bY^0\x\gA\into\sG$, define $\tlB(\ba):=B(\dot0,\ba)$;
then $B(\dot\fz,\ba)=\tlB(\shift{\fz}(\ba))$ for all $\fz\in\ZD$.
If $B$ and $\tlB$ are related in this way, and $C_n$ is related to
$\tlC_n$ as in part (c) for $n=1,2$, then
it is easy to check that $C_1$ is cohomologous to $C_2$ via $B$
iff  $\tlC_1$ is cohomologous to $\tlC_2$ via $\tlB$.
\ethmprf

\example{\label{X:higher.cocycle}
 Let $(\sG,+)$ be an abelian group, and
let $\tlC:\ZD\x\gA\into\sG$ be a dynamical cocycle,
[e.g. Examples \ref{X:cocycle}(c,d,e)], and define $C\in\Zeq^1(\gA,\sG)$ 
from $\tlC$ as in Lemma \ref{three.cocycles}(c).

Let $\tlzet:=(\fz_0\leadsto\fz_1\leadsto\cdots\leadsto\fz_N)$ be
a trail in $\ZD$,
and let $\zeta:=\sum_{n=0}^{N-1} (\onecell{\fz_n}{\fz_{n+1}})$
be the 1-chain `representing' $\tlzet$. 
For all $n\in\CC{1...N}$, let $\fz'_n := \fz_{n}-\fz_{n-1}$,
Then
\beqn
\label{chain.cocycle}
C(\zeta,\ba)\quad =
\quad \sum_{n=1}^{N} C(\onecell{\fz_{n-1}}{\fz_{n}}; \ \ba)
\quad=\quad\sum_{n=1}^{N} \tlC[\fz'_n,\shift{\fz_{n-1}}(\ba)]
\quad\eeequals{(*)}\quad \tlC(\tlzet,\ba),
\eeqn
where $(*)$ is the additive version of eqn.(\ref{trail.cocycle}).}

To generalize the notion of `locally determined' cocycles from
\S\ref{S:cocycle.intro}, we need some notation.
For any $i\in\CC{1...D}$, recall (from \S\ref{S:canonical.cell.cplx}) that $\partial_i^{-} \ \CornerTile{0}{\ }$ is a $(D-1)$-dimensional face of the cube $\CornerTile{0}{\ }$ which is orthogonal to
the $i$th axis.  Also, for any $r\in\Natur$, recall (from \S\ref{S:Wang.representation}) that  $\partial_i^-\, \dB(r) := \CC{-r...r}^{i-1} \x \CO{-r...r} \x
\CC{-r...r}^{D-i-1}$.
For any $d\in\CC{0...D}$,
let $\bY^d_0$ be the set of all $d$-cells in $\bY^d$ contained
in the half-open unit cube $\CO{0,1}^D$.  Let $k:=D-d$.
If $x\in\bY_0^d$, then  $x=
\partial_{i_1}^{-}\, \CornerTile{0}{\ } \, \intsct \,
\partial_{i_2}^{-}\, \CornerTile{0}{\ } \, \intsct 
\cdots\intsct \, \partial_{i_k}^{-}\,\CornerTile{0}{\ }$, for some 
$\{i_1,\ldots,i_k\}\subset\CC{1...D}$ (see \S\ref{S:canonical.cell.cplx}).
We then define $\partial_x \, \dB(r) := 
\partial_{i_1}^{-}\, \dB(r)\,\intsct\,
\partial_{i_2}^{-}\, \dB(r)\,\intsct
\cdots\intsct\,\partial_{i_k}^{-}\,\dB(r)$;  this is the $d$-dimensional
`face' of $\dB(r)$ corresponding to $x$.  For any $\fz\in\ZD$, 
we  define
$\partial_x\,\dB(\fz,r):=\partial_x \, \dB(r) + \fz$.

For any $y\in\bY^d$, there is a unique $x\in\bY^d_0$ and
$\fz\in\ZD$ so that $y=\Upsilon^\fz_d(x)$.
An equivariant cocycle $C\in\Zeq^d(\gA,\sG)$ is {\dfn locally determined} 
with {\dfn radius} $r>0$ if, for each $x\in\bY_0^d$,
there is some {\dfn local rule}  $c_{x}:\gA_{\partial_x \dB(r)}\into\sG$
such that, for any $y\in\bY^d$ and $\ba\in\gA_{(r)}$,
if $y = \Upsilon^\fz_d(x)$ for $\fz\in\ZD$, then 
$C(y,\ba) \ := \ c_{x}(\ba_{\partial_x\dB(\fz,r)})$.  

\example{(a) Let $d=1$. If $C\in\Zeq^1$ and $\tlC\in\Zdyn^1$ are related as in
Lemma \ref{three.cocycles}(c), then
$C$ is locally determined iff $\tlC$ is locally determined in the
sense of \S\ref{S:cocycle.intro}.

(b) The cocycle in Example \ref{X:ice.cube.cocycle} is locally determined, with radius $0$.

(c) If $\sG$ is discrete, then every continuous
equivariant cocycle is locally determined.}

Fix $d\in\CC{0..D}$.  If $C\in\Zeq^d(\gA,\sG)$ is locally determined
with radius $r$, and $\ba\in\tlgA$, then $C(y,\ba)$ is well-defined
for any $y\in\bY^d\intsct\unflawed(\ba)$.  A {\dfn $d$-cycle} is a
$d$-chain $\zeta\in\Zahl[\bY^d]$ such that $\partial_d(\zeta)=0$; let
$\sZ_d=\sZ_d(\bY;\Zahl)$ denote the group of {\dfn $d$-cycles}.  We
say that $\ba$ has a {\dfn $C$-pole of range $r$} if there is some
$d$-cycle $\zeta\in\sZ_d[\unflawed(\ba)]$ such that $C(\zeta,\ba)\neq
0$.  We say that $\ba$ has a {\dfn projective $C$-pole} if
$\ba$ has a {\dfn $C$-pole of range $r$} for all large enough $r\in\Natur$.
We say that $\ba$ has a {\dfn projective $(\sG,d)$-pole}
if $\ba$ has a projective $C$-pole for some $C\in\Zeq^d(\gA,\sG)$.
We will call this a ``$d$-pole'' (resp. ``$\sG$-pole'') if $\sG$
(resp. $d$) is either arbitrary or clear from context.

\example{\label{X.d.pole}
(a) If $d=1$, then a 1-pole is a pole in the sense of
\S\ref{S:pole} (apply 
Example \ref{X:higher.cocycle}).

(b) Let $C:\bY^2\x\gQ\into\Zahl$ be as in 
Example \ref{X:ice.cube.cocycle}, and let $\bq\in\tl\gQ$ be
the configuration in Figure \ref{fig:ice.cube}(D), having a
codimension-three defect.  Let $\beta\in\Zahl[\bY^3]$
be the three-chain consisting
of the the twenty-seven cubes containing the twenty-seven balls in 
Figure \ref{fig:ice.cube}(D),
and let $\zeta:=\partial_3(\beta)$.  Then $C(\zeta,\bq)=6\neq 0$,
so this is a projective $(\Zahl,2)$-pole.}

The next result is analogous to  Theorem \ref{cohom.defect}(a):

\Proposition{\label{projective.d.pole.essential}}
{
 If $\ba\in\tlgA$ has a projective  $d$-pole, 
then $\ba$ has an essential defect.
}
\bthmprf (by contradiction)
  Suppose $\ba$ had a removable defect; we will show that $\ba$ has no
projective $d$-poles.  Let $\ba'\in\gA$ and
suppose $\ba'$ agrees with $\ba$ on $\unflawed(\ba)$ for some $r\in\Natur$.

\Claim{For all $R\geq r$, $\ba$ has no $d$-poles of range $R$.}
\bclaimprf
Let $(\sG,+)$ be a group and
suppose $C\in\Zeq^d(\gA,\sG)$ has range $R$.
Let  $\zeta\in\sZ_d[\unflawed[R](\ba)]$ be any $d$-cycle.
Then $\zeta=\partial_{d+1} \beta$ for some $\beta \in 
\Zahl[\bY^{d+1}]$ (because $\bY=\RD$ has trivial homology in all dimensions).
Thus $C(\zeta,\ba) \ \eeequals{(*)}\ C(\zeta,\ba')=C(\partial_{d+1}\beta,\ba')=\del_d C(\ba')[\zeta]
\ \eeequals{(\dagger)} \ 0$.  Here, $(*)$ is because $\ba$ and $\ba'$
agree on $\unflawed[R](\ba)$, while $(\dagger)$ is 
 because $\del_d C(\ba')=0$ because $\ba'\in\gA$
and $C$ is a cocycle on $\gA$.
\eclaimprf
  Claim 1 holds for all $\sG$ and all large enough $R$, so $\ba$ has no
projective $d$-pole.
\ethmprf

\Remarks{(a) In \S\ref{S:pole}, it was necessary to
define 1-poles in terms of {\em trails} in
$\pi_1[\unflawed(\ba)]$ (rather than {\em one-chains} in
$\Zahl[\bY^1]$), because \S\ref{S:pole} dealt with potentially {\em
nonabelian} cocycles [e.g. Example \ref{X:cocycle}(f)], where the
(ordered) product in eqn.(\ref{trail.cocycle}) is well-defined, but
where the (unordered) sum in eqn.(\ref{chain.cocycle}) is not.

(b)  If $\ba$ has a  projective $C$-pole, then for every large
$r\in\Natur$, there is some $\zeta_r\in\sZ_d[\unflawed(\ba)]$ such
that $C(\zeta_r,\ba)\neq 0$. We can further assume that for every $R\geq r$,
 the cycle $\zeta_R$ is homologous to $\zeta_r$ in $\unflawed(\ba)$.  We could
then define `projective $C$-poles' by treating $C$ as a function on the inverse
limit group $\sH_d[\unflawed[\oo](\ba)]$ in the obvious way, but we
will restrain ourselves.
}

\subsection{Invariant Cohomology for  Subshifts of Finite Type:
\label{S:invariant.cohom}}

The goal of this section is to determine when $d$-poles are persistent
under a cellular automaton.  To do this, we will introduce another
cohomology group $\Hinv^d(\gA,\sG)$, which acts as a bridge from the
group $\sH^d(\gA,\sG)$ of \S\ref{S:projective.homotopy}  to the group
$\Heq^d(\gA,\sG)$ of \S\ref{S:equiv.cohom}.  We will then prove:

\Theorem{\label{persistent.dpole}}
{
 Let $\gA\subset\AZD$ be an SFT.
 Let $\Phi\colon\!\AZD\!\rightarrow\!\AZD$ be a CA with $\Phi(\gA)\subseteq\gA$.
\bthmlist
  \item For all $d\in\CC{0...D}$, $\Phi$ induces endomorphisms
$\Hinv^d\Phi:\Hinv^d(\gA,\sG)\into\Hinv^d(\gA,\sG)$.

  \item Suppose $\Hinv^d\Phi$ is an epimorphism.
\bitem
  \item[{\rm[i]}] If $\ba\in\tlgA$ has a projective
$(\sG,d)$-pole, and $\bb:=\Phi(\ba)$, and 
$\Ext{\sH_{d-1}[\unflawed(\bb)],\sG}=0$ for all large $r\in\Natur$,
then $\bb$ also has a  projective $(\sG,d)$-pole.

  \item[{\rm[ii]}] In particular, any projective $1$-pole or $D$-pole is 
$\Phi$-persistent. 

  \item[{\rm[iii]}] If $\sG$ is the additive group of a field
{\rm (e.g. $\sG=\Zahlmod{p}$ for $p$ prime)}, then all projective
$\sG$-poles are $\Phi$-persistent.
\eitem
\ethmlist
}
(Here, $\sH_{d-1}[\unflawed(\bb)]$ is homology with coefficients in $\Zahl$.
See Appendix \S\ref{S:homo.alg} for the definition of Ext, and
for other homological background for what follows.)
 Suppose $\gA$ has radius $R>0$, and fix $r\geq R$.  Let
$\bX_r=(\bX_r^0,\ldots,\bX_r^D)$ be the cellular complex induced by
the radius-$r$ Wang representation of $\gA$, as in
\S\ref{S:projective.homotopy}.  For any
$\fz\in\ZD$, we define a self-homeomorphism $\Xi^\fz:\bX_r\into\bX_r$
such that, for any $\fx\in\ZD$ and $\ba\in\gA_{(r)}$, if
$\fy=\fx+\fz$, then $\Xi^\fz\restr{}
:\CornerTile{\fx}{\ba}\into\CornerTile{\fy}{\ba}$ is a homeomorphism
(see also `property (e)' in \S\ref{S:wang.homo}).
Also, $\Xi^\fz$ is a cellular map; i.e.  for all $d\in\CC{0...D}$, if
$\Xi_d^\fz:=\Xi^\fx\restr{\bX_r^d}$, then
$\Xi_d^\fz:\bX_r^d\into\bX_r^d$.  Thus, $\Xi^\fz$ induces a
(contravariant) automorphism $\Xi_\fz^d:\sC^d(\bX_r,\sG)
\into\sC^d(\bX_r,\sG)$ defined by $\Xi_\fz^d(C) \ = \
C\circ\Xi_d^\fz$.   A cochain
$C\in\sC^d(\bX_r,\sG)$ is {\dfn $\Xi$-invariant} if $C =
C\circ\Xi_d^\fz$ for all $\fz\in\ZD$.  Let $\Cinv^d=\Cinv^d(\bX_r,\sG)$
be the subgroup of $\Xi$-invariant cochains.  Then
$\del_{d}[\Cinv^d] \subseteq \Cinv^{d+1}$, because
 $\del_{d}\circ \Xi_\fz^d \ = \
\Xi_\fz^{d+1}\circ\del_d$ for all $d\in\CC{0...D}$.   Thus
we get a chain complex $\bCinv = (\Cinv^0 \stackrel{\del_0}{\into}
\Cinv^1 \stackrel{\del_1}{\into}
\Cinv^2 \stackrel{\del_2}{\into}\cdots)$, and for each
$d\in\CC{0...D}$, we define 
 $\Hinv^d(\bX_r,\sG) := \sH^d(\bCinv)$
(see Appendix \S\ref{S:homo.alg}).

If $\zeta_r:\bX_{r+1}\into\bX_r$ is the surjection defined
in \S\ref{S:projective.homotopy}, then for each $d\in\CC{0...D}$, 
$\zeta_r$ induces a contravariant homomorphism $\zeta_r^d:
\Cinv^d(\bX_r,\sG) \into \Cinv^d(\bX_{r+1},\sG)$ defined by
$\zeta_r^d(C):= C\circ\zeta_r$.
(This works because $\zeta_r\circ\Xi_{d}^\fz=\Xi_d^\fz\circ\zeta_{r}$  for all $\fz\in\ZD$, so if $C$ is $\Xi$-invariant then
so is $C\circ\zeta_r$).  These homomorphisms together comprise
a chain map $\bzet_r:\bCinv(\bX_{r})\into\bCinv(\bX_{r+1})$,
which yields homomorphisms
$\sH^d\zeta_r:\Hinv^d(\bX_r,\sG)\into \Hinv^d(\bX_{r+1},\sG)$
for all $d\in\CC{0...D}$.
We then define the $d$th {\dfn invariant cohomology group} to be the
direct limit:
\begin{eqnarray}
\label{inv.cohom.defn}
\lefteqn{  \Hinv^d(\gA,\sG) \  :=  }\qquad\quad\\
\nonumber
&&   \dirlim
\lb( \Hinv^d(\bX_R,\sG) \xrightarrow{\sH^d\zeta_R}
\Hinv^d(\bX_{R+1},\sG) \xrightarrow{\sH^d\zeta_{R+1}}
\Hinv^d(\bX_{R+2},\sG) \xrightarrow{\sH^d\zeta_{R+2}}
\cdots\rb)
\end{eqnarray}
\ignore{There is a natural embedding $\Hinv^d(\gA,\sG)\hookrightarrow \sH^d(\gA,\sG)$.
To see this, compare definitions (\ref{tile.cohom.defn}) and 
(\ref{inv.cohom.defn}),
and recall that  $\Hinv^d(\bX_r,\sG)\subseteq \sH^d(\bX_r,\sG)$ for
all $r\geq R$.}

 Fix $d\in\CC{0...D}$, and
let $\Ceq^d(\gA,\sG)$ be as in \S\ref{S:equiv.cohom}. For any $r\geq R$, 
let $\rCeq^d:=\rCeq^d(\gA,\sG)$ be the set of locally determined 
equivariant $d$-dimensional cochains of radius $r$ on $\gA$. 
Then $\rCeq[1]^d\subseteq\rCeq[2]^d\subseteq\rCeq[3]^d\subseteq\cdots\subseteq
 \Ceq^d(\gA,\sG)$.  Let $\rCeq[\oo]^d := \Union_{r=1}^\oo \rCeq^d$.

Now fix $r\geq R$, and observe that
$\del_d(\rCeq^d) \subseteq \rCeq^{d+1}$ for each $d\in\Natur$.  Thus, we
get a chain complex 
$\rbCeq := (\rCeq[r]^0 \stackrel{\del_0}{\into}\rCeq[r]^1 \stackrel{\del_1}{\into}\cdots)$.  For any $d\in\CC{0...D}$, let
$\rHeq^d(\gA,\sG) := \sH^d(\rbCeq)$ be the $d$th 
{\dfn radius-$r$ equivariant cohomology group}.
Likewise, $\del_d(\rCeq[\oo]^d) \subseteq \rCeq[\oo]^{d+1}$,
yielding a chain complex 
$\rbCeq[\oo] := (\rCeq[\oo]^0 \stackrel{\del_0}{\into}\rCeq[\oo]^1 \stackrel{\del_1}{\into}\cdots)$.  For any $d\in\CC{0...D}$, we define
$\rHeq[\oo]^d(\gA,\sG) := \sH^d(\rbCeq[\oo])$.  This section's other main result is:

\Theorem{\label{three.cohomology.theorem}}
{
Let $\gA\subset\AZD$ be an SFT, and let $(\sG,+)$ be an abelian group.
\bthmlist
\item There is a natural isomorphism $\Hinv^d(\gA,\sG)\cong\rHeq[\oo]^d(\gA,\sG)$.

\item If $\sG$ is discrete, then $\rHeq[\oo]^d(\gA,\sG) = \Heq^d(\gA,\sG)$,
so $\Hinv^d(\gA,\sG)\cong\Heq^d(\gA,\sG)$.

In particular, $\Hinv^1(\gA,\sG)\cong \Hdyn^1(\gA,\sG)$.
\ethmlist
}
This, in turn, will follow from:

\Proposition{\label{Heq.vs.Hinv}}
{
Let $\gA\subset\AZD$ be an SFT of radius $R>0$.
Let $(\sG,+)$ be an abelian group.
For any $r\geq R$ and  $d\in\CC{0...D}$,
there is a canonical isomorphism
$\Hinv^d(\bX_r,\sG)\cong\rHeq^d(\gA,\sG)$.
}

To prove Proposition \ref{Heq.vs.Hinv}, let $\bY=(\bY_0,\ldots,\bY_D)$
be the canonical cell complex for $\RD$ (see
\S\ref{S:canonical.cell.cplx}), and let $\sC^d(\bY,\sG)$ be its cubic
cohomology group (see \S\ref{S:cubic.cohom}).  If
$C\in\sC^d(\bX_r,\sG)$ is any cochain on the tile complex $\bX_r$,
then $C$ induces a function $\tlC:\gA\into\sC^d(\bY,\sG)$ defined by
$\tlC(\ba)=C\circ\varsigma^r_\ba$, for all $\ba\in\gA$ (where
$\varsigma^r_\ba:\bY\into\bX_r$ is the continuous section of $\Pi_r$ induced by
$\ba$, as in \S\ref{S:projective.homotopy}).  Proposition
\ref{Heq.vs.Hinv} then follows immediately from the next result:

\Lemma{\label{equiv.vs.inv}}
{
\bthmlist
\item  $C$ is a $\Xi$-invariant cochain on $\bX_r$ iff $\tlC$ 
is an equivariant cochain of radius $r$ on $\gA$.  

\item This defines an isomorphism
$\psi_r^d:\Cinv^d(\bX_r,\sG)\ni C\mapsto \tlC \in \rCeq^d(\gA,\sG)$.  

\item $C\in\Zinv^d(\bX_r,\sG)$ iff $\tlC\in \rZeq^d(\gA,\sG)$,
and $C\in\Binv^d(\bX_r,\sG)$ iff $\tlC\in \rBeq^d(\gA,\sG)$.
\ethmlist
}
\bthmprf
{\bf(a)}
 \  If $\ba\in\gA$ and $\fz\in\ZD$,  then for any $\zeta\in\Zahl[\bY^d]$,
\beqn
\label{Heq.vs.Hinv.e1}
\begin{array}{rcl}
\tlC(\Upsilon_d^{\fz}(\zeta),\ba)
&=&
C\circ\varsigma^r_\ba\circ\Upsilon_d^{\fz}(\zeta)\\
\And
\tlC(\zeta,\shift{\fz}(\ba))
& = &  C\circ\varsigma^r_{\shift{\fz}(\ba)}(\zeta) 
\ \ \ \eeequals{(\diamond)} \  \ \
C\circ\Xi_d^{-\fz}\circ\varsigma^r_\ba\circ\Upsilon_d^{\fz}(\zeta)
\end{array}
\eeqn
Here, $(\diamond)$ is because 
$\varsigma^r_{\shift{\fz}(\ba)} =  \Xi^{-\fz}\circ\varsigma^r_\ba\circ\Upsilon^{\fz}$.
Thus,
\beq
\statement{$\tlC\in\Ceq^d(\gA,\sG)$}
&\iiiff{(*)}&
\statement{$\forall \ \zeta\in\Zahl[\bY^d]$, $\ba\in\gA$, and $\fz\in\ZD$, \ 
$\tlC(\zeta,\shift{\fz}(\ba))
\ = \ \tlC(\Upsilon_d^{\fz}(\zeta),\ba)$}
\\
&\iiiff{(\dagger)}&
\statement{$\forall\  \fz\in\ZD, \ 
C\circ\Xi_d^{-\fz} \ = \ C$}
\iiiff{(\ddagger)}
\statement{$C\in\Cinv^d(\bX_r,\sG)$}.
\eeq
$(*)$ is defining eqn.(\ref{cochain.eqn}) from \S\ref{S:equiv.cohom}.
$(\ddagger)$ is the definition of $\Cinv^d(\bX_r,\sG)$.
Then `$\seilpmiii{(\dagger)}$' is by
eqn.(\ref{Heq.vs.Hinv.e1}),
while `$\iiimplies{(\dagger)}$' is by
eqn.(\ref{Heq.vs.Hinv.e1}) and the fact that
$\Union_{\ba\in\gA} \varsigma^r_\ba(\bY^d) = \bX_r^d$.

$\tlC$ has radius $r$ because $C$ is defined on the $d$-cells
of $\bX^d_r$, each of which is `labelled' by a block
in $\gA_{\partial_y \dB(r)}$ for some $y\in\bY_0^d$
(see \S\ref{S:equiv.cohom}).   Hence, for any $y\in\bY_0^d$,
$\tlC(\ba,y) = C\circ\varsigma^r_\ba(y)$ depends only
on $\ba_{\partial_y \dB(r)}$.

{\bf(b)} {\em $\psi^d$ is injective:}  If $\tlC_1=\tlC_2$,
then $C_1\circ\varsigma^r_\ba=C_2\circ\varsigma^r_\ba$ for all
$\ba\in\gA$.  But $\Union_{\ba\in\gA} \varsigma^r_\ba(\bY^d) = \bX_r^d$,
so this means $C_1(x)=C_2(x)$ for all $x\in\bX_r^d$, which means
$C_1=C_2$.

{\em $\psi^d$ is surjective:}  Let $C'\in \rCeq^d(\gA,\sG)$ be
a radius-$r$ equivariant cochain;  we seek some $C\in \Cinv^d(\bX_r,\sG)$
so that $\tlC=C'$.  

\Claim{Let $x\in\bX_r^d$,  and let $y:=\Pi_r^d(x)\in\bY^d$.

{\rm(a)} \  There exists $\ba\in\gA$  such that $x=\varsigma^r_\ba(y)$.

{\rm(b)} \ If $\ba'\in\gA$ is another element with $\varsigma^r_{\ba'}(y) =x$, then 
 $C'(y,\ba) = \  C'(y,\ba')$.
}
\bclaimprf 
(a) \ Suppose  $x$ is a $d$-dimensional face of $\CornerTile{\fz}{\bb}$
for some $\bb\in\gA_{(r)}$ and  $\fz\in\ZD$.  Let
 $\dF:=\partial_y \, \dB(\fz,r)$ (defined in \S\ref{S:equiv.cohom}),
and find $\ba\in\gA$ such that $\ba_{\dF}=\bb_{\dF}$.
Then $\varsigma^r_\ba(y)=x$.

(b) \ $C'(y,\ba) \ \eeequals{(*)} \ c_y(\ba_\dF) \ \eeequals{(\dagger)} \ 
 c_y(\ba'_\dF) \ \eeequals{(*)} \  C'(y,\ba')$.
$(*)$ is because $C'$ is locally determined with radius $r$.
$(\dagger)$ is because $\ba_{\dF}=\bb_{\dF}=\ba'_{\dF}$, 
because $\varsigma^r_\ba(y) =x = \varsigma^r_{\ba'}(y)$.
\eclaimprf
Define $C\in\sC^d(\bX_r,\sG)$ as follows:
for any $x\in\bX^d_r$, \  $C(x) := C'(y,\ba)$, where $\ba$ and $y$
are as in Claim 1(a).  Then $C(x)$ is well-defined independent
of the choice of $\ba$, by Claim 1(b).

\Claim{$C\in\Cinv^r(\gA,\sG)$, and \  $\tlC=C'$.}
\bclaimprf
 Let $y\in\bY^d$ and let $\ba\in\gA$.
If $x=\varsigma^r_\ba(y)$, then $x$, $y$, and $\ba$ are related
as in Claim 1, so $\tlC(y,\ba) \ = \ 
C\circ\varsigma^r_\ba(y) \ = \ C(x) \ := \ C'(y,\ba)$, as desired.
Thus, $\tlC = C'$.  But then (a) implies that $C\in\Cinv^r(\gA,\sG)$
(because $C'$ is equivariant).
\eclaimprf

{\bf(c)} This follows from the fact that $\widetilde{\del_d C} = \del_d \tlC$
for any $C\in \sC^d(\bX^r,\sG)$.  This, in turn, is because
$\partial_{d+1} \circ \varsigma^r_\ba(\zeta) \ = \ 
\varsigma^r_\ba \circ \partial_{d+1}(\zeta)$, for any $\ba\in\gA$ and $\zeta\in
\Zahl[\bY^{d+1}]$.
\ethmprf

\bthmprf[Proof of Theorem \ref{three.cohomology.theorem}:]
 {\bf(a)} \
For each $d\in\CC{0...D}$,
recall that $\rCeq[1]^d\subseteq\rCeq[2]^d\subseteq\rCeq[3]^d\subseteq\cdots$.
The inclusion maps 
$\rCeq[1]^d\stackrel{\iota_1 }{\hookrightarrow}\rCeq[2]^d\stackrel{\iota_2 }{\hookrightarrow}\rCeq[3]^d\stackrel{\iota_3 }{\hookrightarrow}\cdots$ define a sequence of chain maps
$(\rbCeq[1] \xrightarrow{\biota_1 } \rbCeq[2] \xrightarrow{\biota_2 } \rbCeq[3] \xrightarrow{\biota_3 }\cdots)$,
which yields a sequence of cohomology homomorphisms
$(\rHeq[1]^d \xrightarrow{\sH^d \iota_{1}} \rHeq[2]^d \xrightarrow{\sH^d \iota_{2}} \rHeq[3]^d \xrightarrow{\sH^d \iota_{3}}\cdots)$.
Also, $\rbCeq[\oo]=\D\lim_{\goto} \ (\rbCeq[1] \xrightarrow{\biota_1 } \rbCeq[2] \xrightarrow{\biota_2 } \cdots)$,
because for each $d\in\CC{0...D}$,
$\rCeq[\oo]^d \ = \ \D \Union_{r=1}^\oo \rCeq^d \ = \
\lim_{\goto} \ (\rCeq[1]^d\stackrel{\iota_1 }{\hookrightarrow}\rCeq[2]^d\stackrel{\iota_2 }{\hookrightarrow}\cdots)$.  Thus, 
Lemma \ref{cohomology.limit.lemma} (in Appendix \S\ref{S:homo.alg})
says 
\beqn
\label{infinite.radius.equivariant.cohomology}
\mbox{For all $d\in\CC{0...D}$,} \quad
\rHeq[\oo]^d(\gA,\sG)
\ \ = \ \ \D\lim_{\goto} \  (\rHeq[1]^d \xrightarrow{\sH^d \iota_{1}} \rHeq[2]^d \xrightarrow{\sH^d \iota_{2}} \rHeq[3]^d \xrightarrow{\sH^d \iota_{3}}\cdots).
\eeqn
  Suppose $\gA$ has radius $R>0$.  Then
Proposition \ref{Heq.vs.Hinv},
eqn.(\ref{inv.cohom.defn}), and eqn.(\ref{infinite.radius.equivariant.cohomology}) together yield a commuting ladder with isomorphism rungs,
which yields an isomorphism of direct limits, as shown:
\[
\Array{
\Hinv^d(\bX_R,\sG) &\xrightarrow{\sH^d\zeta_R}&
\Hinv^d(\bX_{R+1},\sG) &\xrightarrow{\sH^d\zeta_{R+1}}&
\Hinv^d(\bX_{R+2},\sG) &\xrightarrow{\sH^d\zeta_{R+2}}&
\cdots& \Hinv^d(\gA,\sG)\\
\psi_R^d\Longdownarrow & & \psi_{R+1}^d\Longdownarrow & & \psi_{R+2}^d\Longdownarrow & &  &\Longdownarrow\\
\rHeq[R]^d(\gA,\sG) &\xrightarrow{\sH^d \iota_{R}}&
\rHeq[R+1]^d(\gA,\sG) &\xrightarrow{\sH^d \iota_{R+1}}&
\rHeq[R+2]^d(\gA,\sG) &\xrightarrow{\sH^d \iota_{R+2}}&
\cdots& \rHeq[\oo]^d(\gA,\sG).
}
\]
{\bf(b)} If $\sG$ is discrete, then every continuous $\sG$-valued cocycle
is locally determined, hence $\rCeq[\oo]^d(\gA,\sG)=\Ceq^d(\gA,\sG)$
for every $d\in\CC{0...D}$.  This yields an equality
$\rbCeq[\oo](\gA,\sG)=\bCeq(\gA,\sG)$ of chain complexes, so that
$\rHeq[\oo]^d(\gA,\sG)=\Heq^d(\gA,\sG)$ for every $d\in\CC{0...D}$.
This, together with part (a),
implies that $\Hinv^d(\gA,\sG)\cong\Heq^d(\gA,\sG)$.  Then
Theorem \ref{two.cocycles.are.one} implies that
$\Hinv^1(\gA,\sG)\cong \Hdyn^1(\gA,\sG)$.
\ethmprf

We now turn to  Theorem \ref{persistent.dpole}.
If $\ba\in\tlgA$, then for all $r\in\Natur$, 
let $\varsigma^r_\ba:\unflawed[r](\ba)\into\bX_r$ be the continuous section
of $\Pi_r$ from 
\S\ref{S:projective.homotopy}, which induces a  homomorphism
$\rHinv^d\ba\colon \Hinv^d(\bX_r,\sG)\into  \sH^d[\unflawed[r](\ba),\sG]$,
defined by $\rHinv^d\ba(\undC):= \underline{C\circ\varsigma^r_\ba}$.
These homomorphisms converge to a direct limit homomorphism
$\Hinv^d\ba\colon \Hinv^d(\gA,\sG)\into  \sH^d[\unflawed[\oo](\ba),\sG]$,
by an argument analogous to Theorem \ref{homotopy.comm.square}(a).

\Proposition{\label{d.pole.iff.equiv.cohom.defect}}
{
Let $\gA\subset\AZD$ be an SFT. Let $(\sG,+)$ be abelian.
Let $\ba\in\tlgA$.
\bthmlist
  \item For any $r\in\Natur$, \ 
 $(\ba$ has a $d$-pole of range $r) \ \implies \ (\rHinv^d\ba$ is nontrivial$)$.

  If $\Ext{\sH_{d-1}[\unflawed(\ba)],\sG}=0$, then `$\seilpmi$' is also true.

  \item $(\ba$ has a projective $d$-pole$) \ \implies \ 
(\Hinv^d\ba$ is nontrivial$)$.

  If $\Ext{\sH_{d-1}[\unflawed(\ba)],\sG}=0$ for all large $r\in\Natur$, then `$\seilpmi$' is also true.
\ethmlist
}
\bthmprf {\bf(a)} `$\implies$'
Let $\unflawed := \unflawed(\ba)$.
If $\ba$  has a $d$-pole of range $r$, then there exists 
$C'\in\rZeq^d(\gA,\sG)$ 
and a $d$-cycle $\zeta\in\Zahl[\bY^d\intsct\unflawed]$ such
that $C'(\zeta,\ba)\neq0$.   
 Lemma \ref{equiv.vs.inv}(b,c) yields 
$C\in \Zinv^d(\bX_r,\sG)$ with $\tlC=C'$.  Then
 $C\circ\varsigma^r_\ba (\zeta) =\tlC(\zeta,\ba) = C'(\zeta,\ba) \neq 0$.
Hence, the cocycle $C\circ\varsigma^r_\ba$ cannot be nullhomologous on
$\unflawed$.  
But $\underline{C\circ\varsigma^r_\ba} = \rHinv^d\ba(\undC)$.
Hence $\rHinv^d\ba$ is nontrivial.

{\bf(a)} \ `$\seilpmi$'  
Suppose $\rHinv^d\ba$ is nontrivial; hence there exists
$\undC\in \Hinv^d(\bX_r,\sG)$ such that $\undC':=\rHinv^d\ba(\undC)$
is a nontrivial cohomology class in $\sH^d(\unflawed,\sG)$.

\Claim{
There is some $d$-cycle $\zeta\in\sZ_d(\unflawed,\Zahl)$
such that $C'(\zeta)\neq0$.}
\bclaimprf
For any  $C'\in\sZ^d(\unflawed,\sG)$ and $\zeta\in\sZ_d(\unflawed,\Zahl)$,
the value of $C'(\zeta)$ depends only on the cohomology class of $C'$
and the homology class of $\zeta$.  Thus, if $\sH_d:=\sH_d(\unflawed,\Zahl)$,
then $\undC'$
defines a function $\sH_d\into\sG$.
This construction yields a
homomorphism $\sH^d(\unflawed,\sG)\into\Hom{\sH_d, \sG}$.
The Universal Coefficient Theorem \cite[Theorem 3.2]{Hatcher}
says that the kernel of this homomorphism is
$\Ext{\sH_{d-1},\sG}$.  Hence, if
$\Ext{\sH_{d-1},\sG}=0$, then any nontrivial
cohomology class $\undC'\in\sH^d(\unflawed,\sG)$ defines
a nontrivial element of $\Hom{\sH_d, \sG}$,
which means $C'(\zeta)\neq0$ for some
$\zeta\in\sZ_d(\unflawed,\Zahl)$.
\eclaimprf
If $\undC':=\rHinv^d\ba(\undC)$, then $C'=C\circ\varsigma^r_\ba$, so
Claim 1 means that  
$\tlC(\zeta,\ba)\neq 0$, where $\tlC\in\rHeq^d(\gA,\sG)$ is defined as
prior to Lemma \ref{equiv.vs.inv}.  Thus, $\ba$ has a $d$-pole of range $r$.

{\bf(b)} \ ($\ba$ has a projective $d$-pole) $\iff$ 
($\ba$  has a $d$-pole of range $r$ for all large $r\in\Natur$)
$\iiimplies{(*)}$ ($\rHinv^d\ba$ is nontrivial for all large $r\in\Natur$)
$\iff$ ($\Hinv^d\ba$ is nontrivial).  Here ``$\iiimplies{(*)}$'' is by part (a),
and becomes a ``$\iiiff{(*)}$'' if  $\Ext{\sH_{d-1}[\unflawed(\ba)],\sG}=0$
 for all large  $r\in\Natur$.
\ethmprf

\bthmprf[Proof of Theorem \ref{persistent.dpole}:]  {\bf(a)}
Suppose $\Phi$ has radius $q$.
Fix $r\in\Natur$, and recall that $\Phi$ is $\shift{}$-commuting.
Hence, for all $d\in\CC{0...D}$, in the proof of Proposition
\ref{CA.induced.homotopy.map}(a), the induced cellular map
$\Phi_*:\bX^d_{r+q}\into\bX^d_r$ is $\Xi_d$-commuting, so the induced
homomorphism $\sC_r^d\Phi:\sC^d(\bX_r,\sG) \into \sC^d(\bX_{r+q},\sG)$ is
$\Xi^d$-commuting.  Thus, $\sC_r^d\Phi$ restricts to a map
$\sC_r^d\Phi:\Cinv^d(\bX_r,\sG) \into \Cinv^d(\bX_{r+q},\sG)$.  This
yields a chain map $\sC_r\Phi:\bCinv(\bX_r,\sG) \into
\bCinv(\bX_{r+q},\sG)$, which yields cohomology homomorphisms
$\sH_r^d\Phi:\Hinv^d(\bX_r,\sG) \into \Hinv^d(\bX_{r+q},\sG)$ for all
$d\in\CC{0...D}$.  Thus, taking the direct limit (as in Proposition
\ref{CA.induced.homotopy.map}) yields a homomorphism
$\Hinv^d\Phi:\Hinv^d(\gA,\sG)\into \Hinv^d(\gA,\sG)$.

{\bf(b)}[i] If $\ba\in\tlgA$ and $\bb:=\Phi(\ba)$, then part (a) and
an argument analogous to Theorem \ref{homotopy.comm.square}(c) yield a
commuting square
\[
 \begin{array}{ccc}
\sH^k[\unflawed[\oo](\ba),\sG] & \xleftarrow{\ \sH^k\iota \ } & \sH^k[\unflawed[\oo](\bb),\sG]\\
{\scriptstyle \Hinv^k\ba} \longuparrow & & \longuparrow {\scriptstyle \Hinv^k\bb} \\
\Hinv^k(\gA,\sG) & \xleftarrow{\Hinv^k\Phi} & \Hinv^k(\gA,\sG)\\
\end{array}
\]
Thus, if $\Hinv^k\Phi$ is surjective and 
$\Hinv^k\ba$ is nontrivial, then $\Hinv^k\bb$ must also be nontrivial.
Then Proposition \ref{d.pole.iff.equiv.cohom.defect}(b) says:
if $\ba$ has a projective $d$-pole, 
and $\Ext{\sH_{d-1}[\unflawed(\bb)],\sG}=0$ for all large $r\in\Natur$,
then $\bb$ has a projective $d$-pole.

{\bf(b)}[ii]  follows from (b)[i] because 
$\sH_{0}[\unflawed(\bb)]=\Zahl^K$, where $K$ is the number of connected components
of $\unflawed(\bb)$ \cite[Proposition 2.7]{Hatcher},
so $\Ext{\sH_{0}[\unflawed(\bb)],\sG}=0$.
 Also, for each $r\in\Natur$,
$\bG_r(\bb)$ is homotopic to an (orientable) $D$-dimensional submanifold
of $\RD$, so $\sH_{D-1}[\unflawed(\bb)]$ is torsion-free \cite[Corollary 3.28]{Hatcher}, so $\Ext{\sH_{D-1}[\unflawed(\bb)],\sG}=0$.

{\bf(b)}[iii]  follows from (b)[i] because if $\sG$ is the additive group
of a field, then $\Ext{\sH,\sG}=0$ for any group $\sH$.
\ethmprf


\subsection{Appendix on Homological Algebra: \label{S:homo.alg}}
  If $\sH$ is any abelian group,
we can always write $\sH \cong \sF_1/\sF_2$, where $\sF_1$ is a free abelian group
and $\sF_2$ is a subgroup.   Let $\sG$ be another abelian group.
Let $\sH^*:=\Hom{\sH,\sG}$, \  $\sF_1^*:=\Hom{\sF_1,\sG}$, and
$\sF_2^*:=\Hom{\sF_2,\sG}$, and observe that the short
exact sequence $0\into\sF_2\stackrel{i}{\inject}\sF_1\stackrel{q}{\surject}\sH\into0$
induces a sequence $\sF_2^*\stackrel{i^*}{\longleftarrow}\sF_1^*
\stackrel{q^*}{\longleftarrow}\sH$ where 
$i^*(\phi):=\phi\circ i$ and $q^*(\phi):=\phi\circ q$.  Now, $i^*\circ
q^* = 0$ because $q\circ i=0$ (by definition); hence
$\image{q^*}\subseteq\ker(i^*)$.
We thus define
$\Ext{\sH,\sG}:=\ker(i^*)/\image{q^*}$.  This definition is
independent of the choice of `free resolution' $(\sF_2,\sF_1)$.
Furthermore, if $\sH$ is a finitely generated abelian group, so that
$\sH \cong \Zahl^R \dirsum \Zahlmod{n_1}\dirsum\cdots\dirsum
\Zahlmod{n_K}$ (for some $R$ and $n_1,\ldots,n_K$), then
$\Ext{\sH,\sG} \cong \Dirsum_{k=1}^K (\sG/n_k\sG)$.  In particular, if
$\sH\cong \Zahl^R$, or if $(\sG,+)$ is the additive group of a field,
then $\Ext{\sH,\sG}=0$.  See \cite[p.195 of \S3.1]{Hatcher}. 

A {\dfn chain complex} is an infinite sequence of abelian groups
and homomorphisms 
$\sC^0\stackrel{\del^0}{\into} \sC^1\stackrel{\del^1}{\into} 
\sC^2\stackrel{\del^2}{\into} \cdots$ such that $\del^{n+1}\circ\del^n=0$
for all $n\in\Natur$.   We represent this structure as 
$\bC:=\{\sC^n,\del^n\}_{n=0}^\oo$.  
(In fact, for our purposes, only the groups $\sC^0,\ldots,\sC^D$ are
nontrivial;  however it is both conventional and convenient
to develop the theory for infinite chain complexes.)
If $\sZ^n:=\ker(\del^n)$ and $\sB^n:=\image{\del^{n-1}}$, then
$\sB^n\subseteq\sZ^n$.  We define $\sH^n(\bC):=\sZ^n/\sB^n$,
to be the $n$th {\dfn cohomology group} of the chain complex $\bC$
(we formally define $\sB^0:=\{0\}$, so $\sH^0=\sZ^0$).
See \cite[\S2.1]{Hatcher} or \cite[\S IV.2]{Lang}.

 If $\bC_1:=\{\sC_1^n,\del_1^n\}_{n=0}^\oo$ is another chain complex,
then a {\dfn chain map} from $\bC$ to $\bC_1$ is a sequence of
homomorphisms $\bphi:=\{\phi^n:\sC^n\into\sC^n_1\}_{n=0}^\oo$ such that
$\del_1^n\circ\phi^n = \phi^{n+1}\circ\del^n$ for all $n\in\Natur$.
We indicate this by writing ``$\bphi:\bC\into\bC_1$''.  The set of all
chain complexes and chain maps forms a category $\gC$, and
cohomology yields functors $\sH^n$ from  $\gC$ 
to the category $\gA$ of abelian groups. 
To be precise if $\bphi:\bC\into\bC_1$ is a chain map,
then there is a homomomorphism $\sH^n\bphi:\sH^n(\bC)\into\sH^n(\bC_1)$
defined by $\sH^n\bphi(z+\sB^n)=\phi^n(z)+\sB^n_1$ for any $z\in\sZ^n$.
(Recall that elements of $\sH^n(\bC)$ are cosets of $\sB^n$ in $\sZ^n$.
The function $\sH^n\bphi$ is well-defined because $\phi^n(\sZ^n)\subseteq
\sZ^n_1$ and $\phi^n(\sB^n)\subseteq \sB^n_1$).
See \cite[Prop. 2.9]{Hatcher}.

  Let $\bC_1\stackrel{\biota_1}{\into} \bC_2\stackrel{\biota_2}{\into} \bC_3\stackrel{\biota_3}{\into} \cdots$ be an infinite sequence of chain complexes
and chain maps.  The {\dfn direct limit} is the chain complex 
$\bC:=\{\sC^n,\del^n\}_{n=0}^\oo$, where for each $n\in\Natur$,
$\sC^n := \D\lim_{\goto}  \ (\sC^n_1\stackrel{\iota^n_1}{\into} \sC^n_2\stackrel{\iota^n_2}{\into} \sC^n_3\stackrel{\iota^n_3}{\into} \cdots)$,
and where the maps $\del^n:\sC^n\into\sC^{n+1}$ arise from the commuting grid:
{\[
\begin{array}{rcrcrcrr}
\vdots & &\vdots & &\vdots & &&\vdots \\
{\scriptstyle \del^{n-1}_1}\longdownarrow & &  {\scriptstyle \del^{n-1}_2}\longdownarrow & &  {\scriptstyle \del^{n-1}_3}\longdownarrow & &&  {\scriptstyle \del^{n-1} }\longdownarrow \\
\sC^n_1 & \xrightarrow{\iota^n_1} & \sC^n_2 & \xrightarrow{\iota^n_2} & 
\sC^n_3 & \xrightarrow{\iota^n_3} & \cdots & \sC^n \\
{\scriptstyle \del^n_1}\longdownarrow & &  {\scriptstyle \del^n_2}\longdownarrow & &  {\scriptstyle \del^n_3}\longdownarrow & &&  {\scriptstyle \del^n}\longdownarrow \\
\sC^{n+1}_1 & \xrightarrow{\iota^{n+1}_1} & \sC^{n+1}_2 & \xrightarrow{\iota^{n+1}_2} & 
\sC^{n+1}_3 & \xrightarrow{\iota^{n+1}_3} & \cdots& \sC^{n+1} \\
{\scriptstyle \del^{n+1}_1}\longdownarrow & &  {\scriptstyle \del^{n+1}_2}\longdownarrow & &  {\scriptstyle \del^{n+1}_3}\longdownarrow & &&  {\scriptstyle \del^{n+1} }\longdownarrow \\
\vdots & &\vdots & &\vdots & &&\vdots 
\end{array}
\]
}
Fix $n\in\Natur$.  The sequence of chain maps
$(\bC_1\stackrel{\biota_1}{\into} \bC_2\stackrel{\biota_2}{\into}\cdots)$ induces a sequence of cohomology homomorphisms
$(\sH^n\bC_1\xrightarrow{\sH^n\iota_1} \sH^n\bC_2\xrightarrow{\sH^n\iota_2}  \cdots)$.  

\Lemma{\label{cohomology.limit.lemma}}{
For any $n\in\Natur$, \ 
$\sH^n\bC=\D\lim_{\goto} \ (\sH^n\bC_1\xrightarrow{\sH^n\iota_1} \sH^n\bC_2\xrightarrow{\sH^n\iota_2} \sH^n\bC_3\xrightarrow{\sH^n\iota_3} \cdots)$.
}
\bthmprf For every $r\in\Natur$, we have short exact sequences 
$0\into\sB_r^n \inject\sZ_r^n\surject\sH^n(\bC_r)\into0$,
where $\sB_r^n:=\image{\del^{n-1}_r}\subseteq\sZ_r^n:=
\ker(\del^{n}_r)\subseteq\sC_r^n$.
Also, $\iota_r^n(\sZ_r^n)\subseteq
\sZ^n_{r+1}$ and $\iota_r^n(\sB_r^n)\subseteq \sB^n_{r+1}$.
 If $\sB^n:=
\D\lim_{\goto} \ (\sB_1^n\xrightarrow{\iota^n_1} \sB_2^n\xrightarrow{\iota^n_2} \cdots)$
and $\sZ^n:=
\D\lim_{\goto} \ (\sZ_1^n\xrightarrow{\iota^n_1} \sZ_2^n\xrightarrow{\iota^n_2} \cdots)$,
then $\sB^n=\image{\del^{n-1}}\subseteq\sZ^n=\ker(\del^n)\subseteq\sC^n$.
If $\tlsH^n := \D\lim_{\goto} \ (\sH^n\bC_1\xrightarrow{\sH^n\iota_1} \sH^n\bC_2\xrightarrow{\sH^n\iota_2}\cdots)$, then
these short exact sequences converge to a short exact
sequence $0\into\sB^n{\inject}\sZ^n{\surject}\tlsH^n\into 0$.
\ignore{\cite[Ex.22, Ch.III]{Lang}.}  Thus,  $\tlsH^n\cong\sZ^n/\sB^n$.
But $\sZ^n/\sB^n=\sH^n\bC$ by definition.
\ethmprf

\subsection*{Conclusion:}

  We have developed algebraic invariants which help to explain the
emergence, persistence, and interaction of defects in cellular
automata.  However, many questions remain.

1. Example \ref{X:ice.pole}(c) shows that our set of algebraic invariants
is not yet sufficient to detect all essential defects.  Are there
other algebraic invariants?

2.  Proposition \ref{persistent.pole}(b), Theorem \ref{persistent.gap}(b),
Corollary \ref{persistent.homotopy.defect}
and Theorem \ref{persistent.dpole}(b) all say that homotopic or (co)homological
defects are $\Phi$-persistent if the homotopy/(co)homology homomorphism induced by $\Phi$  is injective/surjective.
When is this the case?

3.  With the exception of a few examples of $\pi_1$ computed in in
\cite{GePr}, there are no explicit computations of the
homotopy/(co)homology groups of \S\ref{S:projective.homotopy} and
\S\ref{S:invariant.cohom}, partly because of the difficulty of taking
the required inverse/direct limits.  This limits the applicability
of the theory. Is there an easy way to compute these
groups?

4. The pole/residue theory of \S\ref{S:pole} suggests an appealing analogy
between two-dimensional symbolic dynamics and complex analysis.
Is there a deeper relationship beyond this analogy?

5.  Conway has shown that the Penrose tiling has exactly 61 distinct
$\shift{}$-homocliny classes of essential codimension-two defects,
by using the fact that any Penrose tiling can be cross-hatched by 
`Ammann bars' \cite[\S10.5, p.566]{GrunbaumShephard}.  Can this method
be extended to some two-dimensional subshifts of finite type?

6.  If $\gA\subset\AZD$, and there is a CA $\Phi$ with
$\Phi^n(\AZD)\subseteq\gA\subseteq\Fix{\Phi}$, then $\gA$ admits no
essential defects.  The converse is also true, when $\gA$ is a
one-dimensional sofic shift with a $\shift{}$-fixed point
\cite{Maass}.  Is the converse true in higher dimensions?  

7. Even when $\gA$ admits essential defects, 
\Kurka and Maass \cite{KuMa00,KuMa02,Kur03,Kur05} have described how a
one-dimensional CA can `converge in measure' to $\gA$ through a gradual
process of defect coalescence/annihilation.  Given a subshift
$\gA\subset\AZD$, is it possible to build a CA which converges to
$\gA$ in this sense?

\breath

Finally, we remark that many of the results here should generalize to
random cellular automata (i.e. CA-valued stochastic processes)
which almost-surely preserve a given subshift.  For example, these include 
zero-temperature (anti)ferromagnet models acting on $\Mono$ or
$\Checker$ with random boundary motions \cite{Elo94},
and random tile-rearranging processes on $\Dom$ \cite{CNPr,CKPr}
or $\Ice$ \cite{Elo99,Elo03,Elo05}.

\paragraph*{Acknowledgements:} 
This paper was written during a research leave at Wesleyan University,
and was partially supported by the Van Vleck Fund.  I am grateful to
Ethan Coven, Adam Fieldsteel, Mike and Mieke Keane, and Carol Wood for
their generosity and hospitality. This research was also supported by
NSERC.  Finally, I gratefully acknowledge the influence of George
Peschke, who first introduced me to algebraic topology.

{
\bibliographystyle{alpha}
\bibliography{bibliography}
}

\end{document}